\documentclass{amsart}
\usepackage{mathtools,amssymb,stackengine}
\usepackage[margin=10pt,font=small,labelfont=bf,justification=Centering]{caption} [2023/03/12]
\usepackage{bm}
\usepackage{url}
\usepackage{amssymb,latexsym,xcolor}
\usepackage{amsmath,amsthm}
\usepackage[utf8]{inputenc}
\usepackage[english]{babel}
\usepackage{amsrefs}
\usepackage{tkz-euclide} 
\usepackage[mathscr]{eucal}
\usepackage[all]{xy}
\usepackage{tikz-cd}
\usepackage{tikz,pgf}
\usetikzlibrary{matrix}
\usetikzlibrary{shapes.geometric,calc}
\usepackage{adjustbox}
\usepackage{stackrel}
\usepackage{hyperref}
\usepackage{extarrows}
\usepackage{pgfkeys}
\usepackage{palatino}

\newtheorem{theorem}{Theorem}[section]
\newtheorem{lemma}[theorem]{Lemma}
\newtheorem{proposition}[theorem]{Proposition}

\theoremstyle{definition}
\newtheorem{definition}[theorem]{Definition}
\newtheorem{example}[theorem]{Example}

\theoremstyle{remark}
\newtheorem{remark}[theorem]{Remark}
\numberwithin{equation}{section}

\newcommand\dimv{\mathbf{dim}\,}
\newcommand\cdnv{\mathbf{cdn}\,}
\newcommand\F{\mathbf{F}}

\newcommand\Omegan{\mathbf{\Omega}}
\usepackage{mathtools,amssymb,stackengine}
\newcommand\Rightarrowtail[2][]{\ensurestackMath{\mathrel{%
  \stackengine{1pt}{%
    \stackengine{0pt}{\rightarrowtail}{\scriptstyle#2}{O}{c}{F}{F}{S}%
  }{\scriptstyle#1}{U}{c}{F}{F}{S}%
}}}

\DeclareMathOperator{\rep}{rep}
\DeclareMathOperator{\ind}{Ind}

\DeclareMathOperator{\supp}{supp}
\DeclareMathOperator{\md}{mod}

\DeclareMathOperator{\Hom}{Hom}
\DeclareMathOperator{\Supp}{supp}
\DeclareMathOperator{\csupp}{csupp}
\def\P{\mathcal{P}}

\newcommand{\cale}{\mathcal{E}}
\newcommand{\zg}{\gamma}

\newcommand{\ses}[3]{0\rightarrow #1\rightarrow #2\rightarrow#3\rightarrow 0}

\newcommand\V{\mathbf{V}}
\newcommand\W{\mathbf{W}}
\newcommand\U{\mathbf{U}}

\newcommand\thv{\mathbf{\boldsymbol \theta}}
\newcommand\kpv{\mathbf{\boldsymbol \kappa}}
\newcommand\nv{\mathbf{\boldsymbol\eta}}
\newcommand\wpv{\mathbf{\boldsymbol w}}
\newcommand\colorb[1]{\color{blue}{\bm #1}}
\newcommand\smrud[1]{%
  \begingroup\setlength\arraycolsep{0pt}\renewcommand{\arraystretch}{0.5}\Huge\ensuremath{\begin{matrix}  #1\end{matrix}\arrow[rd]\arrow[ru]\endgroup}}
\newcommand\sm[1]{%
  \begingroup\setlength\arraycolsep{0pt}\renewcommand{\arraystretch}{0.5}\Huge\ensuremath{\begin{matrix}#1\end{matrix}\endgroup}}
\newcommand\smru[1]{%
  \begingroup\setlength\arraycolsep{0pt}\renewcommand{\arraystretch}{0.5}\Huge\ensuremath{\begin{matrix}#1\end{matrix}\arrow[ru]\endgroup}}
\newcommand\smrd[1]{%
  \begingroup\setlength\arraycolsep{0pt}\renewcommand{\arraystretch}{0.5}\Huge\ensuremath{\begin{matrix}#1\end{matrix}\arrow[rd]\endgroup}}
\newcommand\smr[1]{%
  \begingroup\setlength\arraycolsep{0pt}\renewcommand{\arraystretch}{0.4}\Huge\ensuremath{\begin{matrix}#1\end{matrix}\arrow[r]\endgroup}}

\definecolor{myred}{cmyk}{0.5,0,0.5,0.5}
\definecolor{mylightblue}{RGB}{20,90,220}
\definecolor{mypink}{RGB}{200,90,90}

\hypersetup{colorlinks=true,linkcolor=[rgb]{0.5, 0.0, 0},citecolor=[rgb]{0.5, 0.0, 0.5}, filecolor=magenta,urlcolor=blue}

\title[Stability for socle-projective categories ]{Stability for socle-projective categories of type $\mathbb{A}$}

\author{Kostiantyn Iusenko}
\address{Departamento de Matematica, Univ de S\~{a}o Paulo, Rua do Matão, 1010,   S\~{a}o Paulo, SP, CEP: 05508-090 -– Brazil.}
\email{iusenko@ime.usp.br}

\author{Gabriel Bravo Rios}
\address{Proyecto Curricular de Matemáticas, Universidad Distrital Francisco José de Caldas, Bogot\'a 110110.}
\email{gbravor@udistrital.edu.co}

\author{Robinson-Julian Serna}
\address{Escuela de Matem\'aticas y Estad\'istica, Universidad Pedag\'ogica y Tecnol\'ogica de Colombia, Tunja, Boyac\'a 150001.}
\email{robinson.serna@uptc.edu.co}


\begin{document}

\begin{abstract} 

We extend the notion of stability in the non-abelian category of poset representations (introduced by Futorny and Iusenko) to the category of socle-projective representations of a given $r$-peak poset $\P$. When $\P$ is a poset of type $\mathbb{A}$, we demonstrate in two distinct ways that every indecomposable peak $\P$-space is stable. First, this is shown using a bilinear form associated with the poset. Second, we prove it by observing that a stability function derived from a geometric model ensures that all indecomposable objects are stable. Along the way, we provide a new geometric realization of the category of socle-projective representations, inspired by the work of Schiffler and Serna [\textit{J. Pure Appl. Algebra} \textbf{224} (2020), no.~12, 106436, 23 pp.; MR4101480]. Finally, we establish a connection between the geometric perspective and the bilinear form approach.
\end{abstract}

\maketitle


\section*{Introduction}
The notion of stability for an arbitrary abelian category was developed by Rudakov in \cite{Rud}, where the author formalized the stability conditions for representations of quivers (as introduced by King \cite{King} and Schofield \cite{SC91}). Subsequently, Bridgeland defined stability for an arbitrary triangulated category in \cite{Bridgeland}. This concept has proven to be of significant importance in various applications (see, e.g., \cite{MacriShmidt} and references therein). When attempting to transfer the notion of stability from abelian to additive categories, one encounters certain difficulties, particularly with the existence of quotients in the additive case.\\

Recently, Futorny and Iusenko \cite{futorny} introduced the concept of stability in the non-abelian category of poset representations. (For an alternative approach to (semi)stable representations of posets via quiver representations, see \cite{WY13}.) Their primary goal was to develop a general framework to intrinsically define and study the moduli spaces of posets. In particular, they proved that every indecomposable representation of a poset of finite type is stable with respect to some weight, explicitly constructing this weight using the Euler bilinear form associated with the given poset. For further discussion about bilinear forms associated with posets and their applications, see, for instance, \cites{simsonForm, WY13, CI19} and references therein. It also turned out that such representations are closely related (see \cite{futorny}*{Section 6}) to $\chi$-unitary representations of posets (see \cites{KNR06, KR05, BFKSY13, SY12}).\\

A way to generalize the category of representations of posets is through the category $\mathcal{P}\text{-spr}$ of peak $\mathcal{P}$-spaces (see \cite{simson3}). This category is essentially the category $\text{mod}_{\text{sp}}\, \Bbbk\mathcal{P}$ of finitely generated socle-projective modules over the incidence algebra $\Bbbk\mathcal{P}$ associated with $\mathcal{P}$ over the field $\Bbbk$ (see \cite{simson} for a definition of socle-projective modules and their connection with representations of posets). A specific instance is the classical category of representations of a poset $\mathcal{P}$, as it is equivalent to the category $\text{mod}_{\text{sp}}\, \Bbbk\mathcal{P}^\star$, where $\mathcal{P}^\star$ is the extended poset obtained from $\mathcal{P}$ by adding a maximal point $\star$ (in our terminology, $\mathcal{P}^\star$ is a one-peak poset).  \\

In this paper, we extend the notion of stability introduced in \cite{futorny} to the category $\mathcal{P}\text{-spr}$. In particular, when $\mathcal{P}$ is a poset of type $\mathbb{A}$, we prove (Theorem~\ref{Firsttheorem}) that the indecomposable objects in $\mathcal{P}\text{-spr}$ are stable with a weight given by the bilinear form associated with $\mathcal{P}$ (as in \cite{futorny}). Moreover, using the fact that when $\mathcal{P}$ is a poset of type $\mathbb{A}$, the category $\mathcal{P}\text{-spr}$ is embedded in the category $\text{rep}\, Q$ of representations of a certain Dynkin quiver $Q$ of type $\mathbb{A}$, we show (Theorem~\ref{categoricalequivalence}) a categorical equivalence between the category $\text{mod}_{\text{sp}}\, \Bbbk\mathcal{P}$ and a certain subcategory of a category of line segments in a polygon $P(Q)$ defined in \cite{BGMS}.  This construction provides us with a stability function for which all indecomposable objects in $\text{mod}_{\text{sp}}\, \Bbbk\mathcal{P}$ are stable. This categorical equivalence induces an isomorphism of translation quivers; therefore, we obtain a geometric realization, in the spirit of \cite{schifflerserna}, of the category $\text{mod}_{\text{sp}}\, \Bbbk\mathcal{P}$ using the geometric model of BGMS for the category of representations of the quiver $Q$ given in \cite{BGMS}. Finally, we establish a connection (Section~\ref{section4}) between the geometric perspective explored and the bilinear form approach.  \\

The paper is structured as follows: In Section~\ref{section1}, we introduce the notion of stability in the category $\mathcal{P}\text{-spr}$. Section~\ref{section2} addresses the stability arising from a bilinear form associated with the poset $\mathcal{P}$, specifically when $\mathcal{P}$ is of type $\mathbb{A}$. Section~\ref{section3} presents a geometric realization of the category $\text{mod}_{\text{sp}}\, \Bbbk\mathcal{P}$ using the category of line segments introduced in \cite{BGMS}. Finally, in Section~\ref{section4}, we investigate the stability of socle-projective representations of posets of type $\mathbb{A}$ stemming from the geometric approach discussed in Section~\ref{section3}.  
\\

\textbf{Acknowledgements.} The first-named author acknowledges Mark Kleiner for stimulating discussions about stable representations of posets during the initial stages of the manuscript.
K.I. was partially supported by FAPESP grant 2018/23690-6 and by CNPq Universal Grant 405540/2023-0. The third author was supported by Universidad Pedagógica y Tecnológica de Colombia and Minciencias (Conv. 934).

\section{Category of socle-projective representations and stability} \label{section1}

\subsection{Peak posets and their representations.}
We denote by $(\P,\preceq)$ a finite partially ordered set (poset) with respect to the partial order $\preceq$. We  write $x\prec y$ if $x\preceq y$ and $x\neq y$. For the sake of simplicity, we write $\P$ instead of $(\P,\preceq)$. Let $\max\P$ (respectively $\min\P$) be the set of all maximal (respectively minimal) points of $\P$. A poset $\P$ is called an \textit{r-peak poset} if $\left| \max\P\right|=r$.\\

According to \cite{simson3}, the category $\P$-spr of \textit{peak $\P$-spaces} (or \textit{socle-projective representations} of $\P$) over the field $\Bbbk$ is defined as follows. The objects of $\P$-spr are systems $\U=(U_x)_{x\in\P}$ of finite-dimensional $\Bbbk$-vector spaces $U_x$ such that $U_x$ is a $\Bbbk$-subspace of the \textit{ambient space} $U^{\bullet}=\bigoplus_{z\in\max\P}U_{z}$ for each $x\in\P$, $\pi_z(U_x)=0$ for $x\npreceq z\in\max\P$, and $\pi_y(U_x)\subseteq U_y$ for $x\preceq y$ in $\P$, where $\pi_{y}$ is the composition of the natural maps
$$U^{\bullet}\Rightarrowtail{p_y} U^{\bullet}_y\xhookrightarrow[]{i_y} U^{\bullet},$$ and $$U^{\bullet}_y=\bigoplus_{y\preceq z\in\max\P}U_z.$$

By a \textit{morphism} $f:\U\rightarrow \U'$ in $\P$-spr, we mean a system of $\Bbbk$-linear maps $$f=(f_z:U_z\rightarrow U'_{z})_{z\in\max\P}$$ such that for all $x\in\P$
$$\bigoplus_{z\in\max\P}f_z (U_x)\subseteq U'_{x}.$$ Observe that $f$ is an \textit{isomorphism} if for all $x\in\P,$ $$\bigoplus_{z\in\max\P}f_z (U_x)= U'_{x}.$$ One easily checks that the category $\P$-spr is additive with the unique decomposition property. Moreover, $\P$-spr is a Noetherian category closed under taking kernels and extensions, and it has Auslander-Reiten sequences as well as source maps and sink maps \cite{simson3}*{page 78}.

\subsection{Proper morphisms and proper subrepresentations.}
To define (semi)stable objects in the category $\P$-spr, we first develop the notion of proper subobjects in $\P$-spr. A morphism $f:\U\to \U'$ in $\P$-spr is called \textit{proper} if $f(U_x)=U'_{x}\cap f(U^{\bullet})$ for all $x\in \P$. Let $\U$ be a peak $\P$-space and let $K_z$ be a subspace of $U_z$ for all $z\in\max\P$. The subspace $K=\bigoplus_{z\in\max\P} K_z$ of $U^{\bullet}$ is called an \textit{admissible subspace} of $U^{\bullet}$. We show that the system $\U_K=(U^K_x)_{x\in \P}$, where $U^{K}_x=K\cap U_x$ for each $x\in\P$, is a peak $\P$-space. Observe that if $z\in\max\P$, then $U_z^K=K_z$ and $U_K^\bullet=K$. 

\begin{proposition} Given a peak $\P$-space $\U$ and an admissible subspace $K$ of $U^{\bullet}$, the system $\U_K$ is a peak $\P$-space.
\end{proposition}

\begin{proof} It is clear that $U^K_x$ is a subspace of $K$ for each $x\in \P$. Moreover, given $s,t \in \P$ such that $s \preceq t$, we are going to show that $\widetilde\pi_t(U_s^K)\subseteq U^K_t$, where $\widetilde\pi_t$ is the restriction of $\pi_t$ to $K$. Since $\pi_t(U_s)\subseteq U_t$, then $\widetilde\pi_{t}(U_{s}^{K})=\pi_{t}(U_{s}\cap K)\subseteq \pi_{t}(U_{s})\cap K\subseteq U_{t}\cap K=U_{t}^{K}$. Now, supposing that $x \npreceq z$ and $z \in \max \P$, we get $\pi_{z}(U_{x}^{K})=\pi_{z}(U_{x}\cap K)\subseteq \pi_{z}(U_{x})\cap \pi_{z}(K)=0 \cap \pi_{z}(K)=0$. Thus, $\pi_{z}(U_{x}^{K})=0$.
\end{proof} 

\begin{proposition}\label{quotientincats} If $\V=(V_{x})_{x \in \P}$ is a peak $\P$-subspace of $\U=(U_{x})_{x \in \P}$, whose ambient space $K$ is an admissible subspace of $U^{\bullet}$, and the inclusion $i: \V \hookrightarrow \U $ is a proper morphism, then $\V$ is isomorphic to $\U_K$.
\end{proposition}
\begin{proof} As $V^\bullet =K=\bigoplus_{z\in\max\P }K_z$, we get $V_z=K_z$. Using the identity map $1_z: V_z\to K_z$, we define the morphism $\phi=\bigoplus_{z\in\max\P}1_z: \V \to \U_K$. Thus, $\phi(V_{x})=i(V_x) =U_{x} \cap i(K)=U_{x}\cap K=U_{x}^{K}$. Hence, $\phi$ is an isomorphism.
\end{proof}
Generally, given a peak $\P$-space $\U$ and a peak $\P$-subspace $\V$ of $\U$, the quotient $\U/\V$ does not necessarily belong to $\P$-spr. Nevertheless, in the case when $\V=\U_K$ is a proper peak $\P$-subspace, we have $\U/\U_K \in \P\text{-spr}$ as the following proposition claims.

\begin{proposition} Given a peak $\P$-space $\U$ and $K$ an admissible subspace of $U^{\bullet}$, then $\U/\U_{K}$ is a peak $\P$-space. 
\end{proposition}
\begin{proof} We suppose that $\U=(U_x)_{x\in\P}$ is a peak $\P$-space and that $K=\bigoplus_{z\in\max\P}K_z$ with $K_z \leq U_z$ for all $z\in\max\P$. We are going to show that the system $$\check\U=\U/\U_K=(\check{U}_x)_{x\in\P},$$ where $\check{U}_x=U_x/U^K_x$ for all $x\in\P$, is a peak $\P$-space. Note that the ambient space is given by $$\check\U^\bullet=\bigoplus_{z\in\max\P}U_z/K_z.$$ First, it is enough to check that for all $x\in\P\setminus \max\P$, $\check{U}_x$ is a subspace of $\check\U^\bullet$. To do that, we suppose that $\max\P=\{z_1,\dots,z_k\}$ and that $\rho_i:U^{\bullet}\to U_{z_i}$ is the natural projection. We define a function $\alpha_i:\check{U}_x \to U_{z_i}/K_{z_i}$ such that $\alpha_i(u+U^K_x)=\rho_i(u)+K_{z_i}$ for all $i=1,\dots, k$. Note that $\alpha_i$ is well-defined; if $u_1+U_x^K=u_2+U_x^K$, then $u_1-u_2\in U_x\cap K$, i.e., $u_1-u_2\in K$, then $\rho_i(u_1-u_2)=\rho_i(u_1)-\rho(u_2)\in K_{z_i}$. Thus, $\alpha_i(u_1+U_x^K)=\rho_i(u_1)+K_{z_i}=\rho(u_2)+ K_{z_i}=\alpha_i(u_2+U_x^K)$. Moreover, because $\rho_i$ is $\Bbbk$-linear, we have that $\alpha_i$ is $\Bbbk$-linear. Thus, the function $\alpha: \check{U}_x\to \check\U^\bullet$ given by $\alpha(u+ U_x\cap K)=\sum_{i=1}^k\alpha_i(u+ U_x\cap K)e_i$, where $e_i$ is the vector of $\Bbbk^k$ with one in the $i$-th component and all others zero, is $\Bbbk$-linear. Moreover, note that $\alpha$ is injective; if $\alpha(u_1+U_x\cap K)=\alpha(u_1+U_x\cap K)$, then $\sum_{i=1}^k\alpha_i(u_1+ U_x\cap K)e_i=\sum_{i=1}^k\alpha_i(u_2+ U_x\cap K)e_i$, i.e., $\alpha_i(u_1+ U_x\cap K)=\rho(u_1)+K_{z_i}=\alpha_i(u_2+ U_x\cap K)=\rho(u_2)+K_{z_i}$ for all $i=1,\dots, k$. Thus, $\rho(u_1)-\rho(u_2)=\rho(u_1-u_2)\in K_{z_i}$, i.e., $u_1-u_2\in K$. In addition, because $u_1,u_2$ belong to the $\Bbbk$-space $U_x$, we have that $u_1-u_2\in U_x$. Thus, $u_1-u_2\in U_x\cap K$, i.e., $u_1+U_x\cap K=u_2+U_x\cap K$. Consequently, we can identify the $\Bbbk$-space $\check{U}_x$ with the $\Bbbk$-subspace $\alpha(\check{U}_x)$ of $\check{\U}^\bullet$. In the sequel, we consider those two $\Bbbk$-spaces as the same. In this setting, $\check{U}_x$ is a subspace of $\check\U^\bullet$ for all $x\in\P\setminus\max\P$. \\

Now, we suppose that $x\prec y$ in $\P$. By hypothesis, we have that $\pi_y(U_x)\subseteq U_y$, where $\pi_y$ is the composition of the natural projection $U^{\bullet}\to U^{\bullet}_y$ and the natural injection $U^{\bullet}_y\to U^{\bullet}$. Let $\check{\pi}_y$ be the composition of the projection $\check\U^\bullet\to \check\U^\bullet_y$ and the injection $\check\U^\bullet_y\to \check\U^\bullet$, where $$\check\U^\bullet_y=\bigoplus_{\begin{smallmatrix} z\in\max\P \\ z\succeq y \end{smallmatrix}}U_z/K_z.$$ We have to prove that $\check{\pi}_y(\check{U}_x)\subseteq \check{U}_y$. Given an element $\alpha(w+U^K_x)=\sum_{i=1}^k\alpha_i (w+U_x^K)e_i$ in $\check{U}_x=\alpha(\check{U}_x)$, it is clear that $\check{\pi}_y(\sum_{i=1}^k\alpha_i (w+U_x^K)e_i)=\sum_{i=1}^k\alpha_i (w+U_x^K)\check{e}_i$, where $\check{e}_i=e_i$ if $z_i\succeq y$ and $\check{e}_i=(0,\dots,0)$ otherwise. But, $\sum_{i=1}^k\alpha_i (w+U_x^K)\check{e}_i=\sum_{i=1}^k(\rho_i(w)+K_{z_i})\check{e}_i=\pi_y(w)+\sum_{i=1}^k(K_{z_i})\check{e}_i\in \check{U}_y$ . \\

Finally, we suppose that $x\nprec z_i$. Then by hypothesis, $\pi_{z_i}(U_x)=0$. Given an element $\alpha(w+U_x^K)$ in $\check{U}_x=\alpha(\check{U}_x)$, then $\check{\pi}_{z_i}(\alpha(w+U_x^K))=\alpha_i(w+U_x^K)e_i$, where $\alpha_i(w+U_x^K)=\rho_i(w)+K_{z_i}$. Thus, $\alpha_i(w+U_x^K)e_i=\pi_{z_i}(w)+K_{z_i}e_i$. But, the hypothesis implies that $\pi_{z_i}(w)=0$, then $\check{\pi}_{z_i}(\alpha(w+U_x^K))=0$.
\end{proof}

In what follows, the peak $\P$-subspaces of a peak $\P$-space $\U$ of the form $\U_K$, where $K$ is an admissible subspace of $\U^\bullet$, will be called \textit{proper peak $\P$-subspaces} of $\U$.

\subsection{(Semi)stable representations.}\label{subsection:semistable representations}

Recall that the \textit{dimension vector} $\dimv\U$ of a peak $\P$-space $\U$  is a $\mathbb{Z}$-function on $\P$ given by $(\dimv \U)_x=\dim U_x$, that is, the dimension vector of $\U$ is an element of $\mathbb{Z}^{\P}$. The map $\U \longmapsto \dimv \U$ gives rise to an isomorphism between the Grothendieck group $K_0(\P\text{-spr})$ of the category $\P$-spr and $\mathbb{Z}^{\P}$ \cite{simson3}*{Lemma 2.2}. Following \cite{Rud}*{Section 3} one of the most usual ways to define stability (in an abelian category $\mathcal A$) is via ratio of two additive functions (i.e. functions on  $K_0(\mathcal A)$). Proceeding in a similar way, let $\thv$ and $\kpv$ be two additive functions on $\mathbb{Z}^\P$ and let $\kpv(\U)>0$ for any nonzero object $\U \in \P\text{-spr}$. We call the ratio 
$$
    \mu(\U)=\thv(\U)/\kpv(\U),
$$
as the $(\thv:\kpv)$-\textit{slope} of $\U$ and define $\U$ to be $\mu$-\textit{stable} (resp. $\mu$-\textit{semistable}) if for each proper subrepresentation $\U_K$ of $\U$ we have that 
$$
     \mu(\U_K)<\mu(\U) \qquad (\mbox{resp.}\  \mu(\U_K)\leq \mu(\U)).
$$

Fixing a weigh $\thv \in \Hom(\mathbb{Z}^{\P},\mathbb{Z})$ we say that $\U\in \P\text{-spr}$ is \textit{$\thv$-stable} (resp. $\thv$-\textit{semistable})  if $\thv(\U)=0$ and 
$$\thv(\U_K)<0\qquad (\mbox{resp.}\ \leq)$$ for any proper peak $\P$-subspace 
$\U_K$ of $\U$.\\

In fact, fixing a basis in $\Hom(\mathbb{Z}^{\P},\mathbb{Z})$, we will regard $\thv$ and $\kpv$ as the vectors $\thv=(\theta_x)_{x\in \P}$ and $\kpv=(\kappa_x)_{x\in \P}$ in $\mathbb{Z}^{\P}$, so that the $(\thv,\kpv)$-slope of $\U$ is given by
$$
	\mu(\U)=\dfrac{\sum_{x\in \P} \theta_x \dim U_x}{\sum_{x\in \P} \kappa_x \dim U_x}.
$$
The relationship between the $\thv$-stability and the $\mu$-stability is shown in the following lemma. 
\begin{lemma} Given a slope $\mu=\frac{\thv}{\kpv}$ and $\U$ a nonzero objet in $\P$-spr. Then
\begin{itemize}
\item[(i)] If $\U$ is $\thv$-stable (resp. $\thv$-semistable) then $\U$ is $\mu$-stable (resp. $\mu$-semistable).
\item[(ii)] If $\U$ is $\mu$-stable (resp. $\mu$-semistable) then $\U$ is $\widehat{\thv}$-stable (resp. $\widehat{\thv}$-semistable), where $\widehat{\thv}=\kpv(\U)\thv-\thv(\U)\kpv$. 
\end{itemize}
\end{lemma}
\begin{proof} We prove only the case of  stability. The argument for semistability follows a similar line of reasoning. The $\thv$-stability of $\U$ implies that $\mu(\U)=0$ and for any proper peak $\P$-subspace $\U_K$ of $\U$, we have that $\thv(\U_K)<0$ and $\kpv(\U_K)>0$, i.e., 
$\mu(\U_K)<0$. Then $\U$ is $\mu$-stable. To prove (ii), it is clear that $\widehat{\thv}(\U)=0$. Moreover, for a proper peak $\P$-subspace $\U_K$ of $\U$, we have that $\mu(\U_K)=\frac{\thv(\U_K)}{\kpv(\U_K)}<\mu(\U)=\frac{\thv(\U)}{\kpv(\U)}$. Since $\kpv$ is positive on any nonzero object in $\P$-spr, $\widehat{\thv}(\U_K)=\kpv(\U)\thv(\U_K)-\thv(\U)\kpv(\U_K)<0$. 
\end{proof}

The following is known as See-Saw property, it is proved using the same arguments given  in \cite{Hille}*{Lemma  2.1}.
 
\begin{proposition} \label{seesaw}
Let $\ses{\W}{\V}{\U}$ be an exact sequence of peak $\P$-spaces and let $\mu$ be a slope. Then the following conditions are equivalent:
\begin{enumerate}
	\item[(1)]	$\mu(\W)\leq \mu(\V)$,
	\item[(2)]  $\mu(\W)\leq \mu(\mathbf{U})$,
	\item[(3)]  $\mu(\V)\leq \mu(\mathbf{U})$.
\end{enumerate}
\end{proposition}

Moreover, following the same ideas of \cite{futorny}*{Propositions 2.4 and 2.5}, we have the version of the Harder-Narasimhan filtration and the Jordan-Hölder Filtration in this case. 
\begin{proposition}[Harder-Narasimhan filtration] For any peak $\P$-space $\V=(V_x)_{x\in\P}$  there is a unique filtration (of vector subspaces)
$$
	0=K^0\subset K^1\subset \dots \subset K^{h}=V^\bullet ,	
$$
which induces a filtration of $\V$
\begin{equation*} 
	0=\V^{0}\subset \V^{1}\subset \dots \subset \V^{h}=\V,	
\end{equation*}
where $\V^{i}=\V_{K^i}$,  such that:
\begin{itemize}
\item[(1)] $\V^{i}/\V^{i-1}$ are $\mu$-semistable, and
\item[(2)] $\mu(\V^{i}/\V^{i-1})>\mu(\V^{i+1}/\V^{i})$ for all $i=1,\dots,h-1$.
\end{itemize} 
\end{proposition}

\begin{proposition}[Jordan-H\"{o}lder filtration] \label{JordanHolder}
For any $\mu$-semistable $\V=(V_x)_{x\in\P}$, there exists a filtration (of vector subspaces)
$$
0=K^0\subset K^1\subset \dots \subset K^{h}=V^\bullet,	
$$
which induces a filtration of $\V$
\begin{equation*} 
0=\V^{0}\subset \V^{1}\subset \dots \subset \V^{h}=\V,	
\end{equation*}
where $\V^{i}=\V_{K^i}$, satisfying the following conditions:
\begin{itemize}
\item[(1)] $\V^{i}/\V^{i-1}$ are $\mu$-stable, and
\item[(2)] $\mu(\V^{i}/\V^{i-1})=\mu(\V^{i+1}/\V^{i})$ for all $i=1,\dots,h-1$.
\end{itemize}	 
\end{proposition}

\textbf{Positive Stability}: We define a peak $\P$-space $\V$ to be \textit{positively stable} if there exists a form $\thv\in \mathbb{Z}^{\P}$ such that $\theta_x>0$ for all $x\in \P\setminus\max\P$ and $\V$ is $\thv$-stable. It's important to note that if a peak $\P$-space is $\thv$-stable, it doesn't necessarily imply that $\theta_x>0$ for all $x\in \P\setminus\max\P$.
\begin{example}\label{non-positive-weight}Let $\P$ be the three peak poset given by the Hasse diagram

\iftrue

\begin{center}
\begin{tikzpicture}
\node (1) at (0,0) {$\star$};
\node [left  of=1, node distance=0.15cm] (l2)  {$_{_2}$};
\node [right  of=1, node distance=0.6cm] (2)  {$\star$};
\node [right  of=2, node distance=0.18cm] (l2)  {$_{_5}$};

\node [below  of=1, node distance=0.6cm] (1')  {$\circ$};
\node [left  of=1', node distance=0.15cm] (l1)  {$_{_1}$};

\node [below  of=2, node distance=0.6cm] (2')  {$\circ$};
\node [right  of=2', node distance=0.18cm] (l4)  {$_{_4}$};

\node [ below  of=1', node distance=0.6cm] (1'')  {$\circ$};
\node [left  of=1'', node distance=0.15cm] (l3)  {$_{_3}$};
\node [below  of=2', node distance=0.6cm] (2'')  {$\circ$};
\node [right  of=2'', node distance=0.18cm] (l2)  {$_{_6}$};
\node [right  of=2', node distance=0.6cm] (3)  {$\star$};
\node [right  of=3, node distance=0.18cm] (l2)  {$_{_7}$};

\draw [shorten <=-2pt, shorten >=-2pt] (1) -- (1');
\draw [shorten <=-2pt, shorten >=-2pt] (2) -- (2');

\draw [shorten <=-2pt, shorten >=-2pt] (1') -- (1'');
\draw [shorten <=-2pt, shorten >=-2pt] (2') -- (2'');
\draw [shorten <=-2pt, shorten >=-2pt] (1'') -- (2');
\draw [shorten <=-2pt, shorten >=-2pt] (2'') -- (3);
\end{tikzpicture}
\end{center}

\fi

and the peak $\P$-space $\U$ whose dimension vector is $\dimv \U=(1,1,1,1,1,1,0)$. The proper peak $\P$-subspaces are $\V$ and $\W$, where $\dimv\V=(1,1,0,0,0,0,0)$ and $\dimv\W=(0,0,0,1,1,1,0)$ (see Example \ref{propersubspaces}). Then $\U$ is $(1,-2,2,1,-1,-1,0)$-stable but the form is not positive. 

\end{example}
We need the following extension of stability. Given a $\thv$-stable peak $\P$-space $\V$ of a poset $\P$ and any $\widetilde \P$-space $\widetilde \V$  such that $\P$ is a peak subposet of $\widetilde \P$ and $\widetilde \V\big\vert_{\P}=\V$, define
  $\widetilde \theta_x=\theta_x$ if $x\in \P$ and $\widetilde \theta_x=0$ otherwise. Obviously $\widetilde \V$ is $\widetilde \thv$-stable. We prove even a stronger connection.
\begin{proposition} If	$\V$ be a positively stable peak $\P$-space  with form $\thv$. Any peak $\widetilde \P$-space $\widetilde \V$ such that 
	$\P$ is a peak-subposet of  $\widetilde \P$ and $\widetilde \V\big\vert_{\P}=\V$ 
is positively stable with some form $\widetilde \thv$.
\end{proposition}

\begin{proof}
By induction, it is enough to prove the statement for the case when $\P=\widetilde \P \setminus  \{\tilde p\}$.  If  $\tilde{p}$ is a maximal point in $\widetilde \P$, we can consider $\widetilde\theta_i=\theta_i$ for all $i\in \P$ and 
$\widetilde\theta_{\tilde{p}}=0$. Now, let us consider the case when $\tilde{p}$ is a non-maximal point in $\widetilde \P$. Let  $\V=(V_x)_{x\in \P}$ be a positively stable peak $\P$-space with form $\thv$. It is clear that  $\widetilde \V$ can be viewed as $\widetilde \V=(V_{\tilde p},V_x)_{x\in \P}$.  Defining $\thv'=(\dim V_{\tilde p}\dim \V^{\bullet}+1)\cdot \thv$, we have that $\V$ is $\thv'$-stable. Moreover,  for any proper peak $\P$-subspace $\W$ de $\V$ we have that $\thv(\dimv \W)\leq -1$. Thus, $$\thv'(\dimv \W)=(\dim V_{\tilde p}\dim \V^{\bullet}+1)\cdot \thv(\dimv \W)\leq-\dim V_{\tilde p} \dim\V^{\bullet}-1.$$
Set $\nu=-\dim V_{\tilde p} \dim\V^{\bullet} -1$ and define a form $\widetilde \theta$ as  follows: 
$$\widetilde \theta_z=-\nu \theta_z-\dim V_{\tilde p}, \qquad
  \widetilde \theta_x=\left\{\begin{array}{c}
  		\theta'_x,\quad  x\neq \tilde p\\ 
  		\dim \V^{\bullet},\quad  x=\tilde p.
  \end{array}\right.
$$
Then

\begin{align*}
\widetilde \thv(\dimv \widetilde \V)&= \sum_{z\in\max \widetilde\P}\widetilde \theta_z\dim \widetilde V_z + \sum_{x\in \widetilde\P\setminus \max \widetilde\P}\widetilde \theta_x\dim \widetilde V_x\\ &=-\nu\sum_{z\in\max\P}\theta_z\dim V_z-\dim V_{\tilde p}\dim V^{\bullet} + \sum_{x\in\P\setminus\max\P}\theta_x\dim V_x + -\dim V^{\bullet}\dim V_{\tilde p}\\&=\thv'(\dimv \V)=0,
\end{align*}

and for any proper peak $\P$-subpace $\widetilde \W=(W_{\tilde p},W_x)_{x\in \P}$ we have 
\begin{align*}
	\widetilde \thv(\dimv \widetilde \W)&=\thv'(\dimv \W)+\dim V^{\bullet}\dim{ W_{\tilde p}}-\dim V_{\tilde p}\dim{ W^{\bullet}}\\
	&\leq -\dim V_{\tilde p} \dim\V^{\bullet}-1 +\dim V^{\bullet}\dim{W_{\tilde p}}-\dim V_{\tilde p}\dim{ W^{\bullet}}
\end{align*}
Since $\dim V^{\bullet}\dim W_{\tilde p}\leq \dim V^{\bullet}\dim V_{\tilde p}$ and $\dim V_{\tilde p}\dim W^{\bullet}\geq 0$ then  $\dim V^{\bullet}\dim W_{\tilde p}-\dim V_{\tilde p}\dim W^{\bullet}\leq \dim V^{\bullet}\dim V_{\tilde p}$. Thus, $\widetilde \thv(\dimv \widetilde \W) \leq -\dim V_{\tilde p} \dim\V^{\bullet}-1+\dim V_{\tilde p}\dim V^{\bullet}=-1 <0. $
\end{proof}

\section{Stability arising from of a bilinear form} \label{section2}

We recall  that a full subposet $\P'$ of $\P$  is said to be a \textit{peak-subposet} if $\max\P'\subseteq\max\P$. Moreover, a finite connected poset $\P$ is said to be \textit{poset of type $\mathbb{A}$} if $\P$ does not contain any of the following posets as a peak-subposet:\\

\begin{center}
\begin{tabular}{|l|l|l|l|}
\hline 
\begin{tikzpicture}
\node (top) at (-0.7,0.5) {$\mathcal{R}_1$};

\node (top) at (0,0) {$\star$};
\node [below of=top, node distance=0.7cm](center) {$\circ$};
\node [left  of=center, node distance=0.7cm] (left)  {$\circ$};
\node [right of=center, node distance=0.7cm] (right) {$\circ$};

    \draw [blue,  thick, shorten <=-2pt, shorten >=-2pt] (top) -- (left);
    \draw [blue, thick, shorten <=-2pt, shorten >=-2pt] (top) -- (right);
    \draw [blue, thick, shorten <=-2pt, shorten >=-2pt] (top) -- (center);\end{tikzpicture} &
    
    \begin{tikzpicture} 
    \node (top) at (0,0.5) {$\mathcal{R}_2$};
    \node (top1) at (0,0) {$\star$};
    \node [right of=top1, node distance=0.7cm](top2)  {$\star$};
    \node (m) at (0.35,-0.35) {$\circ$};
    \node (bellow) at (0.35,-0.7) {$\circ$};
     \draw [blue,  thick, shorten <=-2pt, shorten >=-2pt] (top1) -- (m);
    \draw [blue, thick, shorten <=-2pt, shorten >=-2pt] (top2) -- (m);
    \draw [blue, thick, shorten <=-2pt, shorten >=-2pt] (m) -- (bellow);
    
   \end{tikzpicture}&
   
    \begin{tikzpicture} 
\node (top) at (0,0.5) {$\mathcal{R}_3$};   
   \node (top3) at (0,0){$\star$};    
    \node [right of=top3, node distance=0.7cm](top4) {$\star$};
        \node [right of=top4, node distance=0.7cm] (top5) {$\star$};
        \node [below of=top4, node distance=0.7cm] (be) {$\circ$};
 \draw [blue, thick, shorten <=-2pt, shorten >=-2pt] (be) -- (top3);
 \draw [blue, thick, shorten <=-2pt, shorten >=-2pt] (be) -- (top4);
 \draw [blue, thick, shorten <=-2pt, shorten >=-2pt] (be) -- (top5);\end{tikzpicture} & 
 
 \begin{tikzpicture}
 \node (1) at (0,0) {$\star_1$};
\node (A) at (0.8,0.5) {$\mathcal{R}_{4,n}$, $n\geq 0$};
\node [right of=1, node distance=0.7cm](2) {$\star_2$};

\node [right of=2, node distance=0.7cm](3) {$\star_3$};
\node [right of=3, node distance=0.7cm] (4)  {$\cdots$};
\node [right of=4, node distance=0.7cm] (5) {$\star_{n+2}$};

\node [below of=1, node distance=0.7cm] (1')  {$\circ$};
\node [below  of=2, node distance=0.7cm] (2')  {$\circ$};
\node [below  of=3, node distance=0.7cm] (3')  {$\circ$};
\node [below  of=4, node distance=0.7cm] (4')  {$\dots$};
\node [below  of=5, node distance=0.7cm] (5')  {$\circ$};

   \draw [blue,  thick, shorten <=-2pt, shorten >=-2pt] (1) -- (1');
    \draw [blue, thick, shorten <=-2pt, shorten >=-2pt] (2') -- (2);
    
\draw [blue,  thick, shorten <=-2pt, shorten >=-2pt] (3) -- (3');
    \draw [blue, thick, shorten <=-2pt, shorten >=-2pt] (5) -- (5');    
    
    \draw [blue, thick, shorten <=-2pt, shorten >=-2pt] (1) -- (2');
    
 \draw [blue, thick, shorten <=-2pt, shorten >=-2pt] (2) -- (3');
     \draw [blue, thick, shorten <=-2pt, shorten >=-2pt] (1') -- (5);
\end{tikzpicture} \\ 
\hline 
\end{tabular} 

\end{center}
For instance, the poset given in Example \ref{non-positive-weight} is a three-peak poset of type $\mathbb{A}$. It was shown in \cite{schifflerserna} that when $\P$ is a poset of type $\mathbb{A}$,  the category $\P$-spr is of finite representation type, this is, the number of isoclasses of indecomposable peak $\P$-spaces is finite. Moreover,  any indecomposable peak $\P$-spaces $\U$ is such that $U_x\simeq\Bbbk$ for each $x\in\Supp \U$, where $\Supp \U=\{x\in\P\mid U_x\neq 0\}$ is called the \textit{support} of $\U$. \\

When $\P$ is a poset of type $\mathbb{A}$, an easy way to calculate all proper peak $\P$-subspaces of a peak $\P$-space $\U$ is provided by the following two lemmas.

\begin{lemma}\label{bijectionproper}
Let $\P$ be a poset of type $\mathbb{A}$, and let $\U$ be a peak $\P$-space. There is a bijection between the proper subsets of $\Supp\U \cap \max\P$ and the proper peak $\P$-subspaces of $\U$.
\end{lemma}

\begin{proof}
Given a proper subset $I$ of $\Supp\U \cap \max\P$, we define 
$$ 
K_I = \bigoplus_{z \in \max\P} K_z 
$$
such that $K_z = \Bbbk$ if $z \in I$ and $K_z = 0$ otherwise. Since $\P$ is a poset of type $\mathbb{A}$, and $U_x = \Bbbk$ for each $x \in \Supp\U$, it follows that $K_I$ is an admissible subspace of $U^{\bullet}$. Hence, $\U_{K_I}$ is a proper peak $\P$-subspace of $\U$. Denote by $\mathcal{X}$ the set of proper subsets of $\Supp\U \cap \max\P$. We define a function $\phi: \mathcal{X} \to \P$-spr such that $\phi(I) = \U_{K_I}$. Clearly, $\phi$ is injective because $K_I = K_J$ implies $I = J$. Moreover, if $\V$ is a proper peak $\P$-subspace of $\U$, then $\V$ has the form $\U_K$, where $K = V^{\bullet}$. The set $I = \{z \in \max\P \mid V_z \neq 0\}$ satisfies $\phi(I) = \V$. Thus, $\phi$ is bijective, and the proof is complete.
\end{proof}

For a point $a \in \P$, we denote by $a_{\vartriangle}$ the corresponding principal ideal, that is, the set $\left\{x \in \mathscr{P} \mid x \preceq a\right\}$. If $X \subseteq \P$, we denote by $X_{\vartriangle}$ the set $\bigcup_{a \in X} a_{\vartriangle}$.

\begin{lemma} \label{Supp} 
If $\U_{K_I}$ is defined as in Lemma \ref{bijectionproper}, then 
$$
\Supp (\U_{K_I}) = I_{\vartriangle} \setminus (I^{c})_{\vartriangle}.
$$
\end{lemma}

\begin{proof}
Let $x$ be an element in $\Supp(\U_{K_I})$, that is, $(\U_{K_I})_x = U_x \cap K_I \neq 0$. Then there exists an element $y = \sum_{z \in \max\P} y_z e_z \neq 0$ in $U^{K_I}_x$. Hence, $y_{z_i} \neq 0$ for some $z_i \in \max\P$. As a result, $K_{z_i} = \Bbbk$, which implies $z_i \in I$. Thus, $x \in I_{\vartriangle}$. On the other hand, if $x \in (I^c)_{\vartriangle}$, then there exists $z_j \notin I$ such that $x \preceq z_j$. In this case, $U_x = \langle u \rangle$, where $p_{z_i}(u) = p_{z_j}(u) \neq 0$ and $p_z(u) = 0$ for all $z \in \max\P \setminus \{z_i, z_j\}$. Since $\P$ is a poset of type $\mathbb{A}$, there are at most two maximal points greater than $x$. However, $u \notin K_I$ because $p_{z_j}(u) \neq 0$, leading to a contradiction. Thus, $x \in I_{\vartriangle} \setminus (I^{c})_{\vartriangle}$. Now, suppose that $x \in I_{\vartriangle} \setminus (I^{c})_{\vartriangle}$. This means $x \npreceq z'$ for any $z' \in (I^c)_{\vartriangle}$. Furthermore, $\pi_{z'}(U_x) = 0$ for each $z' \in (I^c)_{\vartriangle}$. Hence, $U_x = \langle u \rangle$ with $p_{z'}(u) = 0$ for all $z' \in (I^c)_{\vartriangle}$. Since $u \in K_I$, we conclude that $U_x$ is a subspace of $K_I$. Consequently, $U_x \cap K_I = \Bbbk$, implying $x \in \Supp(\U_{K_I})$.
\end{proof}

\begin{example} \label{propersubspaces}
Let $\P$ and $\U$ be as in Example \ref{non-positive-weight}. Then $\Supp \U = \{1, 2, 3, 4, 5, 6\}$, and the proper subsets of $\Supp\U \cap \max\P = \{2, 5\}$ are: $\{2\}$ and $\{5\}$. Thus, Lemma \ref{Supp} implies that the two proper peak $\P$-subspaces of $\U$ have supports $\{1, 2\}$ and $\{4, 5, 6\}$.
\end{example}

We recall from \cite{simson3} that to any peak $\P$-space $\U = (U_x)_{x \in \P}$, one associates the \textit{coordinate vector} $\cdnv\U$ defined as 
$$
(\cdnv\U)_{x} =
\begin{cases}
\dim U_x,  & \text{if } x \in \max\P, \\
\dim \left(U_x \big/ \sum_{y \prec x} \pi_{x}(U_y)\right), & \text{if } x \in \P \setminus \max\P.
\end{cases}
$$
By the \textit{coordinate support} of $\U$, we mean the peak-subposet 
$$
\csupp (\U) = \{ x \in \P \mid (\cdnv\U)_{x} \neq 0 \}
$$
of $\P$. Moreover, the peak $\P$-space $\U$ is said to be \textit{sincere} if it is indecomposable and $\csupp \U = \P$. Furthermore, if there exists a sincere $\P$-space, we say that $\P$ is a \textit{sincere poset}. The classification of all sincere $r$-peak posets of finite representation type (posets $\P$ with only a finite number of nonisomorphic indecomposable $\P$-spaces) was given by M. Kleiner \cite{kleiner75} for the case $r = 1$, and by J. Kosakowska \cites{justina,justina1,justina2} for the remaining cases. Particularly, a poset $\P$ of type $\mathbb{A}$ is a sincere poset if and only if $\P$ is isomorphic to one of the following posets: \vspace{0.2cm}

\begin{table}[ht]  \label{table1}
\begin{adjustbox}{max totalsize={1.0\textwidth}{1.0\textheight},center}
\begin{tabular}{|l|l|l|}
\hline 
\begin{tikzpicture}
 \node (1) at (0,0) {$z_1$};
 
\node (A) at (1.55,0.5) {$\mathcal{S}^{(r)}_{1}$};
\node [right of=1, node distance=0.7cm](2) {$z_2$};

\node [right of=2, node distance=0.7cm](3) {$z_3$};
\node [right of=3, node distance=0.7cm] (4)  {$\cdots$};
\node [right of=4, node distance=0.7cm] (5) {$z_{r}$};

\node [below  of=2, node distance=0.7cm] (2')  {$x_{1}$};
\node [below  of=3, node distance=0.7cm] (3')  {$x_{2}$};
\node [below  of=4, node distance=0.7cm] (4')  {$\dots$};
\node [below  of=5, node distance=0.7cm] (5')  {$x_{r-1}$};

   
    \draw [blue, thick, shorten <=-2pt, shorten >=-2pt] (2') -- (2);
    
\draw [blue,  thick, shorten <=-2pt, shorten >=-2pt] (3) -- (3');
    \draw [blue, thick, shorten <=-2pt, shorten >=-2pt] (5) -- (5');    
    
    \draw [blue, thick, shorten <=-2pt, shorten >=-2pt] (1) -- (2');
    
 \draw [blue, thick, shorten <=-2pt, shorten >=-2pt] (2) -- (3');
    
\end{tikzpicture} & 

\begin{tikzpicture}
 \node (1) at (0,0) {$z_1$};
 
\node (A) at (1.55,0.5) {$\mathcal{S}^{(r)}_{2}$};
\node [right of=1, node distance=0.7cm](2) {$z_2$};

\node [right of=2, node distance=0.7cm](3) {$z_3$};
\node [right of=3, node distance=0.7cm] (4)  {$\cdots$};
\node [right of=4, node distance=0.7cm] (5) {$z_{r}$};

\node [below  of=2, node distance=0.7cm] (2')  {$x_{1}$};
\node [below  of=3, node distance=0.7cm] (3')  {$x_{2}$};
\node [below  of=4, node distance=0.7cm] (4')  {$\dots$};
\node [below  of=5, node distance=0.7cm] (5')  {$x_{r-1}$};
\node [right of=5', node distance=0.7cm] (6')  {$x_{r}$};

   
    \draw [blue, thick, shorten <=-2pt, shorten >=-2pt] (2') -- (2);
    
\draw [blue,  thick, shorten <=-2pt, shorten >=-2pt] (3) -- (3');
    \draw [blue, thick, shorten <=-2pt, shorten >=-2pt] (5) -- (5');    
    
    \draw [blue, thick, shorten <=-2pt, shorten >=-2pt] (1) -- (2');
    
 \draw [blue, thick, shorten <=-2pt, shorten >=-2pt] (2) -- (3');
    
\draw [blue, thick, shorten <=-2pt, shorten >=-2pt] (6') -- (5);

\end{tikzpicture}

 & 
 
\begin{tikzpicture}
 \node (1) at (0,0) {$z_1$};
 
\node (A) at (1.55,0.5) {$\mathcal{S}^{(r)}_{3}$};
\node [right of=1, node distance=0.7cm](2) {$z_2$};

\node [right of=2, node distance=0.7cm](3) {$z_3$};
\node [right of=3, node distance=0.7cm] (4)  {$\cdots$};
\node [right of=4, node distance=0.7cm] (5) {$z_{r}$};

\node [below of=1, node distance=0.7cm] (1')  {$x_{0}$};
\node [below  of=2, node distance=0.7cm] (2')  {$x_{1}$};
\node [below  of=3, node distance=0.7cm] (3')  {$x_{2}$};
\node [below  of=4, node distance=0.7cm] (4')  {$\dots$};
\node [below  of=5, node distance=0.7cm] (5')  {$x_{r-1}$};
\node [right of=5', node distance=0.7cm] (6')  {$x_{r}$};

   \draw [blue,  thick, shorten <=-2pt, shorten >=-2pt] (1) -- (1');
    \draw [blue, thick, shorten <=-2pt, shorten >=-2pt] (2') -- (2);
    
\draw [blue,  thick, shorten <=-2pt, shorten >=-2pt] (3) -- (3');
    \draw [blue, thick, shorten <=-2pt, shorten >=-2pt] (5) -- (5');    
    
    \draw [blue, thick, shorten <=-2pt, shorten >=-2pt] (1) -- (2');
    
 \draw [blue, thick, shorten <=-2pt, shorten >=-2pt] (2) -- (3');
    
\draw [blue, thick, shorten <=-2pt, shorten >=-2pt] (6') -- (5);

\end{tikzpicture} 
 
  \\ 
\hline 
\end{tabular} 
\end{adjustbox}
\caption{Sincere posets of type $\mathbb{A}$.} \label{sincereposets}
\end{table}
for some $r \geq 1$. In this case, the sincere peak $\P$-space $\U$ is such that $U_x$ is the $\Bbbk$-space generated by the vector given by the sum of the standard vectors $e_z$ in $\U^{\bullet}$, where $z$ runs over all of the maximal points in $\P$ greater than $x$. In other words, $\U_x \simeq \Bbbk$. \\

The lists of sincere posets and their sincere representations are important because we can obtain all indecomposable objects in $\P$-spr for a given poset $\P$ by lifting all sincere objects of all sincere peak-subposets $\mathcal{S}$ of $\P$ via the well-known \textit{subposet induced functor} $T_{\mathcal{S}} : \mathcal{S}\text{-spr} \to \P\text{-spr}$ \cite{dowbor}, which assigns to the peak $\mathcal{S}$-space $\U$ the peak $\P$-space $\widehat{\U}$, where $\widehat{U}_x$ is defined by
$$
    \widehat{U}_{x} =
    \begin{cases}
    U_x,  & \text{if } x \in \max \mathcal{S}, \\
    0, & \text{if } x \notin (\max \mathcal{S})_{\vartriangle}, \\
    \sum_{x \succeq y \in \mathcal{S}} \pi_x(U_y), & \text{if } x \in (\max \mathcal{S})_{\vartriangle} \setminus \max \mathcal{S}.
    \end{cases}
$$
In other words, up to isomorphism, any indecomposable peak $\P$-space $\U$ is the image $T_{\mathcal{S}}(\V)$ of a sincere $\mathcal{S}$-space, where $\mathcal{S}$ is a sincere peak-subposet of $\P$. In this case, $\mathcal{S}$ is the poset determined by $\csupp \U$.

\begin{proposition} \label{thetalifting} 
Let $\U$ be a sincere peak $\mathcal{S}$-space, where $\mathcal{S}$ is a peak-subposet of a type $\mathbb{A}$ poset $\P$. If $\U$ is $\thv$-stable, then the peak $\P$-space $T_{\mathcal{S}}(\U)$ is $\widetilde{\thv}$-stable for the form
$$
\widetilde{\theta}_x =
\begin{cases}
\theta_x, & \text{if } x \in \mathcal{S}, \\
0, & \text{if } x \notin \mathcal{S}.
\end{cases}
$$
\end{proposition}

\begin{proof}
The definition of $\widetilde{\thv}$ implies that if $T_{\mathcal{S}}(\U) = \widehat{\U}$, then 
$$
\widetilde{\thv}(\widehat{\U}) = \sum_{x \in \mathcal{S}} \theta_x \dim \widehat{U}_x.
$$
By definition, $\dim \widehat{U}_x = \dim U_x$ if $x \in \max \mathcal{S}$. Let $x \in \mathcal{S} \setminus \max \mathcal{S}$. The $\P$-space $\widehat{\U}$ is such that $\widehat{U}_x = \pi_x(U_x)$. Since $\U$ is sincere and $\mathcal{S}$ is of type $\mathbb{A}$, $x$ is a minimal point in $\mathcal{S}$. It is known from \cite{schifflerserna} that $U_x$ is generated by the vector $u_x = \sum_{x \prec z \in \max \mathcal{S}} e_z$. Since $\max \mathcal{S} \subseteq \max \P$, we have that $\pi_x(u_x) = u_x$. Therefore, $\pi_x(U_x) = U_x$, which implies that $\dim \widehat{U}_x = \dim U_x$. Thus, 
$$
\widetilde{\thv}(T_{\mathcal{S}}(\U)) = \thv(\U) = 0.
$$

Let $\V = (V_x)_{x \in \P}$ be a proper peak $\P$-subspace of $\widehat{\U}$. Note that 
$$
\widetilde{\thv}(\V) = \sum_{x \in \mathcal{S}} \theta_x \dim V_x.
$$
Clearly, the ambient space $\widehat{U}^{\bullet} = U^{\bullet}$, so there exists a proper subset $I$ of $\max \mathcal{S}$ such that $\V = \widehat{\U}_{K_I}$, where $K_I$ is the admissible $\Bbbk$-subspace of $U^{\bullet}$ associated with $I$ (see Lemma \ref{bijectionproper}). Suppose that $\W$ is the peak $\P$-space $T_{\mathcal{S}}(\U_{K_I})$. We will prove that $\dim W_x = \dim V_x$ for each $x \in \mathcal{S}$. In fact, let $x$ be a point in $\mathcal{S}$. If $x$ is a maximal point in $\mathcal{S}$, then $\widehat{U}_x = U_x$ and $W_x = U_x \cap K_I$. Therefore, 
$$
V_x = \widehat{U}_x \cap K_I = U_x \cap K_I = W_x.
$$
In the case when $x \notin \max \mathcal{S}$, we have that $\widehat{U}_x = \sum_{x \succeq y \in \mathcal{S}} \pi_x(U_y)$. Since $x$ is a minimal point in $\mathcal{S}$, we have $\widehat{U}_x = \pi_x(U_x) = U_x$. Therefore,
$$
V_x = \widehat{U}_x \cap K_I = \pi_x(U_x) \cap K_I = \pi_x(U_x \cap K_I) = W_x.
$$
Thus, $\widetilde{\thv}(\V) = \widetilde{\thv}(\W) = \widetilde{\thv}(T_{\mathcal{S}}(\U_{K_I})) = \thv(\U_{K_I}) < 0$.
\end{proof}

We recall that if $\P$ is a poset with $n$ points, the \textit{incidence matrix} $C_{\P}$ of $\P$ is an integral square $n \times n$ matrix $C_{\P} = [c_{xt}]_{x,t \in \P} \in \mathbb{M}_{\P}(\mathbb{Z})$, with
$$
c_{xt} = \begin{cases}
1, & \mbox{if } x \preceq t, \\
0, & \mbox{if } x \npreceq t.
\end{cases}
$$
Given a poset $\P$, the following bilinear form $b_{\P}$ plays a fundamental role (see \cites{simsonForm, CI19, kasjan1996tame} and references therein):
$$
b_{\P}: \mathbb{Z}^{\P} \times \mathbb{Z}^{\P} \rightarrow \mathbb{Z},
$$
$$
b_{\P}(\alpha, \beta) = \alpha \cdot C_{\P}^{-1} \cdot \beta^{tr} = \sum_{x \in \P} \alpha_x \beta_x + \sum_{t \prec x \in \P} c^{-}_{xt} \alpha_x \beta_t,
$$
where $c^{-}_{xt}$ is the $(x,t)$ entry of the matrix $C_{\P}^{-1} \in \mathbb{M}_{\P}(\mathbb{Z})$, the inverse of $C_{\P}$. In \cite{futorny}, the authors showed that when $\P$ is a one-peak poset, every indecomposable $\P$-space is stable with the weight described in \eqref{stableformula} (see \cite{futorny}*{Proposition 4.3}). We prove that such a weight also serves for the category of peak $\P$-spaces when $\P$ is an $r$-peak poset of type $\mathbb{A}$.

\begin{proposition}\label{thetasincere}
Let $\U$ be a sincere peak $\P$-space of a sincere poset $\P$ of type $\mathbb{A}_n$. Then  $\U$ is $\thv$-stable with the form
 \begin{equation}\label{stableformula}
\thv(\W) = b_{\P}(\dimv \U, \dimv \W) - b_{\P}(\dimv \W, \dimv \U).
\end{equation} 
\end{proposition}

\begin{proof}
To prove that sincere peak $\P$-space $\U$ of $\P$ is stable with the form \eqref{stableformula}, we first consider $\W_j$, the system of $\Bbbk$-vector spaces whose dimension $\dimv \W_j$ is the vector $e_j$ in the standard basis $\{e_1, \dots, e_n\}$ of $\Bbbk^n$. Note that $\W_j$ is a peak $\P$-space if and only if $j \in \max \P$. Thus, we have:

\begin{equation}\label{formulaproof}
\begin{array}{l}
\thv(\W_{x_{i}}) = \dimv \W_{z_{i}} + \dimv \W_{z_{i+1}}, \\
\thv(\W_{z_{i}}) = -\dimv \W_{x_{i-1}} - \dimv \W_{x_{i}},
\end{array}
\end{equation}

where $z_{i}$ is a maximal point and $x_{i}$ is a minimal point of $\P$, for $1 \leq j \leq n$ and $n \geq 1$ (see Table \ref{sincereposets}). Suppose that $\P$ is a sincere poset of type $\mathbb{A}$. If $\P$ is equal to $\mathcal{S}_{1}^{(r)}$, the incidence matrix $C_{\P}$ of $\P$ is
$$ C_{\P}= \left [
      \begin{array}{c|c}
      I_{r} & 0\\
      \hline
      J_{r-1}(1)_{\downarrow}\ & I_{r-1}
   \end{array}
   \right ] $$
   
where $J_{r-1}(1)_{\downarrow}$ is the $r-1 \times r$ matrix 

$$\left[
\begin{smallmatrix}
 J_{r-1}(1) & \begin{smallmatrix}
   \mathrm{0}\\ \mathrm{\vdots}\\  \mathrm{0}\\ \mathrm{1} \end{smallmatrix} \end{smallmatrix}\right]$$ with $J_{r-1}(1)$  a Jordan block, if $\dimv \U$ is the vector $$(\dim U_{z_{1}}, \dots, \dim U_{z_{r}}, \dim U_{x_{1}}, \dots \dim U_{x_{r-1}}),$$ we have that $$\thv(\W_{j})=b_{\P}(\dimv \U,\dimv \W_{j})-b_{\P}( \dimv \W_{j},\dimv \U),$$ in this case, if $\W_{j}=\W_{x_{i}}$,
\begin{equation}
\begin{array} {r l}
\thv(\W_{x_{i}})&= \dimv \U \cdot C_{\P}^{-1} \cdot (e_{x_{i}})^{tr}-e_{x_{i}} \cdot C_{\P}^{-1} \cdot (\dimv \U)^{tr}\\
&=\dimv \W_{x_{i}}-(-\dimv \W_{z_{i}}- \dimv \W_{z_{i+1}}+\dimv \W_{x_{i}})\\
&=\dimv \W_{z_{i}}+ \dimv \W_{z_{i+1}},
\end{array}
\end{equation} 
with $r+1\leq j \leq 2r-1$ and  $1\leq i \leq r-1$. In the same way , if $\W_{j}=\W_{z_{i}}$,

\begin{equation}
\begin{array} {r l}
\thv(\W_{z_{i}})&= \dimv \U \cdot C_{\P}^{-1} \cdot (e_{z_{i}})^{tr}-e_{z_{i}} \cdot C_{\P}^{-1} \cdot (\dimv \U)^{tr}\\
&=\dimv \W_{z_{i}}-\dimv \W_{x_{i-1}}- \dimv \W_{x_{i}}-\dimv \W_{z_{i}}\\
&=-\dimv \W_{x_{i-1}}- \dimv \W_{x_{i}},
\end{array}
\end{equation}

with $1\leq j \leq r$, $1\leq i \leq r$, and $\dimv \W_{x_{0}}=\dimv \W_{x_{r}}=0$. The same arguments can be used to obtain the proof in the cases when $\P$ is equal to $\mathcal{S}^{(r)}_{2}$ or  $\P$ is equal to $\mathcal{S}^{(r)}_{3}$. Note that, the incidence matrix of $\mathcal{S}^{(r)}_{2}$ is 

$$ \left [
      \begin{array}{c|c}
      I_{r} & 0\\
      \hline
      J_{r}(1)& I_{r}
   \end{array}
   \right ] $$ 
and incidence matrix of $\mathcal{S}^{(r)}_{3}$ is of the form 

$$ \left [
      \begin{array}{c|c}
      I_{r} & 0\\
      \hline
      (J_{r}(1)\downarrow)^{tr}& I_{r+1}
   \end{array}
   \right ] $$ 
In particular, if $\P= \mathcal{S}^{(r)}_{2}$ (resp. $\P= \mathcal{S}^{(r)}_{3}$), we have that $\dimv \W_{x_{0}}=\dimv \W_{z_{r+1}}=0$ (resp. $\dimv \W_{z_{0}}=\dimv \W_{z_{r+1}}=0$).\\

To prove that $\U$ is stable, we consider that $\U$ is the sincere peak $\P$-space , it is clear that $\thv(\U)=0$ by using identity \eqref{formulaproof}. Now, we suppose $\W$ is a 
proper peak $\P$-subspace of $\U$,  according Lemma \ref{bijectionproper}, there is a set $I$ with $\emptyset \neq I \subset \Supp \U \cap \max \P $ such that $\W=\U_{K_{I}}=(U_{x}^{K_{I}})_{x \in \P}$, where

\begin{equation}
U_{x}^{K_{I}} = 
   \begin{cases} 
         \Bbbk & \mbox{if } x \in I \cup A, \\ 
      0 & \mbox{otherwise, }  
   \end{cases}
\end{equation} 
 and $A=\lbrace x_{i} \in \min \P \text{ }| \text{ } z_{i}, z_{i+1} \in I \rbrace$. Note that, when $\P\simeq \mathcal{S}^{(r)}_{2}$ if $z_{r} \in I$  then $x_{r} \in A$, and when  $\P\simeq \mathcal{S}^{(r)}_{3}$ if $z_{1} \in I$ (resp. $z_{r} \in I$) then  $x_{0} \in A$ (resp. $x_{r+1} \in A$).  Given a set $ X=\lbrace x_{m_{1}},x_{m_{1}+1},\dots, x_{m_{t}}\rbrace$ of consecutive points in $A$ we define the set $Z_{X}=\bigcup_{x \in X} x^{\blacktriangledown}$, where $x^{\blacktriangledown}=\{y\in\P\mid y\succ x \}$. Thus, we establish a  system $\W^{X}$ of $\Bbbk$-vector spaces whose dimension vector $\dimv \W^X$ is
$\sum_{k \in X \cup Z_X } e_{k}$ (note that $\W^{X}$ is not necessarily a peak $\P$-space), then   $\dimv \W$ is the sum of all dimensions $\dimv \W^X$, where $X$  runs over all of subsets of consecutive points in $A$. Hence, it is enough to prove that $\thv(\W^X)<0$ for each set of consecutive points $X$ in $A$. To do that, we consider the following possibilities:\\

(i) If $\P=\mathcal{S}^{(r)}_{1}$, given a set $X$ as above, then  $Z_X=\lbrace z_{m_{1}},z_{m_{1}+1}, \dots, z_{m_{t}}, z_{m_{t}+1} \rbrace$   is such that if $z_{1}$ or $z_{r}$ belongs to $Z_X
$ then

\begin{equation*}
 \displaystyle \thv(\W^X)= \displaystyle \sum_{z_{m_{j}} \in Z_X} \thv(\W_{z_{m_{j}}}) + \sum_{ x_{m_{j}} \in X} \thv(\W_{x_{m_{j}}})=-(2(m_{t}+1)-1)+2(m_{t})=-1, 
\end{equation*}

if $z_{1}, z_{r}\notin Z_X$, we have

\begin{equation*}\label{case1stable}
\displaystyle \thv(\W^X)=  \displaystyle \sum_{z_{m_{j}} \in Z_X} \thv(\W_{z_{m_{j}}}) + \sum_{ x_{m_{j}} \in X} \thv(\W_{x_{m_{j}}})=-(2(m_{t}+1))+2(m_{t})=-2.   
\end{equation*}

(ii) If we consider $\P=\mathcal{S}^{(r)}_{2}$, the situation is similar to (i), but the point $x_{r}$ does not belong to $\mathcal{S}^{(r)}_{1}$, then if $x_{r} \in X$ or $x_{1} \in X$

\begin{equation*}
\displaystyle \thv(\W^X)= \displaystyle \sum_{z_{m_{j}} \in Z}  \thv(\W_{z_{m_{j}}}) + \sum_{ x_{m_{j}} \in X} \thv(\W_{x_{m_{j}}})=-(2(m_{t}))+2(m_{t}-1)=-1,
\end{equation*}
otherwise $\thv(\W^X)=-2$  as the case (i).\\

(iii) If $\P \simeq  \mathcal{S}^{(r)}_{3}$, the cases when $x_{0}$ (resp. $x_{r+1}$) belongs  to $X$ or $x_{0}, x_{r+1} \notin X$  satisfy the same conditions of the previous case, then $\thv(\W^X)<0$.
\end{proof}

\begin{theorem} \label{Firsttheorem}
Let $\P$ be a poset of type $\mathbb{A}$. Every indecomposable peak $\P$-space is $\thv$-stable with some weight $\thv$.
\end{theorem}

\begin{proof}
Any indecomposable peak $\P$-space $\U$ is (up to isomorphism) the image $T_{\mathcal{S}}(\V)$ of a sincere $\mathcal{S}$-space, where $\mathcal{S}$ is a sincere peak-subposet of $\P$. Thus, it is enough to apply Propositions \ref{thetalifting} and \ref{thetasincere}.
\end{proof}

\begin{example} Continuing with Example \ref{propersubspaces}, the peak $\P$-space $\widehat{\U}$ is such that $$\widehat{U}^{\bullet}=\widehat{U}_2\oplus \widehat{U}_5\oplus \widehat{U}_7= \Bbbk\oplus\Bbbk\oplus 0,\; \widehat{U}_1=\left<\widehat{u}\right>,\;  \widehat{U}_4=\widehat{U}_6=\left<\widehat{v}\right>, \text{ and } \widehat{U}_3=\left<\widehat{u}+\widehat{v}\right>, $$
where $\widehat{u}=(1,0,0)$ and $\widehat{v}=(0,1,0)$. Here, we will consider the order in the vector $\dimv \widehat \U=(\dim U_2,\dim U_5,\dim U_7,\dim U_1,\dim U_4,\dim U_3,\dim U_6)$. Note that, $\mathcal{S}=\csupp \widehat{\U}$ is given by the diagram

\begin{center}
\begin{tikzpicture}
\node (1) at (0,0) {$\star$};
\node [left  of=1, node distance=0.15cm] (l2)  {$_{_2}$};
\node [right  of=1, node distance=0.6cm] (2)  {$\star$};
\node [right  of=2, node distance=0.18cm] (l2)  {$_{_5}$};

\node [below  of=1, node distance=0.6cm] (1')  {$\circ$};
\node [left  of=1', node distance=0.15cm] (l3)  {$_{_3}$};
\node [below  of=2, node distance=0.6cm] (2')  {$\circ$};
\node [right  of=2', node distance=0.18cm] (l2)  {$_{_6}$};

\draw [shorten <=-2pt, shorten >=-2pt] (1) -- (1');
\draw [shorten <=-2pt, shorten >=-2pt] (2) -- (2');
\draw [shorten <=-2pt, shorten >=-2pt] (1') -- (2);

\end{tikzpicture}, 
\end{center}
which is a sincere peak-subposet  of $\P$. Moreover, the $\mathcal{S}$-space $\U$ is

$${U}^{\bullet}={U}_2\oplus {U}_5= \Bbbk\oplus\Bbbk,\; U_6=\left<v\right>, \text{ and } U_3=\left<u+v\right>, $$
with $u=(1,0)$ and $v=(0,1)$.  In this case, we calculate the proper peak $\mathcal{S}$-subspaces using Lemmas \ref{Supp} and \ref{bijectionproper},   those are $\W$ and $\W'$, where $\Supp \W=\{2\}$ and $\Supp \W'=\{5,6\}$, that is, $\dim\W=(1,0,0,0)$ and $\dim\W'=(0,1,0,1)$. Then, it is clear that $\thv(\U)=0$ and by \eqref{stableformula} we have:

$$\thv(\W)= \left(\begin{smallmatrix}
1 & 1&1&1 
\end{smallmatrix}\right)\left(\begin{smallmatrix}
 1&0&0&0\\0&1&0&0\\-1&-1&1&0\\ 0&-1&0&1
\end{smallmatrix}\right)\left(\begin{smallmatrix}
 1\\0\\0\\0
\end{smallmatrix}\right)-\left(\begin{smallmatrix}
1 & 0&0&0 
\end{smallmatrix}\right)\left(\begin{smallmatrix}
 1&0&0&0\\0&1&0&0\\-1&-1&1&0\\ 0&-1&0&1
\end{smallmatrix}\right)\left(\begin{smallmatrix}
 1\\1\\1\\1
\end{smallmatrix}\right)
=-1$$  
and 

$$\thv(\W')= \left(\begin{smallmatrix}
1 & 1&1&1 
\end{smallmatrix}\right)\left(\begin{smallmatrix}
 1&0&0&0\\0&1&0&0\\-1&-1&1&0\\ 0&-1&0&1
\end{smallmatrix}\right)\left(\begin{smallmatrix}
 0\\1\\0\\1
\end{smallmatrix}\right)-\left(\begin{smallmatrix}
0 & 1&0&1 
\end{smallmatrix}\right)\left(\begin{smallmatrix}
 1&0&0&0\\0&1&0&0\\-1&-1&1&0\\ 0&-1&0&1
\end{smallmatrix}\right)\left(\begin{smallmatrix}
 1\\1\\1\\1
\end{smallmatrix}\right)
=-1$$
Replacing $\dimv\W$ by the standard vectors $e_1,\dots,e_4\in\Bbbk^4$ in \eqref{stableformula}, we have that  $\thv=(-1,-2,2,1)$. Thus,  $\widetilde \thv= (-1,-2,0,0,2,1)$.
\end{example}

\section{A new geometric model for socle-projective categories (type $\mathbb{A}$)}  \label{section3}

\definecolor{myred}{cmyk}{0.5,0,0.5,0.5}
\definecolor{mylightblue}{cmyk}{0.24,0.06,0,0.19}
\def\myxscale{0.5}
\def\myyscale{0.6}

\def\myboundaryedge{blue}
\def\myedgesize{thick}
\tikzset{->-/.style={decoration={
			markings,
			mark=at position #1 with {\arrow{stealth}}},postaction={decorate}}}

The present section is dedicated to construct,  in the spirit of \cite{schifflerserna}, a geometric realization of the socle-projective categories for posets of type $\mathbb{A}$ using the geometric model of BGMS \cite{BGMS}. \\

According \cite{schifflerserna}, we recall the posets of type $\mathbb{A}$. There is a deep relation between Dynkin  quivers and posets of type $\mathbb{A}$ because each poset $\P$ of type $\mathbb{A}$ is the poset associated to a quiver $Q^F$ obtained from a Dynkin quiver $Q$ of type $\mathbb{A}$ by adding a set of new arrows $F$ as follows.\\

A set $F=\{\alpha_1,\dots,\alpha_t\}$ of new arrows for $Q$ is called an \textit{alien set} for $Q$ if the following conditions  hold:
\begin{itemize}
\item[(a)] for all $\alpha\in F$, there exists a sink vertex   $z$ (that is, there is no arrow starting at $z$) in $Q$ such that the vertices $s(\alpha)$  (starting of $\alpha$)  and  $t(\alpha)$ (target  of $\alpha$)  belong to the support of the indecomposable injective representation $I(z)$ at vertex $z$;
\item[(b)] for all $\alpha\in F$, $t(\alpha)$  is not a source vertex  in $Q$ (that is, there is an arrow ending at $t(\alpha)$) unless it is either vertex $1$ or vertex $n$ in the  diagram 
\begin{center}
\begin{tikzcd}[arrows=-,row sep= tiny, column sep = normal]
 \underset{_1}{\circ} \arrow[r] & \underset{_2}{\circ}\arrow[r,dotted, thick] & \underset{_{n-1}}{\circ} \arrow[r] &
  \underset{_{n}}{\circ},  
 \end{tikzcd}
\end{center}
associated to $Q$;
\item [(c)]  for all $\alpha\in F$,  the  arrow $\alpha$ is the unique path from $s(\alpha)$ to  $t(\alpha)$ in $Q^F$, where $Q^F$  is  the quiver such that $Q^F_0=Q_0$ (set of vertices) and $Q^F_1=Q_1\cup F$ (set of arrows). 
 \item[(d)] the quiver $Q^F$ is acyclic.
\end{itemize}
The arrows in an alien set for $Q$ will be called \textit{alien arrows}.

\begin{example}\label{alien set}
Let $Q$ be the quiver  $\xymatrix@=1.5em{1\ar[r]&2&3\ar[l]\ar[r]&4\ar[r]&5&6\ar[l]\ar[r]&7}$. The set $$F=\lbrace\alpha:3\rightarrow 1,\hspace{0.1cm}\beta:6\rightarrow 4\rbrace$$   is an  alien set  for $Q$. Moreover,  the quiver $Q^F$ is equal to

\begin{center}
\begin{tikzcd}[row sep= tiny, column sep = small]
1\arrow[dr]&&3\arrow[dl]\arrow[ll, "\alpha",swap]\arrow[dr]&&&&\\
&2&&4\arrow[dr]&&6\arrow[ll,"\beta",swap]\arrow[dl]\arrow[dr]&\\
&&&&5&&7\\ 
\end{tikzcd}
\end{center}
\end{example}
Recall \cite{schifflerserna} that the poset $\P_{Q^F}$ associated to $Q^{F}$, given by $x\preceq y$ in $\P_{Q^F}$ if and only if there exists a path from $x$ to $y$ in $Q^F$, is a poset of type $\mathbb{A}$. For instance, the poset associated to $Q^F$ in Example \ref{alien set} is the three-peak poset defined in Example \ref{non-positive-weight}. Furthermore, it was proved in \cite{schifflerserna}*{Proposition 3.8} that
a poset $\P$ is of type $\mathbb{A}$ if and only if there exists a Dynkin quiver $Q$ of type $\mathbb{A}$ and an  alien set $F$ for $Q$ such that $\P=\P_{Q^F}$.

\subsection{The geometric model of BGMS} \label{section3.1}

In this subsection, we recall some concepts  regarding the construction of the polygon $P(Q)$ associated to a Dynkin quiver of type $\mathbb{A}_{n+2}$, and its corresponding category of line segments introduced in \cites{BGMS}.

\subsubsection*{The polygon $P (Q)$}  In \cite{Reading} the authors defined the polygon $P(Q)$ induced by a Dynkin quiver $Q$ of type $\mathbb{A}$ as follows. Suppose that $Q$ is a quiver of type $\mathbb{A}_{n+2}$ so that  its vertices are labeled from left to right in linear order. Let $[n+1]$ be the set $\{1, 2, \dots, n+1 \}$. First, use $Q$ to divide the set $[n+1]$ into two sets, the upper-barred integers $\overline{[n+1]}$ and the lower-barred integers $\underline{[n+1]}$, as follows. Let $\overline{[n+1]}$ be the set of all vertices $i$ such that $i \rightarrow (i+1)$ is in $Q$ and let $\underline{[n+1]}$ be the set of all vertices $i$ such that $i \leftarrow (i+1)$ is in $Q$.\\

Next, associate to $Q$ an $(n + 3)$-gon $P(Q)$ with vertex labels $0, 1, 2, \dots , n+2$. Draw the vertices $0, 1, 2, \dots , n+2$ in order from left to right so that:

\begin{enumerate} 
    \item [(1)] the vertices $0$ and $n+2$ are placed on the same horizontal line $L$;
    \item [(2)] the upper-barred vertices are placed above $L$;
    \item [(3)] the lower-barred vertices are placed below $L$.
\end{enumerate}

\begin{example} 
Let $Q$ be the quiver in Example  \ref{alien set}. The following figure shows the polygon associated the quiver $Q$ where $\overline{[n+1]}$ is equal to the set $\{ 1,2,4,6 \}$ and the  $\underline{[n+1]}$ is the set $\{ 2,5\}$. 
\end{example}
\begin{figure}[ht]
\begin{center}

\begin{tikzpicture}[xscale=0.4,yscale=0.4]
\tikzstyle{every node}=[font=\footnotesize]
	\node (0) at (0,0) [fill,circle,inner sep=0.3] {};
		\node (1) at (0.8,1.4) [opacity=1,fill,circle,inner sep=0.3] {};
		\node (3) at (2.2,2) [fill,circle,inner sep=0.3] {};
		\node (4) at (3.8,2) [fill,circle,inner sep=0.3]{}; 
		\node (6) at (5.2,1.4) [fill,circle,inner sep=0.3] {};
		\node (7) at (6.4,0) [fill,circle,inner sep=0.3] {};
		\node (2) at (1.4,-1.2) [fill,circle,inner sep=0.3] {};
		\node (5) at (4.6,-1.2) [fill,circle,inner sep=0.3] {};
		\draw (0) -- (1) -- (3) -- (4) -- (6) -- (7)-- (5) -- (2)--(0);
		
		\draw (0) node[left] {$0$};
		\draw (1) node[above] {$1$};
		\draw (2) node[below] {$2$};
		\draw (3) node[above] {$3$};
		\draw (4) node[above] {$4$};
		\draw (5) node[below] {$5$};
		\draw (6) node[above] {$6$};
		\draw (7) node[right] {$7$};

\end{tikzpicture}
\end{center}
\caption{Polygon $P(Q)$ given by Example \ref{alien set}} \label{fan}

\end{figure}
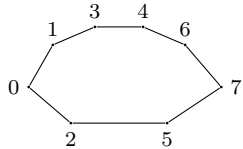

\subsubsection*{The Category of line segments} Given a polygon $P(Q)$, $\cale$ is the set of all line segments $\gamma(i,j)$ of $P (Q)$ defined as follows
$$ \cale = \{ \gamma(i,j) \mid 0 \leq i < j \leq n+2 \} $$ where $\gamma(i,j)$ denotes the line segment between the vertices $i$ and $j$. For any vertex $\ell$ of $P (Q)$, denote by $R(\ell)$ and by $R^{-1}(\ell)$, respectively, the clockwise and counterclockwise neighbor of $\ell$ on the boundary of $P(Q)$.
\begin{definition}\cite{BGMS}*{Definition 4.2} \label{definition1}
  Let $\gamma(i,j) \in \cale$. The line segment $\gamma(i,R^{-1}(j))$ is called a \textit{pivot} of $\gamma(i,j)$ if it lies in the set $\cale$, that is, if $i< R^{-1}(j)$. Similarly, the line segment $\gamma(R^{-1}(i),j)$ is called a \textit{pivot} of $\gamma(i,j)$ if it lies in the set $\cale$, that is, if $R^{-1}(i)<j$.
\end{definition}

The authors in \cite{BGMS} defined a translation quiver $(\Gamma_{P (Q)}, R)$ \cite{BGMS}. The vertices of the quiver are the line segments in $\cale$. There is an arrow $\gamma(i,j) \rightarrow \gamma(i',j')$  if and only if $\gamma(i',j')$ is a pivot of $\gamma(i,j)$, thus if and only if $i'<j'$ and $(i', j') = (R^{-1}(i), j)$ or $(i', j') = (i, R^{-1}(j))$. Finally, the translation $R$ is defined by
\begin{equation*}
   R(\gamma(i,j))= \left\{ \begin{array}{ll}
\gamma(R(i), R(j)) &   \mbox{if } R(i) < R(j), 
\\ 0 &  \mbox{otherwise.} 
\end{array}
\right.
\end{equation*}
Thus $R$ acts on $\cale$ by rotation.

\begin{definition}\cite{BGMS}*{Definition 4.4}
 Let $\mathcal{C}_{P(Q)}$ be the mesh category of the translation quiver $(\Gamma_{P (Q)}, R)$. $\mathcal{C}_{P(Q)}$ is called the \textit{category of line segments} of $P (Q)$.
\end{definition}

We recall that, for an arrow $\alpha$ in $Q$, the \textit{hook} (resp. \textit{cohook}) of $\alpha$ is given by the maximal path of $Q$ starting (resp. ending) at $x = s(\alpha)$ (resp. $y=t(\alpha)$) that does not use $\alpha$  \cite{Butler}.\\

In \cite{BGMS} the authors defined a $\Bbbk$-lineal additive functor  
\begin{equation}\label{functorF}
\mathbf{F}: \mathcal{C}_{P(Q)} \rightarrow \ind Q,    
\end{equation}
where $\ind Q$ denotes the indecomposable objets in the category $\rep Q$ of representations of $Q$. On objects this functor acts as  $$\F(\gamma(i,j))=M(i+1,j),$$
 where $M(i,j)$  denotes the indecomposable representation supported on the vertices between $i$ and $j$ for each $1 \leq i\leq j\leq n+2$, that is, $M(i,j)=(M_\ell,\varphi_\alpha)$ with
$$M_\ell=\left\{\begin{array}{ll}
\Bbbk&\textup{if $i\leq \ell\leq j$;}\\
0&\textup{otherwise;}\end{array}\right.$$
and   $\varphi_\alpha=1$, whenever $M_{s(\alpha)}$ and $M_{t(\alpha)}$ are nonzero, and $\varphi_\alpha=0$, otherwise.\\

To define $\F$ on  morphisms, it suffices to define it on the pivots introduced in Definition~\ref{definition1}. Define $\F\left(\zg(i,j)\to \zg(R^{-1}(i),j)\right) $ to  be the irreducible morphism $M(i+1,j)\to M(R^{-1}(i) +1 , j)$ given by \begin{itemize}
\item [(a)] adding the hook corresponding to the boundary edge $\zg(R^{-1}(i),i)$, if $R^{-1}(i)<i$;
\item[(b)]  removing the cohook corresponding to the boundary edge $\zg(i,R^{-1}(i))$, if $i<R^{-1}(i)$.

\end{itemize}
Similarly,  $\F\left(\zg(i,j)\to \zg(i,R^{-1}(j))\right) $ is the irreducible morphism $M(i+1,j)\to M(i+1,R^{-1}(j) )$ given by \begin{itemize}
\item [(a)] adding the hook corresponding to the boundary edge $\zg(j,R^{-1}(j))$, if $j<R^{-1}(j)$;
\item[(b)]  removing the cohook corresponding to the boundary edge $\zg(R^{-1}(j),j)$, if $R^{-1}(j)<j$.
\end{itemize}

\begin{theorem}\cite{BGMS}*{Theorem 4.6}
 The functor $\F$ is an equivalence of categories
 $$\mathcal{C}_{P(Q)} \rightarrow \ind Q$$
 In particular,
 \begin{enumerate}
 \item $\F$ induces an isomorphism of translation quivers 
 $\left(\Gamma_{P(Q)}, R \right) \rightarrow \left(\Gamma_{\rep Q}, \tau \right)$; 
\item  $\F$ induces bijections
 $$\begin{array}{rcl}
 \{ \textup{line segments in $P(Q)$}\}&\to& \ind Q;\\
\{ \textup{pivots  in $P(Q)$}\} &\to& \{\textup{irreducible morphisms in $Q$}\};
\end{array}$$
\item the rotation $R$ corresponds to the Auslander--Reiten translation $\tau$ in the following sense
$$\F\circ R= \tau \circ \F;$$
\item $\F$ is an exact functor with respect to the induced abelian structure on $\mathcal{C}_{P(Q)}$.
\end{enumerate}
\end{theorem}

\subsection{Category of sp-segments} 
Let $Q$ be a quiver of type $\mathbb{A}_{n+2}$ and let $P(Q)$ be its corresponding $(n+3)$-gon.

\begin{definition}
A line segment $\gamma(i,i+1)$ is called a \textit{suitable line segment} if it satisfies any of the following conditions. 
\begin{enumerate}
\item [(a)] $i\in \overline{[n+1]}\cup \lbrace 0 \rbrace \text{ and } i+1 \in \underline{[n+1]}\cup\lbrace n+2 \rbrace$, or
\item [(b)] $i+1 \in \overline{[n+1]}\cup \lbrace 0 \rbrace \text{ and } i \in \underline{[n+1]}\cup\lbrace n+2 \rbrace$.
\end{enumerate}
\end{definition}

\begin{example}
Given the quiver $Q$ in Example \ref{alien set}, the suitable line segments in $P(Q)$ are $\gamma(1,2)$, $\gamma(2,3)$, $\gamma(4,5)$, $\gamma(5,6)$ and $\gamma(6,7)$.

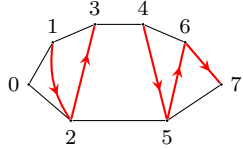
\begin{figure}[ht]
\begin{center}

\begin{tikzpicture}[xscale=0.4,yscale=0.4]
\tikzstyle{every node}=[font=\footnotesize]
	\node (0) at (0,0) [fill,circle,inner sep=0.3] {};
		\node (1) at (0.8,1.4) [opacity=1,fill,circle,inner sep=0.3] {};
		\node (3) at (2.2,2) [fill,circle,inner sep=0.3] {};
		\node (4) at (3.8,2) [fill,circle,inner sep=0.3]{}; 
		\node (6) at (5.2,1.4) [fill,circle,inner sep=0.3] {};
		\node (7) at (6.4,0) [fill,circle,inner sep=0.3] {};
		\node (2) at (1.4,-1.2) [fill,circle,inner sep=0.3] {};
		\node (5) at (4.6,-1.2) [fill,circle,inner sep=0.3] {};
		\draw (0) -- (1) -- (3) -- (4) -- (6) -- (7)-- (5) -- (2) -- (0);

		\draw (0) node[left] {$0$};
		\draw (1) node[above] {$1$};
		\draw (2) node[below] {$2$};
		\draw (3) node[above] {$3$};
		\draw (4) node[above] {$4$};
		\draw (5) node[below] {$5$};
		\draw (6) node[above] {$6$};
		\draw (7) node[right] {$7$};

\draw[->-=0.66,red,\myedgesize] (1) to [bend right=20] (2);		
\draw[->-=0.7,red,\myedgesize] (2) -- (3); 
\draw[->-=0.7, red, \myedgesize] (4) -- (5);
\draw[->-=0.7, red, \myedgesize] (5) -- (6);
\draw[->-=0.7, red, \myedgesize] (6) -- (7);
\end{tikzpicture}
\end{center}
\caption{Set of all suitable line segments in $P(Q)$.} \label{linesegments}

\end{figure}
 
\end{example}
Note that the suitable line segments correspond to simple projective and simple injective objects in $\ind Q$.
\begin{proposition} \label{propositionsegment1.1}
Let $Q$ be a quiver of type $\mathbb{A}_{n+2}$, let $P(Q)$ be its corresponding polygon and let $\F$ be the functor in  Definition \ref{functorF}. The segment $\gamma(i-1,i)$ is a suitable line segment in $P(Q)$ if and only if $\F(\gamma(i-1,i))$ is a simple projective  or a simple injective representation. 
\end{proposition}

\begin{proof}
Suppose  that, $i \in \underline{[n+1]}\cup\lbrace n+2 \rbrace$, then there are two arrows $\alpha$ and $\beta$ in $Q$ such that $\alpha$ goes from $i-1$ to $i$, and $\beta$ starts  at $i+1$ and ends at $i$, thus the vertex $i$ is a sink in $Q$. On the other hand, if $\F(\gamma(i-1,i))$ is a simple  projective representation, there are arrows from $i+1$ to $i$ and from $i-1$ to $i$, i.e., $i-1 \in \overline{[n+1]}\cup\lbrace n+2 \rbrace$ and $i \in \underline{[n+1]}\cup\lbrace n+2 \rbrace$. Other case is obtained via the same argument.
\end{proof}

\begin{definition} The \textit{bi-fan} $[\gamma(n,m)]_{bf}$ of the line segment $\gamma(n,m)$ in $P(Q)$ is given by the set
$$\{\gamma(x,y)\in P(Q) \mid x=R^{-t}(n)<n \text{  and } m<y=R^{-t'}(m)\  \exists t,t'\in \mathbb{Z}_{\geq 0} \},$$ where
for any vertex $\ell$ in $P(Q)$,  $R^{-t}(\ell)=R^{-1}(R^{-(t-1)}(\ell))$, $R^{t}(\ell)=R(R^{t-1}(\ell))$ with  $t>0$ and $R^{-0}(\ell)=R^{0}(\ell)=\ell$.

\end{definition}
\begin{proposition}
Let $\gamma(i-1,i)$ be a suitable line segment in $P(Q)$ such that $i \in \underline{[n+1]}\cup\lbrace n+2 \rbrace$. Then there is a bijective correspondence between the bi-fan of $\gamma(i-1,i)$ and the objects $X$ in $\ind Q$ such that  $\Hom_{\Bbbk}(M(i,i),X)\neq 0$. 
\end{proposition}

\begin{proof}
Let $\gamma(i-1,i)$ be a suitable line segment in $P(Q)$ such that $i \in \underline{[n+1]}\cup\lbrace n+2 \rbrace$, and let $H_i$ be the set of all of objects $X$ in $\ind Q$ such that $\Hom_{\Bbbk}(M(i,i),X)\neq 0$. We define a map $\varphi: [\gamma(i-1,i)]_{bf} \longrightarrow  H_{i}$, where $\varphi(\gamma(s,m))=\F(\gamma(s,m))=M(s+1,m)$. Firstly, we will prove that $\varphi$ is well defined, suppose that $\gamma(s,m) \in [\gamma(i-1,i)]_{bf}$, then there exist $t,t' \in \mathbb{Z}_{\geq 0}$ such that $R^{-t}(i-1)=s<i-1$ and $R^{-t'}(i)=m> i$, it turns out that

$$\xymatrix @=0.6cm{ M(i,i) \ar[d]_{f}:&  \cdots & 0  \ar[d] \ar@{-}[r] \ar@{-}[l]& \cdots & 0 \ar[r]^{0} \ar[d]_{0} \ar@{-}[l]& \Bbbk \ar[d]^{1_{\Bbbk}}& \ar[l]_{0} \ar[d]^{0} 0 \ar@{-}[r]& \cdots & 0 \ar@{-}[r] \ar[d] \ar@{-}[l]& \cdots \\
M(s+1,m):& \cdots  &\Bbbk \ar@{-}[r] \ar@{-}[l]& \cdots & \Bbbk \ar[r]_{1_{\Bbbk}} \ar@{-}[l]& \Bbbk & \ar[l]^{1_{\Bbbk}} \Bbbk \ar@{-}[r]& \cdots & \Bbbk \ar@{-}[l] \ar@{-}[r]& \cdots }
$$
i.e.,  $\Hom_{\Bbbk}(M(i,i), M(s+1, m))\simeq \Bbbk $, therefore $\varphi(\gamma(s,m))$ belongs to $H_{i}$, we also note that by definition map $\varphi$ is injective and surjective. We are done.
\end{proof}

\begin{example}\label{example bi-fan}
Let $Q$ be the quiver in Example \ref{alien set}  and  let  $\gamma(1,2)$ be a suitable line segment in $P(Q)$ with $2 \in \underline{[6]}\cup\lbrace 7 \rbrace$. For $\gamma(1,2)$, $R^{-1}(1)=0$, $R^{-1}(2)=5$, $R^{-2}(2)=7$, $R^{-3}(2)=6$, $R^{-4}(2)=4$, and $R^{-5}(2)=3$. In this case, the bi-fan of $\gamma(1,2)$ is given by  

$$\left \{ 
      \begin{matrix} 
        \gamma (1,2), \gamma(1,5),\gamma (1,7), \gamma (1, 6), \gamma (1,4), \gamma (1,3), \\
         \gamma (0,2), \gamma(0,5),\gamma (0,7),\gamma (0, 6), \gamma (0,4), \gamma (0,3)
      \end{matrix} 
   \right \}. $$
   
 In addition the bi-fan of  $\gamma(4,5)$ is the set
 
 $$\begin{Bmatrix} 
        \gamma(4,5),  \gamma(4,7), \gamma(4,6),\\
        \gamma(3,5), \gamma(3,7), \gamma(3,6),\\ \gamma(1,5),\gamma(1,7), \gamma(1,6),\\ \gamma(0,5), \gamma(0,7), \gamma(0,6),\\ \gamma(2,5), \gamma(2,7), \gamma(2,6)
      \end{Bmatrix}. $$
Finally, the bi-fan of   $\gamma(6,7)$ is the set \small{$$\lbrace \gamma(6,7), \gamma(4,7), \gamma(3,7), \gamma(1,7), \gamma(0,7), \gamma(2,7), \gamma(5,7)\rbrace.$$}

\end{example}

Recall that a poset $\P$ is called a \textit{chain}  if any two elements of $\P$ are comparable. 
Note that, the sets $C=\underline{[n+1]}\cup\lbrace n+2 \rbrace $  and $D=\overline{[n+1]}\cup \lbrace 0 \rbrace$ with the natural orders are  chains.

\begin{definition} A subchain $C'$  of $C$  such that $\vert C' \vert\geq 2$  is a \textit{principal subchain} of $C$  if the following conditions are true:
\begin{enumerate}
\item [(a)] there is exactly one suitable line segment ending at a vertex $x$ of $C'$ and $x=\min C'$;
\item[(b)] there are at most one suitable line segment starting at a vertex of $C'$. If there exists one suitable line segment starting at a vertex $x$ of $C'$ then $x=\max C'$, otherwise $\max C'=n+2$. 

\end{enumerate}
\end{definition}

\begin{definition} A subchain $D'$  of $D$  such that $\vert D' \vert\geq 2$  is a \textit{principal subchain} of $D$  if the following conditions are true:
\begin{enumerate}
\item [(a)] there is exactly one suitable line segment starting at a vertex $x$ of $D'$ and $x=\max D'$;
\item[(b)] there are at most one suitable line segment ending at a vertex of $D'$. If there exists one suitable line segment ending at a vertex $x$ of $D'$ then $x=\min D'$, otherwise $\min D'=0$. 
\end{enumerate}
\end{definition}

Denote by $\mathfrak{C}$  the set of all principal subchains of $C$, and by $\mathfrak{D}$ the set of all principal subchains of $D$.\\

Let $C'$ and $D'$ be  principal subchains in $\mathfrak{C}$ and $\mathfrak{D}$, respectively. We denote by  $R_{C'}$   the set of the  line segments  $\gamma(s,m)$ in $P(Q)$ such that $s\in C'\setminus \{\max C'\}$. Similarly,   $R^{D'}$ denotes the set of the line segments $\gamma(m,s)$ in $P(Q)$ such that $s\in D\setminus\{\min D'\}$.\\

\begin{definition} Given 
$\underline{P(Q)}=\left(\bigcup_{C'\in\mathfrak{C}} R_{C'} \right)\cup\left(\bigcup_{D'\in\mathfrak{D}}R^{D'} \right)$. The segments in $\overline{P(Q)}\setminus \underline{P(Q)}$  are called the \textit{ line $\star$-segments} in $P(Q)$ where $$\overline{P(Q)}=\bigcup_{\substack{\gamma(i-1,i) \text{ is a suitable} \\ i \in  \underline{[n+1]}\cup\lbrace n+2 \rbrace}} [\gamma(i-1,i)]_{bf}.$$ 
 \end{definition}
 
\begin{example}\label{example principal subchains}

Let $Q$ be the quiver in Example \ref{alien set}, then  $C$ is the chain $\lbrace 2 < 5< 7\rbrace$ and $D$ is the chain $\lbrace 0< 1<3<4<6\rbrace$ associated to the polygon $P(Q)$. In this case $\mathfrak{C}=\emptyset$  and the set $\mathfrak{D}$ contains exactly one principal subchain given by $\lbrace 3<4 \rbrace$, therefore  $\underline{P(Q)}=\lbrace \gamma(0,4), \gamma(1,4),\gamma(2,4), \gamma(3,4) \rbrace$.
\end{example}
 
Recall that the \textit{socle-projective representations} are the representations whose socle is a projective representation. We will see that the line $\star$-segments correspond to the socle-projective representations. 
\begin{lemma} \label{lemmasegment1.1}
Let $Q$ be a quiver of type $\mathbb{A}_{n+2}$, let $P(Q)$ be its corresponding polygon and let $\F$ be the functor in Definition \ref{functorF}. If $\gamma(r,s)$ is in $\underline{P(Q)}$, then $\F(\gamma(r,s))$ is not a socle-projective representation in $\ind Q$.
\end{lemma}

\begin{proof}
Suppose that $\gamma(s,m)$ belongs to $\underline{P(Q)}$. If $\gamma(s,m)$ is in $\bigcup_{C'\in\mathfrak{C}} R_{C'}$, then there exists  a principal subchain $C'$ in $C$ such that, if $\min C'=x$ and $\max C'=x'$,  there is an arrow from $i+1$ to $i$ in $Q$ with $x \leq i \leq x'$, as $s \in C'\setminus \{\max C'\}$, the vector spaces $V_{s+1}$ and $V_{s+2}$ in the representation $\text{soc  }M(s+1, m)$ are isomorphic to $\Bbbk$ and $0$, respectively,  i.e., $M(s+1, s+1)$ is a direct summand of the socle of $M(s+1, m)$ which is not a simple projective representation of $Q$. 
\begin{center}
\begin{tikzpicture}[y=.3cm, x=.3cm,font=\normalsize, scale=1.5]

\draw[-, >=latex,black] (1.5,0.8) -- (0.5,1.2);
\draw[-, >=latex,black] (1.5,0.8) -- (2,0.7);
\draw[-, >=latex,black] (3,0.5) -- (4,0.4);
\draw[-, >=latex,black] (5,0.5) -- (4,0.4);
\draw[-, >=latex,black] (6,0.7) -- (6.5,0.8);
\draw[-, >=latex,black] (6.5,0.8) -- (7.5,1.2);
\draw[->-=0.7,red,\myedgesize] (1.4,3.0) -- (1.5,0.8);
\draw[->-=0.7,red,\myedgesize] (6.5,0.8)--(6.6,3.0) ;
\draw[-, >=latex,black] (1.4,3.0) -- (6.6,3.0);
\draw[-, >=latex,black] (1.4,3.0) -- (0.4,2.5);
\draw[-, >=latex,black] (6.6,3.0) -- (7.6,2.5);
\node[above] at (1.5,0.2) {$_{x}$};
\node[above] at (4,-0.3) {$_{i}$};
\node[above] at (6.5,0.2) {$_{x'}$};
\node[above] at (1.4,2.9) {$_{x-1}$};
\node[above] at (6.6,2.9) {$_{x'+1}$};
\end{tikzpicture}
\end{center}

In the same way, if $\gamma(s,m) \in \bigcup_{D'\in\mathfrak{D}}R^{D'}$, there exists a principal subchain $D'$ in $D$ such that in the quiver $Q$ there is an arrow from $i$ and $i+1$ for $x' \leq i \leq x$ with $x'=\min D'$ and $x=\max D'$, for this case, as $m \in D\setminus\{\min D'\}$, the simple representation $M(m,m)$ is a direct  summand of $\text{soc }M(s+1, m)$  with $M(m,m)$ a representation which is not projective.
\begin{center}
\begin{tikzpicture}[y=.3cm, x=.3cm,font=\normalsize, scale=1.5]
\draw[-, >=latex,black] (11.5,2.6) -- (10.5,2.2);
\draw[-, >=latex,black] (11.5,2.6) -- (12,2.7);
\draw[-, >=latex,black] (13,2.9) -- (14,3);
\draw[-, >=latex,black] (15,2.9) -- (14,3);
\draw[-, >=latex,black] (16,2.7) -- (16.5,2.6);
\draw[-, >=latex,black] (16.5,2.6) -- (17.5,2.2);
\draw[->-=0.7,red,\myedgesize] (11.4,0.4) -- (11.5,2.6);
\draw[->-=0.7,red,\myedgesize] (16.5,2.6) --(16.6,0.4) ;
\draw[-, >=latex,black] (11.4,0.4) -- (16.6,0.4);
\draw[-, >=latex,black] (11.4,0.4) -- (10.4,0.9);
\draw[-, >=latex,black] (16.6,0.4) -- (17.6,0.9);
\node[above] at (11.5,2.5) {$_{x'}$};
\node[above] at (14,2.9) {$_{i}$};
\node[above] at (16.5,2.5) {$_{x}$};
\node[above] at (11.4,-0.3) {$_{x'-1}$};
\node[above] at (16.6,-0.3) {$_{x+1}$};

\end{tikzpicture}
\end{center}

Therefore $\F(\gamma(s,m))$ is not a socle-projective.
 \end{proof} 

\begin{proposition}\label{*segment}
$\gamma(s,m)$ is a line $\star$-segment in $P(Q)$ if and only if $\F(\gamma(s,m))$ is  a socle-projective representation in $\rep Q$.
\end{proposition}
\begin{proof}
Suppose that $\gamma(s,m) \in \overline{P(Q)}\setminus \underline{P(Q)}$, $s \not\in C'\setminus \{\max C'\}$ and $m \not\in  D\setminus\{\min D'\}$ for some principal subchains $C'$ and $D'$ of $C$ and $D$, respectively. Without loss of generality, suppose that $s \geq \max C'$ and $m \leq \min D'$, then we have the following cases:\\

(i) If $\max C'=n+2$, there are no line segments from $n+2$ to another vertex in $P(Q)$.\\

(ii) If $\min D'=0$, there are no line segments from another vertex to $0$.\\

(iii) If $\max C'\neq n+2$ and $\min D'\neq 0$, if $m=s+1$, for any suitable line segment $\gamma(i-1, i)$ in $P(Q)$ with $i \in \underline{[n+1]}\cup \lbrace n+2 \rbrace$, there exist $t, t' \in \mathbb{Z}_{\geq 0}$ such that $R^{-t}(i)=s$ (resp. $R^{-t}(i-1)=s$) and $R^{-t'}(i)=s+1$ (resp. $R^{-t'}(i-1)=s+1$) with $i-1 < s<s+1$ (resp. $s<s+1<i$) (see following left (resp. right) diagram)

\begin{center}
\begin{tikzpicture}[y=.3cm, x=.3cm,font=\normalsize, scale=1.5]

\draw[->-=0.7,red,\myedgesize] (1.0,3.0) -- (1.1,0.8);
\draw[-, >=latex,black] (0,2.6) -- (1.0,3.0);
\draw[-, >=latex,black] (2,3.2) -- (1.0,3.0);
\draw[-, >=latex,black] (0.1,1.2) -- (1.1,0.8);
\draw[-, >=latex,black] (2.1,0.6) -- (1.1,0.8);
\draw[->-=0.7,red,\myedgesize] (4.0,0.8) --(4.1,3.0);
\draw[-, >=latex,black] (5.1,2.6) -- (4.1,3.0);
\draw[-, >=latex,black] (3.1,3.2) -- (4.1,3.0);
\draw[-, >=latex,black] (5.0,1.2) -- (4.0,0.8);
\draw[-, >=latex,black] (3.0,0.6) -- (4.0,0.8);

\node[above] at (1.0,2.9) {$_{i-1}$};
\node[above] at (1.1,0.2) {$_{i}$};
\node[above] at (4.1,2.9) {$_{s+1}$};
\node[above] at (4.0,0.2) {$_{s}$};

\draw[->-=0.7,red,\myedgesize] (9.0,0.8)--(9.1,3.0);
\draw[-, >=latex,black] (8.1,2.6) -- (9.1,3.0);
\draw[-, >=latex,black] (10.1,3.2) -- (9.1,3.0);
\draw[-, >=latex,black] (9.0,0.8)--(8.0,1.2);
\draw[-, >=latex,black] (10.0,0.6) -- (9.0,0.8);
\draw[->-=0.7,red,\myedgesize] (12,3.0)--(12.1,0.8);
\draw[-, >=latex,black] (13.0,2.6) -- (12.0,3.0);
\draw[-, >=latex,black] (11.0,3.2) -- (12.0,3.0);
\draw[-, >=latex,black] (13.1,1.2) -- (12.1,0.8);
\draw[-, >=latex,black] (11.1,0.6) -- (12.1,0.8);

\node[above] at (9.1,2.9) {$_{s+1}$};
\node[above] at (9.0,0.2) {$_{s}$};
\node[above] at (12.0,2.9) {$_{i-1}$};
\node[above] at (12.1,0.2) {$_{i}$};

\end{tikzpicture}
\end{center}

therefore $\gamma(s,m) \not\in \overline{P(Q)}$.\\

If $m> s+1$, there are suitable line segments $\gamma(s_{i_{p}}^{2}, s_{i_{p}}^{1})$ and $\gamma(s_{j_{q}}^{1}$, $s_{j_{q}}^{2})$ with $s_{i_{p}}^{1} \in \overline{[n+1]}\cup \lbrace 0 \rbrace$, $s_{j_{q}}^{2} \in \underline{[n+1]}\cup\lbrace n+2 \rbrace$, $1 \leq p \leq h$ and  $0\leq q \leq h$ such that $\gamma(s,m)$ belongs to the bi-fans $[\gamma(s^2_{i_{r}}, s^{1}_{i_{r}})]_{bf}$ for $1 \leq r \leq h-1$, moreover, $s^{1}_{j_{0}} \leq s \leq s^{2}_{i_{1}}$ and $s^{1}_{i_{h-1}}\leq m \leq s^{2}_{j_{h}}$ (see Figure \ref{fig:gabo}).

\tikzset{every picture/.style={line width=0.2pt}} 

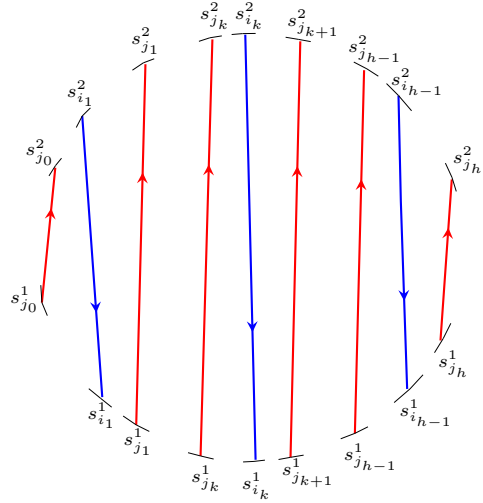
\begin{figure}[ht]
    \centering

\begin{tikzpicture}[x=0.75pt,y=0.75pt,yscale=-0.55,xscale=0.55]
\draw [->-=0.7,red,\myedgesize]   (124,271.5) -- (136.3,147.74) ;
\draw [->-=0.7,blue,\myedgesize]   (161,100.25) -- (179.36,358.26) ;
\draw [->-=0.7,red,\myedgesize]   (210.5,383.75) -- (218.95,52.75) ;
\draw [->-=0.7,blue,\myedgesize]   (310.5,25.66) -- (319.95,415.75) ;
\draw [->-=0.7,red,\myedgesize]   (269.5,412.25) -- (280.35,30) ;
\draw [->-=0.7,red,\myedgesize]   (352,413.5) -- (360.95,31.75) ;
\draw [->-=0.7,red,\myedgesize]   (489.5,306.25) -- (499.86,158.5) ;
\draw [->-=0.7,blue,\myedgesize]   (451,81.5) -- (453.14,178.52) -- (458.93,350.75) ;
\draw [->-=0.7,red,\myedgesize]   (411,391.5) -- (419.45,58.25) ;
\draw    (298,25.5) -- (308,24.66) ;
\draw    (289,26.75) -- (278,28.5) ;
\draw    (268,31.75) -- (278,28.5) ;
\draw    (217,51.5) -- (227,47.75) ;
\draw    (207,58.25) -- (217,51.5) ;
\draw    (161,100.25) -- (168.5,92.84) ;
\draw    (155.5,110.84) -- (161,100.25) ;
\draw    (136.5,145.75) -- (142,139) ;
\draw    (130.5,155.25) -- (136.5,145.75) ;
\draw    (123,255.5) -- (124,271.5) ;
\draw    (124,271.5) -- (129,283.5) ;
\draw    (167,350.5) -- (178,359.5) ;
\draw    (178,359.5) -- (188,367.5) ;
\draw    (198.5,375.59) -- (210.5,383.75) ;
\draw    (210.5,383.75) -- (222.5,389.75) ;
\draw    (258,409.25) -- (269.5,412.25) ;
\draw    (308,24.66) -- (320,24.5) ;
\draw    (269.5,412.25) -- (281.5,415.09) ;
\draw    (308.5,417.75) -- (318.5,417.25) ;
\draw    (341,415.5) -- (352,413.5) ;
\draw    (352,413.5) -- (364,411.5) ;
\draw    (398,396.5) -- (411,391.5) ;
\draw    (411,391.5) -- (423,385.5) ;
\draw    (447,363.25) -- (461.5,350.75) ;
\draw    (461.5,350.75) -- (473.5,337.75) ;
\draw    (492,304.5) -- (499,290.5) ;
\draw    (504,169.5) -- (500,156.5) ;
\draw    (494,143.5) -- (500,156.5) ;
\draw    (451,81.5) -- (463,95.5) ;
\draw    (440,70.5) -- (451,81.5) ;
\draw    (422,57) -- (432.5,63.75) ;
\draw    (410,50.75) -- (422,57) ;
\draw    (347.5,28) -- (359.5,30) ;
\draw    (359.5,30) -- (370.5,32) ;
\draw    (328.5,416.09) -- (318.5,417.25) ;
\draw    (484,317.5) -- (492,304.5) ;

\draw (94,255.4) node [anchor=north west][inner sep=0.3pt]  [font=\scriptsize]  {$s_{j_0}^{1}$};
\draw (146,66.4) node [anchor=north west][inner sep=0.3pt]  [font=\scriptsize]  {$s_{i_1}^{2}$};
\draw (165,362.4) node [anchor=north west][inner sep=0.3pt]  [font=\scriptsize]  {$s_{i_1}^{1}$};
\draw (198,386.4) node [anchor=north west][inner sep=0.3pt]  [font=\scriptsize]  {$s_{j_1}^{1}$};
\draw (205,18.4) node [anchor=north west][inner sep=0.3pt]  [font=\scriptsize]  {$s_{j_1}^{2}$};
\draw (260,419.65) node [anchor=north west][inner sep=0.3pt]  [font=\scriptsize]  {$s_{j_k}^{1}$};
\draw (268,-4.9) node [anchor=north west][inner sep=0.3pt]  [font=\scriptsize]  {$s_{j_k}^{2}$};
\draw (300,-5) node [anchor=north west][inner sep=0.3pt]  [font=\scriptsize]  {$s_{i_k}^{2}$};
\draw (306.5,426.15) node [anchor=north west][inner sep=0.3pt]  [font=\scriptsize]  {$s_{i_k}^{1}$};
\draw (343,413.9) node [anchor=north west][inner sep=0.3pt]  [font=\scriptsize]  {$s_{j_{k+1}}^{1}$};
\draw (348,1.4) node [anchor=north west][inner sep=0.3pt]  [font=\scriptsize]  {$s_{j_{k+1}}^{2}$};
\draw (486.5,314.65) node [anchor=north west][inner sep=0.3pt]  [font=\scriptsize]  {$s_{j_h}^{1}$};
\draw (443,55.4) node [anchor=north west][inner sep=0.3pt]  [font=\scriptsize]  {$s_{i_{h-1}}^{2}$};
\draw (451,359.4) node [anchor=north west][inner sep=0.3pt]  [font=\scriptsize]  {$s_{i_{h-1}}^{1}$};
\draw (404,25.4) node [anchor=north west][inner sep=0.3pt]  [font=\scriptsize]  {$s_{j_{h-1}}^{2}$};
\draw (400,399.9) node [anchor=north west][inner sep=0.3pt]  [font=\scriptsize]  {$s_{j_{h-1}}^{1}$};
\draw (108,120.4) node [anchor=north west][inner sep=0.3pt]  [font=\scriptsize]  {$s_{j_{0}}^{2}$};
\draw (500,127.4) node [anchor=north west][inner sep=0.3pt]  [font=\scriptsize]  {$s_{j_{h}}^{2}$};

\end{tikzpicture}
\caption{Configuration of suitable line segments.}
    \label{fig:gabo}
\end{figure}
The socle of the representation  $M(s+1, m)$ is given by  $\bigoplus_{r=1}^{h-1} M(s^1_{i_{r}},s^1_{i_{r}})$. Proposition \ref{propositionsegment1.1} implies that $M(s^1_{i_{r}},s^1_{i_{r}})$ is a simple projective, therefore  $\F(\gamma(s, m))$ is a socle-projective representation.\\

In the other direction, suppose that $\F(\gamma(s,m))$ is a socle-projective representation, that is, $\text{soc }M(s+1, m)$ is a direct sum of the simple projective representations $M(j_{i},j_{i})$ with  $j_{i}$ a sink in $Q$ and $0 \leq i \leq r$ for some $r \in \mathbb{Z}_{\geq 0}$. By Lemma \ref{lemmasegment1.1}, $\gamma(s,m)$ does not belong to  $\underline{P(Q)}$, then if there is no a suitable line segment to left (or right) of $\gamma(j_{0}-1,j_{0})$ (or $\gamma(j_{r}-1, j_{r})$) in $P(Q)$, we have that $0 \leq s \leq j_{0}-1$  (or $j_{r} \leq m \leq n+2$), in the case that there exists a closest suitable line segment $\gamma(r_{1},r_{2})$ to left (or right) of $\gamma(j_{0}-1,j_{0})$ (or $\gamma(j_{r}-1, j_{r})$) with $r_{2} \in \overline{[n+1]} \cup \lbrace 0 \rbrace$, in the same way $r_{1} \leq s \leq j_{0}-1$ (or $j_{k} \leq m \leq r_{2}$). For any of the two possibilities to left (or right) of $\gamma(j_{0}-1,j_{0})$ (or $\gamma(j_{r}-1, j_{r})$) there exist $t, t' \in \mathbb{Z}_{\geq 0}$ such that $R^{-t}(j_{0}-1)=s$ and $R^{-t'}(j_{r})=m$, thus $\gamma(s,m)$ is a line segment in $[\gamma(j_{i}-1, j_{i})]_{bf}$ for all $j_{i}$. We conclude that $\gamma(s,m)$ is a line $\star$-segment in $P(Q)$.
\end{proof} 


\begin{definition}
Let $\alpha: j_2 \longrightarrow j_1$ be an alien arrow in $Q$ defined in the support of the injective indecomposable $I(i)$. A line segment $\gamma(s,m)$ is an \textit{frozen line segment } in $P(Q)$ associated to $\alpha$ if one of the following conditions hold.

\begin{enumerate}
\item [(a)] If $j_{1}<j_{2}$, then $m\geq j_{2}$ and $j_{1} \leq s <i$,
\item [(b)] if $j_{1}>j_{2}$, then  $s<j_{2}$  and $i \leq m < j_{1}$,
\end{enumerate}
In this case, we say that $\gamma(n,m)$ \textit{is frozen by} $\alpha$.
\end{definition}
Note that $i$ is the end-point  of a suitable line segment $\gamma(i-1,i)$ between $j_{1}$ and $j_{2}$ with $i \in \underline{[n+1]}\cup\lbrace n+2 \rbrace$ .\\

The set of all the frozen line segments in $P(Q)$ associated to $\alpha$ is denoted by $[\alpha]$. We denote by $P(Q^F)$ the set of all the line segments in $P(Q)$ which are not frozen by an alien arrow $\alpha\in F$. In other words,

$$P(Q^F)=\overline{P(Q)}\setminus \bigcup_{\alpha\in F}[\alpha]$$ 

The line segments in $P(Q^F)$ are called \textit{non-frozen by F} line segments in $P(Q)$.
\begin{example}\label{frozen segments}
Let $F$ be the alien set of $Q$ in Example \ref{alien set}. Given that, $\alpha$ is an alien arrow from 3 to 1 in $Q$ and there is a suitable line segment $\gamma(1,2)$ in $P(Q)$ with $2 \in \underline{[6]}\cup\lbrace 7 \rbrace$, then the line segments $\gamma(1,3), \gamma(1,4), \gamma(1,5), \gamma(1,6)$ and  $\gamma(1,7)$ are frozen line segments in $Q(P)$ associated to $\alpha$. In the same way, the set  of all frozen line segments in $P(Q)$ associated to alien arrow $\beta$ is given by  $[\beta]=\lbrace \gamma(4,5), \gamma(4,6)\rbrace$.
\end{example}

\begin{definition}
A line $\star$-segment $\gamma(n,m)$  in $ P(Q^F)$ is called  a \textit{line sp-segment}. We denote by $P_{sp}(Q^F)$ the set of all the line sp-segments in $P(Q)$. 
\end{definition}

\begin{example} \label{example sp-segments}Let $Q^{F}$ be a quiver in Example \ref{alien set}.  The set of line sp-segments is given by 

$$P_{sp}(Q^{F})=[\gamma(1,2)]_{bf} \cup [\gamma(4,5)]_{bf}\cup[\gamma(6,7)]_{bf}\setminus \left( \underline{P(Q)} \cup [\alpha]\cup [\beta]\right),$$ where 
$[\gamma(1,2)]_{bf}$ , $[\gamma(4,5)]_{bf}$, and $[\gamma(6,7)]_{bf}$  are as in Example \ref{example bi-fan}, $\underline{P(Q)}$ is as in Example \ref{example principal subchains},   whereas  $[\alpha]$ and $[\beta]$ are as in Example \ref{frozen segments}.
\end{example}

\textit{Category of sp-segments: }Let $\P$ be a poset of type $\mathbb{A}$ associated to a quiver $Q^F$, where $Q$ is a Dynkin quiver of type $\mathbb{A}$ and $F$ is an alien set for $Q$. We denote  by $\mathcal{C}_{P_{\text{sp}}(Q^F)}$  the full subcategory of $\mathcal{C}_{P(Q)}$ generated by all of the line sp-segments in $P(Q)$.  \\

Since the  irreducible morphisms in  $\mathcal{C}_{P_{\text{sp}}(Q^F)}$  cannot be factorized through  sp-segments, we introduce the notion of a  \textit{sp-pivot} from  $\gamma(s,m)\in P_{\text{sp}}(Q^F)$ to  $\gamma(s',m')\in P_{\text{sp}}(Q^F)$, that is, a composition of pivots introduced in \cite{BGMS} of the form  
$\rho:\gamma(s,m)\rightarrow \gamma(s, R^{-1}(m))\rightarrow \dots \rightarrow \gamma(s,R^{-(t'-1)}(m)) \rightarrow \gamma(s,R^{-t'}(m))=\gamma(s',m') (\text{ or }\rho:\gamma(s,m)\rightarrow \gamma(R^{-1}(s),m)\rightarrow \dots \rightarrow \gamma(R^{-(t-1)}(s),m) \rightarrow \gamma(R^{-t}(s),m)=\gamma(s',m'))$  with $s<m'$ (or $s' <m$) such that $\gamma(s,R^{-r}(m))$  (or $(\gamma(R^{-r})(s),m))$) are not sp-segments  in $P(Q)$ for $1\leq r \leq t'-1$ (or $1 \leq r \leq t-1$) and some $t,t' \in \mathbb{Z}_{\geq 0}$. Note that the irreducible morphisms in  $\mathcal{C}_{P_{\text{sp}}(Q^F)}$ are precisely the  sp-pivots between sp-segments.\\

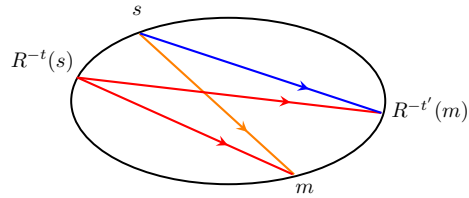
\begin{figure}[ht]
\begin{center}
\tikzset{every picture/.style={line width=0.75pt}} 

\begin{tikzpicture}[x=0.75pt,y=0.75pt,yscale=-0.8,xscale=0.8]

\draw   (209,496.16) .. controls (209,467.16) and (253.32,443.66) .. (308,443.66) .. controls (362.68,443.66) and (407,467.16) .. (407,496.16) .. controls (407,525.15) and (362.68,548.66) .. (308,548.66) .. controls (253.32,548.66) and (209,525.15) .. (209,496.16) -- cycle ;
\draw [->-=0.7,orange,\myedgesize ]   (252,453.66) -- (349.52,542.44) ;
\draw [->-=0.7,red,\myedgesize ]   (212.5,481.16) -- (348.68,542.44) ;
\draw [->-=0.7,red,\myedgesize]   (214,481.51) -- (404.51,503.67) ;
\draw [->-=0.7,blue,\myedgesize ]   (252.5,453.51) -- (404.59,503.67) ;
\draw (246.5,435.41) node [anchor=north west][inner sep=0.3pt]  [xscale=0.8,yscale=0.8]  {$s$};
\draw (349.5,547.41) node [anchor=north west][inner sep=0.3pt]  [xscale=0.8,yscale=0.8]  {$m$};
\draw (170,462.91) node [anchor=north west][inner sep=0.3pt]  [xscale=0.8,yscale=0.8]  {$R^{-t}( s)$};
\draw (410,492.91) node [anchor=north west][inner sep=0.3pt]  [xscale=0.8,yscale=0.8]  {$R^{-t'}( m)$};

\end{tikzpicture}
\end{center}
\caption{Diagram of the mesh relations in $\mathcal{C}_{P_{\text{sp}}(Q^F)}$} \label{figuremeshrelation}
\end{figure}

The mesh relations in $\mathcal{C}_{P_{\text{sp}}(Q^F)}$ is defined in the following way. We suppose that $\gamma(s,m)$ and $\gamma(s',m')$ are $sp$-segments and  that there are two compositions $\gamma(s,m) \xrightarrow{\rho_{1}} \gamma(i, j)  \xrightarrow{\rho_{2}} \gamma(s',m')$ and  $\gamma(s,m) \xrightarrow{\rho_{3}} \gamma(i', j')  \xrightarrow{\rho_{4}} \gamma(s',m')$ of $sp$-pivots with $\gamma(i,j)\neq \gamma(i',j')$. Then   $$ \gamma(s,m) \xrightarrow{\rho_{1}} \gamma(i, j)  \xrightarrow{\rho_{2}} \gamma(s',m') = \gamma(s,m) \xrightarrow{\rho_{3}} \gamma(i', j')  \xrightarrow{\rho_{4}} \gamma(s',m') $$

if and only if $m=j$, $s=i'$, $R^{-t}(s)=i=s'$ and $R^{-t}(m)=j'=m'$ for some $t, t' \in \mathbb{Z}_{\geq 0}$ (see Figure \ref{figuremeshrelation}).\\

\subsection{The categorical equivalence} 
Actually, we are going to show that the category of sp-segments is equivalent to the category of peak $\P$-spaces, where $\P$ is the poset asocciated to the quiver $Q^F$. However, for an easier explanation, we introduce a presentation of the category $\P$-spr using quiver representations. First, we recall that the \textit{incidence algebra} $\Bbbk\P$ of $\P$ can be described as a bound quiver algebra $\Bbbk Q/I$ induced by the \textit{Hasse quiver} $Q$ of $\P$ whose vertices are the points of $\P$ and there is an arrow $\alpha:x\rightarrow y$ for each pair $x,y\in\P$ such that $y$ covers $x$.  The ideal $I$ is generated by all the commutativity relations $\gamma-\gamma'$  with  $\gamma$ and $\gamma'$ parallel paths in $Q$. In this case, the category $\md(\Bbbk\P)$ of the finitely generated $\Bbbk\P$-modules is identified with the well known category $\rep (Q,I)$ of representations  of the bound quiver $(Q,I)$. More precisely, we identify the category $\P$-spr with the subcategory  $\md_{sp}(\Bbbk \P)$ of $\md(\Bbbk\P)$ of  socle-projective modules (see \cites{simson3,simson4,simson5}). Thus, an object $M$ in $\md_{sp}(\Bbbk\P)$ can be  explicitly described  as a collection $M=(M_x,{_y}h_x)_{x,y\in\P}$  of finite-dimensional $\Bbbk$-vector spaces $M_x$, one for each point $x\in\P$, and a collection of $\Bbbk$-linear maps ${_y}h_x: M_x \to M_y$, one for each relation $x\preceq y$ in $\P$, such that

\begin{itemize}
\item[(a)]${_x}h_x$ is the identity of $M_x$  for all $x\in\P$ and ${_w}h_y\cdot {_y}h_x= {_w}h_x$ for all $x\preceq y\preceq w$ in $\P$.
\item[(b)] For all $x\in\P\setminus\max\P$, the $\Bbbk$-subpace $$I_x=\stackbin[\begin{smallmatrix}  x\prec z\in\max\P \end{smallmatrix}]{}{\bigcap}\ker {_z}h_x$$ of $M_x$ is the zero subspace.
\end{itemize}

Note that it is enough to define the linear maps ${_y}h_x$ when $y$ covers $x$, that is,  one for each arrow in the Hasse quiver of $\P$ because if $x\prec y$ but $y$ does not cover $x$ then for any chain $x=x_0\prec x_1\prec\dots \prec x_{l}=y $ in $\P$ such that $x_{i+1}$ covers $x_i$ we have that ${_y}h_x={_y}h_{x_{l-1}}\cdots {_{x_{1}}}h_x$. The condition (a) implies that it is well defined.\\

Let $M=(M_x,{_y}h_x)_{x,y\in\P}$ and $N=(N_x,{_y}h'_x)_{x,y\in\P}$   be two objects in $\md(\Bbbk\P)$. Recall that a \textit{morphism}   $f:M\to N$ in $\md(\Bbbk\P)$ is a collection $f=(f_x)_{x\in\P}$ of linear maps
$f_x:M_x\to N_x$ 
such that for each relation $x\preceq y$ in $\P$, $f_y\circ{_y}h_x={_y}h'_x\circ f_x$.\\

Let $\P$ be a poset of type $\mathbb{A}$ associated to a quiver $Q^{F}$, and let $P(Q)$ be the corresponding polygon to $Q$, where $Q$ is a Dynkin quiver of type $\mathbb{A}$ labeled from left to right and $F$ an alien set for $Q$. 
We define a $\Bbbk$-linear additive functor
$$\Omegan: \mathcal{C}_{P_{\text{sp}}(Q^F)} \rightarrow \md_{sp}(\Bbbk \P)$$
from the category of sp-segments to the category of socle-projective $\Bbbk \P$ modules  such that for any sp-segment $\gamma(s,m)$ in $\mathcal{C}_{P_{\text{sp}}(Q^F)}$,  $\Omegan(\gamma(s,m))=M^{\gamma(s,m)}=(M_x^{\gamma(s,m)},{_y}h^{\gamma(s,m)}_{x})$ where $M^{\gamma(s,m)}$ is  defined by the following identity:
    $$M_x^{\gamma(s,m)}= \left\{ \begin{array}{ll}
\Bbbk &   \mbox{if } s+1 \leq x \leq m, 
\\ 0 &  \mbox{otherwise.} 
\end{array}
\right.$$ and if $x\preceq y\in \P$  then  $${_y}h_{x}^{\gamma(s,m)}= \left\{ \begin{array}{lc}
1_\Bbbk &   \mbox{if }   x,y  \in \lbrace s+1, \dots, m \rbrace ,\\
 0 &  \mbox{otherwise.} 
\end{array}
\right.$$ 

For any sp-pivot $\rho:\gamma(s,m)\longrightarrow \gamma(s', m')$ in $\mathcal{C}_{P_{\text{sp}}(Q^F)}$ where $s'=R^{-t}(s)<m'=m$ or $m'=R^{-t'}(m)>s'=s$,  the morphism
$\Omegan({\rho})=(\Omegan(\rho)_x)_{x\in \P}$  $$(\Omegan(\rho)_x)_{x\in \P}:(M_x^{\gamma(s,m)}, {_y}h^{\gamma(s,m)}_x) \to  (M_x^{\gamma(s',m')}, {_y}h^{\gamma(s',m')}_x)$$ satisfies that
$$\Omegan(\rho)_x=
\begin{cases}
1_\Bbbk,  &\text{if }M^{\gamma(s,m)}_x=M^{\gamma(s',m')}_x=\Bbbk, \\
0,            &\text{otherwise.}
\end{cases}$$

\begin{theorem}\label{categoricalequivalence}
The functor $\Omegan$ is an equivalence of categories between the category of sp-segments and the category of socle-projective $\Bbbk \P$-modules.
\end{theorem}

\begin{proof}
Let $\P$ be a poset of type $\mathbb{A}$, let $\P_{Q}$ be the same poset $\P$ without the alien set $F$ and suppose that \textbf{n+1} is the chain $1 < 2 < \dots < n <n+1$ (note that, \textbf{n+1} is a linearly ordered set given by the vertices labelled in  the quiver $Q$). First, we will prove that the functor is well defined, for this we have to check that $\Omegan(\gamma(s,m))$ satisfies the conditions (a) and (b) above. To prove (a), we suppose that $x \prec y \prec w$ in $\P$ such that $s+1 \leq x < w \leq m$, $y < s+1$ and $y> m$ in \textbf{n+1}, then we have the following cases.\\

(i) If $x \prec y \prec w$ in $\P_{Q}$, then $s+1 \leq x<y<w \leq m$, but this is a contradiction.\\

(ii) If $x\prec y$ in $\P_{Q}$ and $y \nprec w$ in $\P_{Q}$, there exists an alien arrow $\alpha$ from a vertex $y'$ to a vertex $w'$ in $Q$ defined in the support of a injective indecomposable $M(z,z)$, such that $y \preceq y' \prec z$ and $w' \preceq w \prec z$ in $\P_{Q}$. In the case when $m < w'$ in \textbf{n+1}, $\gamma(s,m)$ is frozen by $[\alpha]$, which is a contradiction. If $m \geq w'$ in \textbf{n+1}, $s+1 \leq x < y \leq y' < z < w \leq w' \leq m$, again a contradiction.\\

(iii) If  $y \prec w$ in $\P_{Q}$ and $x \nprec y$ in $\P_{Q}$, there is an alien arrow $\alpha: x'\longrightarrow y'$ associated to a vertex sink $z$, with $x \preceq x' \prec z$ and $y' \preceq y \prec z$ in $\P_{Q}$. When $m < y'$ in \textbf{n+1}, $\gamma(s,m) \in [\alpha]$ but this is a contradiction, if $m \geq y'$ in \textbf{n+1}, we have that $s+1 \leq x \leq x' < z <w <y \leq y' \leq m$ in \textbf{n+1}, which is a contradiction.\\

(iv) If $x \nprec y$ in $\P_{Q}$ and $y \nprec w$ in $\P_{Q}$, there are two aliens arrows $\alpha:x' \longrightarrow y'$ and $\beta:y''\longrightarrow w'$ defined in the support of a injective indecomposable $M(z,z)$, such that $x \preceq x' \prec w' \preceq w \prec z $ in $\P_{Q}$ and $y' \preceq y \preceq y'' \prec z$ in $\P_{Q}$. For the case when $m < y'$ in \textbf{n+1}, it implies that  $\gamma(s,m)$ belongs to the set $[\alpha]$, again a contradiction. When $m \geq y'$, then  $s+1 \leq x \leq x' <w' \leq w < z <y'' \leq y \leq y' \leq m$ in \textbf{n+1}, but this is a contradiction.\\

In this way, if $x \prec y \prec w$ in $\P$ with $s+1 \leq x < w \leq m$ then $s+1 \leq y \leq  m$, therefore ${_w}h_{x}^{\gamma(s,m)}={_w}h_{y}^{\gamma(s,m)}= {_y}h_{x}^{\gamma(s,m)}= 1_{\Bbbk}$.\\

To check (b), suppose that $x \leq s+1$ or $x > m$ in \textbf{n+1} for any $x \in \P\setminus\max\P$, it hold that $\ker {_z}h_x=0$ for each $z\in\max\P$ with $x\prec z$. In the same way, if $s+1 \leq x \leq m$, for any $z\in\max\P$ with $x\prec z$ we have that $M^{\gamma(s,m)}_{x}\simeq M^{\gamma(s,m)}_{z}$ and ${_z}h_{x}^{\gamma(s,m)}=0$, thus  $\bigcup_{z\in\max\P} \ker {_z}h_{x}^{\gamma(s,m)}=0$. Similar arguments can be used to obtain the proof in the case $s+1 \leq w < x \leq m$, therefore  $\Omegan(\gamma(s,m))$ is an object in $\md_{sp}(\Bbbk \P)$.\\

In order to prove that $\Omegan({\rho})$ is a morphism in $\md_{sp}(\Bbbk \P)$, let $\rho:\gamma(s,m)\rightarrow \gamma(s',m')$ be a sp-pivot in $\mathcal{C}_{P_{\text{sp}}(Q^F)}$ and to consider that $x\prec y$ in $\P$ such that $y$ covers $x$, then, we have to check the commutativity of the following diagram.

$$ \xymatrix@=1.3cm{ M_x^{\gamma(s,m)} \ar[r]^{{_y}h_{x}^{\gamma(s,m)}} \ar[d]_{\Omegan({\rho})_{x}}& M_y^{\gamma(s,m)} \ar[d]^{\Omegan({\rho})_y} \\
M_x^{\gamma(s,m)} \ar[r]^{ {_y}h_{x}^{\gamma(s',m')}}&   M_y^{\gamma(s,m)}} $$

For this, we suppose the following critical possibilities where the diagram does not commute with $x < y $ in \textbf{n+1}.\\

(v) If $ M_x^{\gamma(s,m)}= \Bbbk$, $ M_y^{\gamma(s,m)}=0$,  $M_x^{\gamma(s',m')}=\Bbbk$, and $ M_x^{\gamma(s',m')}=\Bbbk$, then $s+1 \leq x \leq m$ and $y >m $ in \textbf{n+1}, if $x \prec y$ in $\P_{Q}$ with $x$ a vertex different to a source in $Q$, $x=m$ and $y=m+1$ in \textbf{n+1}, then there is a maximal point $z$ in $\P_{Q}$ with $y \prec z$ such that $\ker {_z}h_{x}^{\gamma(s,m)}=0$, in this case  $\Omegan(\gamma(s,m))$ does not belong to $\md_{sp}(\Bbbk \P)$, which is a contradiction. If $x\nprec y$ in $\P_{Q}$, there exists an alien arrow $\alpha: x \rightarrow y$ defined in the support of a injective indecomposable $M(z,z)$ with $z$ a maximal point in $\P$, given that $y > m$ then $\gamma(s,m) \in [\alpha]$, again a contradiction. If $x$ is a source in $Q$, $s'+1 \leq x,y \leq m'$ and $s=s'$, if $x \prec y$ in $\P_{Q}$, then $m' \geq y$ in \textbf{n+1}, if $m' \in \overline{[n+1]}$, there is no a $t\geq 0$ such that $R^{-t}(m)=m'$  but this is a contradiction. The same situation occurs when $m \in \underline{[n+1]}$.\\

(vi)  Same arguments are used for the case $s+1 \leq x,y \leq m$, $s'+1 \leq y \leq m'$ and $x > m'$ or $m < s'+1$, i.e.,  the case when $ M_x^{\gamma(s,m)}=\Bbbk$,  $M_y^{\gamma(s,m)}=\Bbbk$, $ M_x^{\gamma(s',m')}=0$ and $ M_x^{\gamma(s',m')}=\Bbbk$ does not occur.\\

Note that, in the other situations that satisfy the hypothesis, the diagram commutes. Same arguments are used for the case when $y < x $ in \textbf{n+1}.\\

Finally,  we  will check the mesh relations. Suppose that $\gamma(s,m)$, $\gamma(s',m')$, $\gamma(i, j)$ and $\gamma(s',m')$  are sp-segments and $\gamma(s,m) \xrightarrow{\rho_{1}} \gamma(i, j)  \xrightarrow{\rho_{2}} \gamma(s',m') = \gamma(s,m) \xrightarrow{\rho_{3}} \gamma(i', j')  \xrightarrow{\rho_{4}} \gamma(s',m')$ for $\gamma(i, j)\neq \gamma(i', j')$ and $ \rho_{i}$ a sp-pivot in $\mathcal{C}_{P_{\text{sp}}(Q^F)}$ with $i=1,2,3,4$ (see Figure \ref{figuremeshrelation}). Assume that $x \in \P$,  and we  check if the following diagram

$$ \xymatrix@=1.3cm{ M_x^{\gamma(s,m)} \ar[r]^{\Omegan(\rho_{1})_{x}} \ar[d]_{\Omegan(\rho_{3})_{x}}& M_x^{\gamma(i,j)} \ar[d]^{\Omegan(\rho_{2})_x} \\
M_x^{\gamma(i',j')} \ar[r]^{\Omegan(\rho_{4})_{x}}&   M_x^{\gamma(s',m')}} $$

commutes. Again, we  take  the two critical possibilities where the diagram does not commute, given that $i<m$ and $j'> s$, if $s+1 \leq x \leq m$, $i'+1\leq x \leq j'$, $s'+1 \leq x \leq m'$ and $x< i+1$ or $x> j$ in \textbf{n+1}, then  $$s+1 < x \leq s' < s'+1 \leq s  \text{ in  \textbf{n+1}}$$
 which is a contradiction. The other case where $s+1 \leq x \leq m $, $i'+1 \leq x \leq j'$, $s'+1 \leq x \leq m'$ and $x<i+1$ or $x>j$ is proved in a similar way, consequently these possibilities do not occur. Therefore the functor $\Omegan$ is well defined.\\
 
In order to prove that $\Omegan$ is is faithful, we need to show that the image of a nonzero morphism between two $sp$-segments is a nonzero morphism in the category $\md_{sp}(\Bbbk \P)$. Suppose that $\rho \in \Hom_{\mathcal{C}_{P_{\text{sp}}(Q^F)}}(\gamma(s,m), \gamma(s',m'))$ is a nonzero morphism in $\mathcal{C}_{P_{\text{sp}}(Q^F)}$, in particular, $\rho$ is a nonzero morphism in $\mathcal{C}_{P(Q)}$. Then there exists a vertex $x$ in $Q$ such that $s+1 \leq x \leq m$ and $s'+1 \leq x \leq m'$ in \textbf{n+1}, by definition of $\Omegan$, $M_{x}^{\gamma(s,m)}=\Bbbk=M_{x}^{\gamma(s',m')}$ and $\Omegan(\rho)_{x}=1_{\Bbbk}$.\\

To prove that $\Omegan$ is full, we suppose that $\Omegan(\rho): \Omegan(\gamma(s,m))\rightarrow \Omegan(\gamma(s',m'))$ is a nonzero morphism in $\mathcal{C}_{P_{\text{sp}}(Q^F)}$, then  $\Omegan(\rho)$ is a $\Bbbk$-linear maps $f=(f_{x})_{x \in \P}$  such that $f_{x}$ is a map $f_{x}:M_{x}^{\gamma(s,m)} \rightarrow M_{x}^{\gamma(s',m')}$ for any $x \in \P$. Let $\overline{f}=(\overline{f}_{x})_{x \in Q_{0}} $ be a morphism  from $\F(\gamma(s,m))$ to $\F(\gamma(s',m'))$ in $\ind Q$  such that $\overline{f}=f$. Note that, if $\gamma(s,m)$ is a $sp$-segment in  $\mathcal{C}_{P_{\text{sp}}(Q^F)}$, the representations $\F(\gamma(s,m))=(\F(\gamma(s,m))_{x}, \varphi^{\gamma(s,m)}_{\alpha})=M(s+1,m)$ in $\ind Q$ and $\Omegan(\gamma(s,m))=(M_{x}^{\gamma(s,m)},{_y}h_{x}^{\gamma(s,m)})$ in $\md_{sp} \Bbbk \P$ have the same $\Bbbk$-vector space  $\F(\gamma(s,m))_{x}=M_{x}^{\gamma(s,m)}$ for $x \in Q_{0}$, and if $\alpha: x \rightarrow y $ is an arrow in $Q$ it is satisfied that $x \prec y$ in $\P$ and $\varphi^{\gamma(s,m)}_{\alpha}={_y}h_{x}^{\gamma(s,m)}$. Moreover,  the diagram

$$ \xymatrix@=1.3cm{ \F(\gamma(s,m))_{x}\ar[r]^{\varphi^{\gamma(s,m)}_{\alpha}} \ar[d]_{\overline{f}_{x}}& \F(\gamma(s,m))_{y} \ar[d]^{\overline{f}_{y}} \\
\F(\gamma(s',m'))_{x} \ar[r]^{\varphi^{\gamma(s',m')}_{\alpha}}&   \F(\gamma(s',m'))_{y}} $$
commutes in $\md \Bbbk Q$ if the diagram  
$$ \xymatrix@=1.3cm{ M_x^{\gamma(s,m)} \ar[r]^{{_y}h_{x}^{\gamma(s,m)}} \ar[d]_{f_{x}}& M_y^{\gamma(s,m)} \ar[d]^{f_{y}} \\
M_x^{\gamma(s',m')} \ar[r]^{{_y}h_{x}^{\gamma(s',m')}}&   M_y^{\gamma(s',m')}} $$ commutes in $\md_{sp} \Bbbk \P$. By using the equivalence of categories $\F:\mathcal{C}_{P(Q)}\rightarrow \ind Q$ in Definition \ref{functorF} \cite{BGMS}, there is a morphism $\rho \in \Hom_{\mathcal{C}_{P(Q)}}(\gamma(s,m), \gamma(s',m'))$ such that $\F(\rho)=\overline{f}$. Given that, $\gamma(s,m)$ and $\gamma(s',m')$ are $sp$-segments in $\mathcal{C}_{P(Q)}$, $\rho$ is a morphism in the full subcategory $\mathcal{C}_{P_{\text{sp}}(Q^F)}$. By definition of $\Omegan$, $\Omegan(\rho)$ is isomorphic to $f$.\\

Finally, we prove that functor $\Omegan$ is dence. Let $N$ be a indecomposable in $\md_{sp}(\Bbbk \P)$, by Lemma 3.10 (b), Theorem 4.4  \cite{schifflerserna} and Theorem \ref{functorF} \cite{BGMS} the support of $N$ is a connected as a subset of the quiver $Q$, i.e., there is a line segment $\gamma(s,m)$ in $\mathcal{C}_{P(Q)}$ such that $\F(\gamma(s,m))$  is a socle-projective module $M(s+1,m)$ in $\ind Q$  with $\supp(M)=\supp(N)$. In order to prove that $\gamma(s,m)$ belongs to $\mathcal{C}_{P_{\text{sp}}(Q^F)}$,  Proposition \ref{*segment} allows us to establish that $\gamma(s,m)$ is a line $\star$-segment in $P(Q)$. Now, we suppose  that $\alpha: x \rightarrow y$ is an alien arrow in $Q$ defined in the support of the injective indecomposable $I(i)$ with $x,y \in \supp(M)$  for some sink $i$ in $Q$, then $s+1 \leq x < i < y \leq m$ (resp. $s+1 \leq y < i < x \leq m$) in \textbf{n+1},  if $x \leq y$ (resp. $ y \leq x$) therefore $\gamma(s.m) \notin [\alpha]$.
\end{proof}

\begin{figure}[ht]
\begin{adjustbox}{scale=0.4,center}

\begin{tikzcd}[ampersand replacement=\&, row sep= normal, column sep = normal]
\&\smrd{ \colorb{7}}\&\&\smrd{ \colorb{6}\\ \colorb{4}\\ \colorb{5}} \&\& \sm{ \colorb{3}\\ \colorb{1}\\ \colorb{2}}\\
\&\&\smrud{&  \colorb{6}\\ \colorb{4}&&  \colorb{7}\\  \colorb{5} }\&\&\smrud{& \colorb{3}&&  \colorb{6}\\ \colorb{1}&&  \colorb{4}\\ \colorb{2}&&  \colorb{5}} \&\\
\smr{ \colorb{5}}\&\smrud{ \colorb{4}\\ \colorb{5}} \& \&\smrud{&  \colorb{3}&&  \colorb{6}\\ \colorb{1}&&  \colorb{4}&&  \colorb{7}\\ \colorb{2}&&  \colorb{5}}\&\&\sm{& \colorb{3}&&  \colorb{6}\\ &&  \colorb{4}\\ &&  \colorb{5}}\\
\&\&\smrud{&  \colorb{3} \\ \colorb{1}&&  \colorb{4}\\ \colorb{2}&&  \colorb{5}} \&\& \smrud{ \colorb{3}&&  \colorb{6}\\ &  \colorb{4}&&  \colorb{7}\\  &  \colorb{5}}\&\\
\&\smru{ \colorb{1}\\ \colorb{2}}\&\&\smru{ \colorb{3}\\ \colorb{4}\\ \colorb{5}}\&\&\sm{ \colorb{6}\\ \colorb{7}}\\
\smru{ \colorb{2}}
\end{tikzcd}
\end{adjustbox}
\caption{Auslander-Reiter quiver of $\md_{sp}(\Bbbk \P)$}\label{AR}
\end{figure}
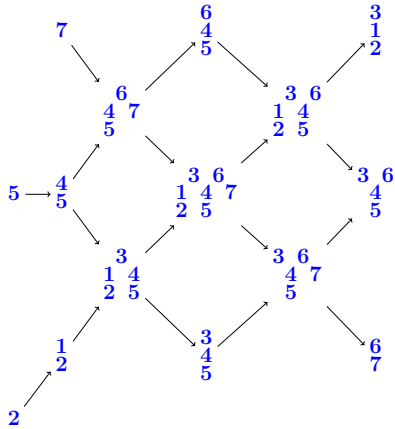

\begin{example}
Given the three-peak poset $\P=\P_{Q^{F}}$ shown in Examples \ref{non-positive-weight} and \ref{alien set}, Figure \ref{AR} illustrates the corresponding Auslander-Reiten quiver of the category of finitely generated socle-projective modules \(\md_{sp}(\Bbbk \mathcal{P})\). On the other hand, Figure \ref{ARQ} displays all the line segments in \(\mathcal{C}_{P(Q)}\). Specifically, the blue line segments correspond to the category of sp-segments \(P_{sp}(Q^{F})\); thus, the Auslander-Reiten quiver of this category is shown in Figure \ref{figure7}, which is structurally identical to the Auslander-Reiten quiver depicted in Figure \ref{AR}. 
\end{example}

\def\pentascale{1}
\def\myyellow{white}
\def\mygray{gray}
\def\myboundaryedge{blue}
\def\myedgecolor{red}
\def\myedgesize{thick}
\tikzset{->-/.style={decoration={
			markings,
			mark=at position #1 with {\arrow{stealth}}},postaction={decorate}}}

\newcommand\mypentagon{
	\begin{scope}[opacity=0.8,scale=0.4]
	
		\node (0) at (0,0) [fill,circle,inner sep=0.3] {};
		\node (1) at (0.8,1.4) [opacity=1,fill,circle,inner sep=0.3] {};
		\node (3) at (2.2,2) [fill,circle,inner sep=0.3] {};
		\node (4) at (3.8,2) [fill,circle,inner sep=0.3]{}; 
		\node (6) at (5.2,1.4) [fill,circle,inner sep=0.3] {};
		\node (7) at (6.4,0) [fill,circle,inner sep=0.3] {};
		\node (2) at (1.4,-1.2) [fill,circle,inner sep=0.3] {};
		\node (5) at (4.6,-1.2) [fill,circle,inner sep=0.3] {};
		\draw (0) -- (1) -- (3) -- (4) -- (6) -- (7)-- (5) -- (2)-- (0);
				\coordinate (northsource) at (7,2);
		\coordinate (northtarget) at (8.5,3.5);
		\coordinate (southsource) at (7,-2);
		\coordinate (southtarget) at (8.5,-3.5);
		
		\draw (0) node[left] {$0$};
		\draw (1) node[above] {$1$};
		\draw (2) node[below] {$2$};
		\draw (3) node[above] {$3$};
		\draw (4) node[above] {$4$};
		\draw (5) node[below] {$5$};
		\draw (6) node[above] {$6$};
		\draw (7) node[right] {$7$};
	\end{scope}
}
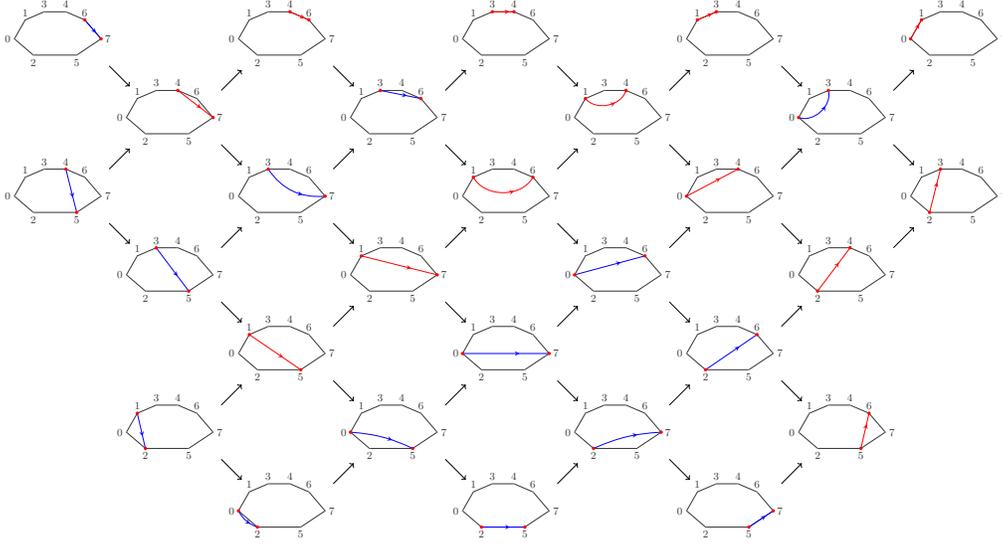
\begin{figure}[ht]
\begin{adjustbox}{scale=0.45, center}

\iftrue

\begin{tikzpicture}
\tikzstyle{every node}=[font=\small]
\matrix[outer sep=0pt,column sep=-5mm,row sep=-5mm,
ampersand replacement=\&]
{
	\mypentagon
	\draw[->-=0.7,blue,\myedgesize] (6) -- (7);
	\draw (6) node [red,fill,circle,inner sep=1pt] {};
	\draw (7) node [red,fill,circle,inner sep=1pt] {};
	\draw[->] (southsource) -- (southtarget); 
	\&
	\&
	\mypentagon{\mygray}
	\draw[->-=0.8,\myedgecolor,\myedgesize] (4) -- (6);
	\draw (4) node [red,fill,circle,inner sep=1pt] {};
	\draw (6) node [red,fill,circle,inner sep=1pt] {};
	\draw[->] (southsource) -- (southtarget);
	\&
	\&
	\mypentagon{\mygray}
	
	\draw[->-=0.8,\myedgecolor,\myedgesize] (3) -- (4);
	\draw (3) node [red,fill,circle,inner sep=1pt] {};
	\draw (4) node [red,fill,circle,inner sep=1pt] {};
	\draw[->] (southsource) -- (southtarget);
	\&
	\&
	\mypentagon{\mygray}
	\draw[->-=0.8,\myedgecolor,\myedgesize] (1) -- (3);
	\draw (1) node [red,fill,circle,inner sep=1pt] {};
	\draw (3) node [red,fill,circle,inner sep=1pt] {};
	\draw[->] (southsource) -- (southtarget);
	\&
	\&
	\mypentagon{\myyellow}
	\draw[->-=0.8,\myedgecolor,\myedgesize] (0) -- (1);
	\draw (0) node [red,fill,circle,inner sep=1pt] {};
	\draw (1) node [red,fill,circle,inner sep=1pt] {};
	\\
	\&
	\mypentagon{\mygray}
	\draw[->-=0.66,\myedgecolor,\myedgesize] (4) -- (7);
	\draw (4) node [red,fill,circle,inner sep=1pt] {};
	\draw (7) node [red,fill,circle,inner sep=1pt] {};
	\draw[->] (northsource) -- (northtarget);
	\draw[->] (southsource) -- (southtarget);
	\&
	\&
	\mypentagon{\mygray} 
	
	\draw[->-=0.66,blue,\myedgesize] (3) -- (6);
	\draw (3) node [red,fill,circle,inner sep=1pt] {};
	\draw (6) node [red,fill,circle,inner sep=1pt] {};
	\draw[->] (northsource) -- (northtarget);
	\draw[->] (southsource) -- (southtarget);
	\&
	\&
	\mypentagon{\mygray}
	
	\draw[->-=0.66,\myedgecolor,\myedgesize] (1) to [bend right=60] (4);
	\draw (1) node [red,fill,circle,inner sep=1pt] {};
	\draw (4) node [red,fill,circle,inner sep=1pt] {};
	\draw[->] (northsource) -- (northtarget);
	\draw[->] (southsource) -- (southtarget);
	\&
	\&
	\mypentagon{\mygray}
	
	\draw[->-=0.66,blue,\myedgesize] (0) to [bend right=60] (3);
	\draw (0) node [red,fill,circle,inner sep=1pt] {};
	\draw (3) node [red,fill,circle,inner sep=1pt] {};
	\draw[->] (northsource) -- (northtarget);
	\draw[->] (southsource) -- (southtarget);
	\&
	\&
	\&
	\&
	\\
	\mypentagon{\myyellow}
	\draw[->-=0.66,blue,\myedgesize] (4) to (5);
	\draw (4) node [red,fill,circle,inner sep=1pt] {};
	\draw (5) node [red,fill,circle,inner sep=1pt] {};
	\draw[->] (northsource) -- (northtarget);
	\draw[->] (southsource) -- (southtarget);
	\&
	\&
	\mypentagon{\mygray}
	\draw[->-=0.66,blue,\myedgesize] (3) to [bend right=30] (7);
	\draw (3) node [red,fill,circle,inner sep=1pt] {};
	\draw (7) node [red,fill,circle,inner sep=1pt] {};
	\draw[->] (northsource) -- (northtarget);
	\draw[->] (southsource) -- (southtarget);
	\&
	\&
	\mypentagon{\mygray}
	\draw[->-=0.66,\myedgecolor,\myedgesize] (1) to [bend right=60] (6);
	\draw (1) node [red,fill,circle,inner sep=1pt] {};
	\draw (6) node [red,fill,circle,inner sep=1pt] {};
	\draw[->] (northsource) -- (northtarget);
	\draw[->] (southsource) -- (southtarget);
	\&
	\&
	\mypentagon{\myyellow}
	\draw[->-=0.66,\myedgecolor,\myedgesize] (0) to (4);
	\draw (0) node [red,fill,circle,inner sep=1pt] {};
	\draw (4) node [red,fill,circle,inner sep=1pt] {};
	\draw[->] (northsource) -- (northtarget);
	\draw[->] (southsource) -- (southtarget);
	\&
	\&
	\mypentagon{\myyellow}
	\draw[->-=0.70,\myedgecolor,\myedgesize] (2) -- (3);
	\draw (2) node [red,fill,circle,inner sep=1pt] {};
	\draw (3) node [red,fill,circle,inner sep=1pt] {};
	\&
	\&
	
	\\
	\&
	\mypentagon{\myyellow}
	\draw[->-=0.66,blue,\myedgesize] (3) to (5);
	\draw (3) node [red,fill,circle,inner sep=1pt] {};
	\draw (5) node [red,fill,circle,inner sep=1pt] {};
	\draw[->] (northsource) -- (northtarget);
	\draw[->] (southsource) -- (southtarget);
	\&
	\&
	\mypentagon{\mygray}
	\draw[->-=0.66,\myedgecolor,\myedgesize] (1) to (7);
	\draw (1) node [red,fill,circle,inner sep=1pt] {};
	\draw (7) node [red,fill,circle,inner sep=1pt] {};
	\draw[->] (northsource) -- (northtarget);
	\draw[->] (southsource) -- (southtarget);
	\&
	\&
	\mypentagon{\myyellow}
	\draw[->-=0.66,blue,\myedgesize] (0) to (6);
	\draw (0) node [red,fill,circle,inner sep=1pt] {};
	\draw (6) node [red,fill,circle,inner sep=1pt] {};
	\draw[->] (northsource) -- (northtarget);
	\draw[->] (southsource) -- (southtarget);
	\&
	\&
	\mypentagon
	\draw[->-=0.66,\myedgecolor,\myedgesize] (2) to (4);
	\draw (2) node [red,fill,circle,inner sep=1pt] {};
	\draw (4) node [red,fill,circle,inner sep=1pt] {};
	\draw[->] (northsource) -- (northtarget);
	\&
	\&
	\&
	\&
	\\
	\&
	\&
	\mypentagon{\myyellow}
	\draw[->-=0.66,\myedgecolor,\myedgesize] (1) to (5);
	\draw (1) node [red,fill,circle,inner sep=1pt] {};
	\draw (5) node [red,fill,circle,inner sep=1pt] {};
	\draw[->] (northsource) -- (northtarget);
	\draw[->] (southsource) -- (southtarget);
	\&
	\&
	\mypentagon{\myyellow}
	\draw[->-=0.66,blue,\myedgesize] (0) to (7);
	\draw (0) node [red,fill,circle,inner sep=1pt] {};
	\draw (7) node [red,fill,circle,inner sep=1pt] {};
	\draw[->] (northsource) -- (northtarget);
	\draw[->] (southsource) -- (southtarget);
	\&
	\&
	\mypentagon{\myyellow}
	\draw[->-=0.66,blue,\myedgesize] (2) to (6);
	\draw (2) node [red,fill,circle,inner sep=1pt] {};
	\draw (6) node [red,fill,circle,inner sep=1pt] {};
	\draw[->] (northsource) -- (northtarget);
	\draw[->] (southsource) -- (southtarget);
	\&
	\&
	\&
	\&
	\&
	\&
	\\
	\&
	\mypentagon
	\draw[->-=0.66,blue,\myedgesize] (1) to (2);
	\draw (1) node [red,fill,circle,inner sep=1pt] {};
	\draw (2) node [red,fill,circle,inner sep=1pt] {};
	\draw[->] (northsource) -- (northtarget);
	\draw[->] (southsource) -- (southtarget);
	\&
	\&
	\mypentagon
	\draw[->-=0.66,blue,\myedgesize] (0) to [bend left=10] (5);
	\draw (0) node [red,fill,circle,inner sep=1pt] {};
	\draw (5) node [red,fill,circle,inner sep=1pt] {};
	\draw[->] (northsource) -- (northtarget);
	\draw[->] (southsource) -- (southtarget);
	\&
	\&
	\mypentagon
	\draw[->-=0.66,blue,\myedgesize] (2) to [bend left=10] (7);
	\draw (2) node [red,fill,circle,inner sep=1pt] {};
	\draw (7) node [red,fill,circle,inner sep=1pt] {};
	\draw[->] (northsource) -- (northtarget);
	\draw[->] (southsource) -- (southtarget);
	\&
	\&
	\mypentagon{\myyellow}
	\draw[->-=0.70,\myedgecolor,\myedgesize] (5) -- (6);
	\draw (5) node [red,fill,circle,inner sep=1pt] {};
	\draw (6) node [red,fill,circle,inner sep=1pt] {};
	
	\\
	\&
	\&
	\mypentagon
	\draw[->-=0.70,blue,\myedgesize] (0) to [bend right=20] (2);
	\draw (0) node [red,fill,circle,inner sep=1pt] {};
	\draw (2) node [red,fill,circle,inner sep=1pt] {};
	\draw[->] (northsource) -- (northtarget);
	\&
	\&
	\mypentagon
	\draw[->-=0.66,blue,\myedgesize] (2) to (5);
	\draw (2) node [red,fill,circle,inner sep=1pt] {};
	\draw (5) node [red,fill,circle,inner sep=1pt] {};
	\draw[->] (northsource) -- (northtarget);
	\&
	\&
	\mypentagon
	\draw[->-=0.70,blue,\myedgesize] (5) to (7);
	\draw (5) node [red,fill,circle,inner sep=1pt] {};
	\draw (7) node [red,fill,circle,inner sep=1pt] {};
	\draw[->] (northsource) -- (northtarget);
	\\
};
\end{tikzpicture}
\fi

\end{adjustbox}
\caption{Auslander-Reiten quiver of the category of segments $\mathcal{C}_{P(Q)}$. The indecomposable objets of the category of sp-segments $\mathcal{C}_{P(Q^F)}$ corresponds to the blue segments}
\label{ARQ}
\end{figure}
\begin{figure}
    \centering
\includegraphics[scale=0.18]{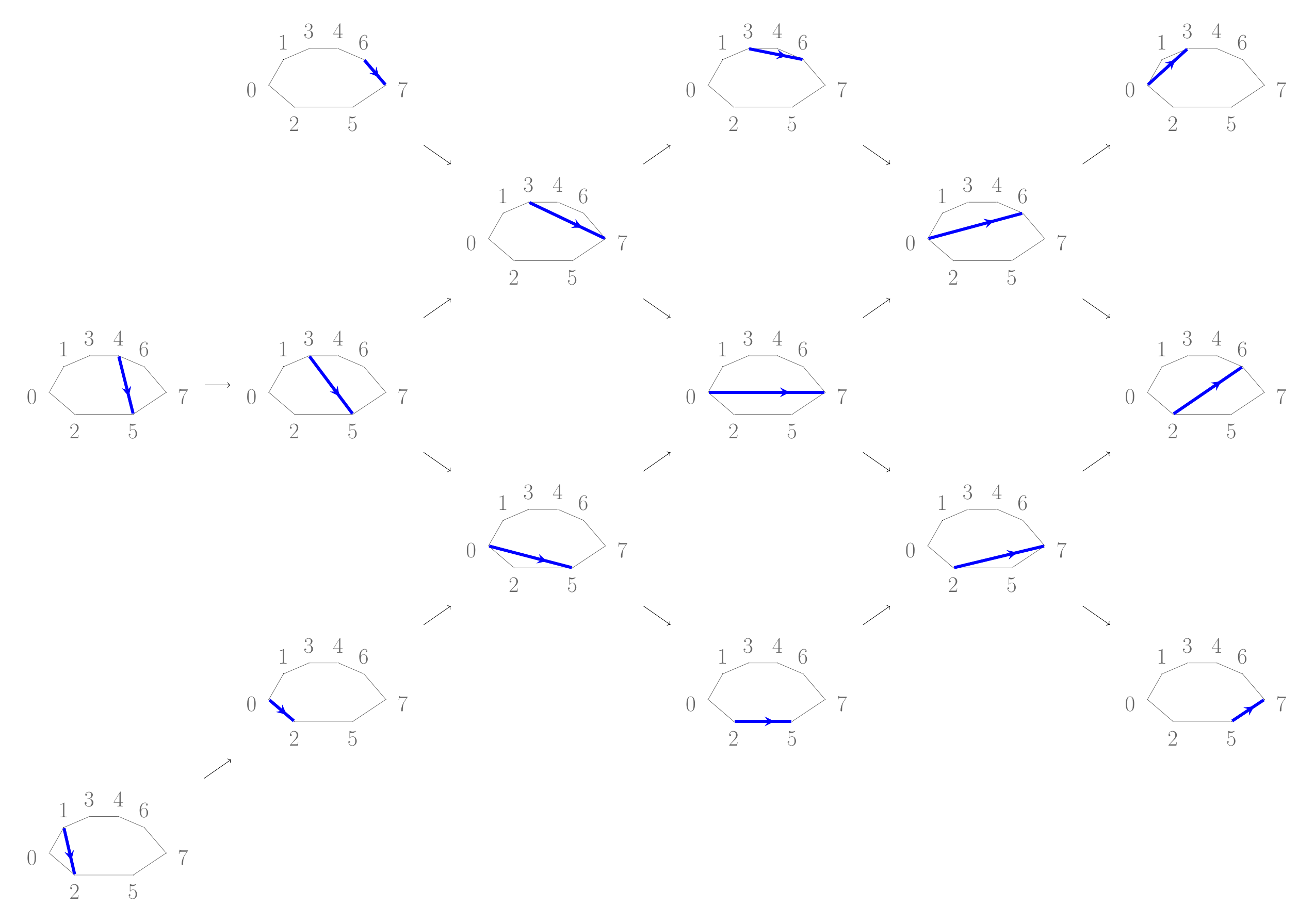}
    \caption{Auslander-Reiten quiver of the category of sp-segments $\mathcal{C}_{P(Q^F)}$} 
    \label{figure7}
\end{figure}
\begin{center}    
\end{center}
\newpage
\section{Stability arising from the geometric model of BGMS}  \label{section4}
In this section, we prove Theorem \ref{theototal}, which establishes a stability function for the category $\md_{sp}(\Bbbk\P)$, derived from the geometric model introduced in Section~\ref{section3}. In fact, Theorem \ref{theototal} is similar to Theorem 5.3 in \cite{BGMS}. Finally, using the geometric model from Section~\ref{section3} and the stability conditions associated with bilinear forms discussed in Section~\ref{section2}, we derive specific stability conditions for each poset of type $\mathbb{A}$.\\

Following the exposition in \cite{BGMS}, we introduce the notion of stability via a central charge associated with the category $\md_{sp}(\Bbbk\P)$. As discussed in Section \ref{subsection:semistable representations}, and since the categories $\md_{sp}(\Bbbk\P)$ and $\P$-spr are equivalent, one checks that  $K_0(\md_{sp}(\Bbbk\P))$ is isomorphic to  $\mathbb{Z}^{\P}$.\\

According to \cite[Definition 5.1]{BGMS} we have the following definition.
\begin{definition}
A \emph{central charge} is a group homomorphism $ Z\colon K_0(\md_{sp}(\Bbbk\P)) \to \mathbb{C} $ such that for all nonzero socle-projective representations $M$, the complex number $ Z[M] $ lies in the strict right half-plane\footnote{\cite{Bridgeland} uses the strict \emph{upper} half-plane}:
$$
\mathbb{H} = \left\{ r\,e^{i\pi\phi} \mid r > 0 \textup{ and } -\frac{1}{2} < \phi < \frac{1}{2} \right\}.
$$

Given a central charge $ Z $, we obtain an associated \emph{stability function}:
$$
\phi(M) = \frac{1}{\pi}\, \textup{arg}(Z[M]),
$$
where $ \textup{arg}(Z[M]) $ denotes the principal argument of the complex number $ Z[M]$. Given a stability function $ \phi$, a nonzero socle-projective  $ M $ is called \emph{$\phi$-stable} if every nonzero proper subrepresentation  $L \subsetneq M$ satisfies $ \phi(L) < \phi(M) $.
\end{definition}

\begin{remark}
 If a socle-projective representation $M $ is $\phi$-stable, then $ M $ is indecomposable. In fact, if $ M = N \oplus L $ in $\md_{sp}(\Bbbk\P)$, then $N$ and $ L $ are proper subrepresentations of $ M $. Thus, $\phi(N) < \phi(M)$ and $\phi(L) < \phi(M)$. Without loss of generality, suppose that
$$
\phi(M) = \frac{1}{\pi} \text{arg } (Z[M])=\frac{1}{\pi} \text{arg } (Z[N] + Z[L]) \leq \frac{1}{\pi} \text{arg } (Z[L]) = \phi(L).
$$
This leads to a contradiction.
\end{remark}

Given a poset $\P$ of type $\mathbb{A}$ such that $\P$ is the poset associated with the quiver $Q^F$ (as in Section \ref{section3}), we take the set $P_{sp}(Q^F) \subseteq \cale$ of all the line sp-segments in $P(Q)$ associated to the poset $\P$ via the functor $\Omega$ given in Theorem \ref{categoricalequivalence}  . According to \cite{BGMS}, since the vertices of $P(Q)$ are points in the plane, we can define a map
$$
\textup{vec} \colon P_{sp}(Q^F) \to \mathbb{C}, \quad \zg(s,m) \mapsto \vec{\zg}(s,m),
$$
where $\vec{\zg}(s,m)$ is the complex number $re^{i\theta}$ given by the vector with the same direction and magnitude as the oriented line segment from point $s$ to point $m$.

In the same way, given that
$\Omega \colon P_{sp}(Q^F) \to \md_{sp}(\Bbbk \P)$ is a bijection between objects, we define $\Phi_n$ as the inverse of $\Omega$, and we denote $\vec{\Phi_n} = \textup{vec} \circ \Phi_n$. Note that if $M = M(s,m)$ is a socle-projective representation in $\md_{sp}(\Bbbk \P)$,
$$
\vec{\Phi_n}(M) = \vec{\zg}(s-1,m) = r(M) \, e^{i\theta(M)},
$$
where $r(M)$ is the length of the vector $\vec{\Phi_n}(M) = \vec{\zg}(s-1,m)$, and $\theta(M)$ is the angle from the positive real axis to the vector $\vec{\Phi_n}(M)$ in the complex plane (see an example in Figure 7 of \cite{BGMS}).

Then $\vec{\Phi_n}$ induces a central charge given by

\begin{equation}
 \label{eq A1} 
 Z([M]) = \vec{\Phi_n}(M) = r(M) \, e^{i\theta(M)}.
\end{equation}
Note that $Z$ maps every socle-projective representation in $\md_{sp}(\Bbbk \P)$ to the strict right half-plane, and thus $Z$ is a central charge. The corresponding stability function $\phi$ is given by
$$\phi(M) = \frac{1}{\pi} \theta(M) \in \textstyle \left(-\frac{1}{2}, \frac{1}{2}\right).$$

\begin{theorem} \label{theototal}
Let $\P = \P_{Q^F}$ be a poset of type $\mathbb{A}$, where $F$ is an alien set for $Q$ as defined in Section~\ref{section3}. Then the function that associates to every indecomposable socle-projective representation $M = M^{\gamma(s,m)}$ the angle of the corresponding oriented line segment $\Phi_n(M) = \gamma(s,m)$ is a stability function, for which every indecomposable socle-projective representation of $\P$ is stable.
\end{theorem}

\begin{proof} 
Let $M$ be a representation in $\md_{sp}(\Bbbk\P)$. Following Lemma \ref{bijectionproper} and the equivalence of categories between $\md_{sp}(\Bbbk\P)$ and $\P$-spr (see \cite{simson3}), any proper subrepresentation of $M$ in $\md_{sp}(\Bbbk\P)$ has the form $M^I = (M^I_x, {}_yh_x)$, where $I$ is a proper subset of $\mathrm{Supp} M \cap \max\P$, $M^I_x = \Bbbk$ for each $x \in I_{\bigtriangledown} \setminus (I^c)_{\bigtriangledown}$ and zero otherwise (see Lemma \ref{Supp}), and if $x \preceq y$ in $\P$, the map ${}_yh_x$ is $1_{\Bbbk}$ whenever possible.

We see that $\phi(M^I) < \phi(M)$. In effect, we say that a set of maximal points $\{z_1, \dots, z_m\}$ in $\P$ is consecutive when $(z_{i-1})_{\vartriangle} \cap (z_{i})_{\vartriangle} \neq \emptyset$ for each $i = 2, \dots, m$. Then, the set $I$ is partitioned into consecutive disjoint subsets $I_1, \dots, I_s$. Then $M^I$ is the direct sum of $M^{I_1} \oplus \dots \oplus M^{I_s}$, where each $M^{I_j}$ is an indecomposable object in $\md_{sp}(\Bbbk\P_{Q^F})$. Note that there exist some $t, t' \in \mathbb{Z}_{\geq 0}$ such that $\Phi_n(M^{I_j})$ is a sp-segment which goes from $R^t(s)$ to $R^{t'}(m)$ in the polygon $P(Q)$, and $\phi(M^{I_j}) < \phi(M)$. Hence, the angle of the sum $\vec{\Phi_n}(M^{I_1}) + \dots + \vec{\Phi_n}(M^{I_s})$ is  $\theta(M^I) < \theta(M)$. 
\end{proof}

\subsection*{Associated total stability function} 
Let $M=(M_i,{_j}h_i)_{i,j\in\P}$ be the sincere socle-projective representation of a sincere peak poset $\P$ of type $\mathbb{A}_n$ and let $\{e_1, \dots, e_n\}$ be the standard basis of $\Bbbk^n$. Then the vector dimension $\dimv M = \sum_{i=1}^n e_iM_i$. Moreover, the weight $\thv$ defined in Proposition \ref{thetasincere} can be written as
$$
\thv = \kpv(M)\wpv - \wpv(M)\kpv,
$$
where $\wpv = \sum_{i=1}^n b_{\P}(\dimv M, e_i) e_i$ and $\kpv = \sum_{i=1}^n b_{\P}(e_i, \dimv M) e_i$. According to Table \ref{sincereposets}, these vectors $\wpv$ and $\kpv$ are as follows:

\begin{itemize}
    \item Case (i) $\P = \mathcal{S}_1^{(r)}$.  
    Observe that $\max\P = \{z_1, \dots, z_r\}$ and $\min\P = \{x_1, \dots, x_{r-1}\}$. Then
    \begin{equation}\label{equation2.6}
        w_{z_{i}} = b_{\P}(\dimv M, e_{z_i}) = 
    \begin{cases} 
        0 & \mbox{if } i=1 \text{ or } r, \\ 
        -1 & \mbox{otherwise,}  
    \end{cases}
    \end{equation}
    $w_{x_{j}} = b_{\P}(\dimv M, e_{x_j}) = -\kappa_{x_{j}} = -b_{\P}(e_{x_j},\dimv M) = 1$ for all $1 \leq j \leq r-1$, and $\kappa_{z_{i}} = b_{\P}(e_{z_i}, \dimv M) = 1$ for all $1 \leq i \leq r$.

    \item Case (ii) $\P = S_2^{(t)}$.  
    Observe that $\max\P = \{z_1, \dots, z_r\}$ and $\min\P = \{x_1, \dots, x_{r}\}$. Then
    \begin{equation}\label{equation2.7}
    w_{z_{i}} = \kappa_{x_{r+1-i}} = 
    \begin{cases} 
        0 & \mbox{if } i=1, \\ 
        -1 & \mbox{otherwise,}  
    \end{cases}
    \end{equation}
    $w_{x_{i}} = \kappa_{z_{i}} = 1$ with $1 \leq i \leq r$.

    \item Case (iii) $\P = \mathcal{S}_3^{(t)}$.  
    It is the case dual to case (i), where the maximal elements are minimal, and vice versa. Then the vectors $\wpv$ and $\kpv$ are a rearrangement of those obtained in case (i).
\end{itemize}

Fixing a natural number $m \geq 1$, consider the function $Z_m: K_0(\md_{sp}(\Bbbk\P)) \to \mathbb{C}$ given by
$$
Z_m[M] = (\wpv + \nv) \cdot \dimv M + ((\kpv + \nv) \cdot \dimv M)i,
$$
where $\nv = \sum_{i=1}^n m e_i$. Since the vectors $\wpv + \nv$ and $\kpv + \nv$ have positive entries, the complex number $Z_m(M)$ is in the upper half plane $\mathbb{H}$. Therefore, $Z_m$ induces a stability function
$$
\phi_m(M) = \frac{1}{\pi} \left( \frac{\pi}{2} - \text{arg}(Z_m[M]) \right).
$$

\begin{theorem} \label{Theorempeak}
    Let $\P$ be a sincere peak-poset of type $\mathbb{A}_n$, $n \geq 3$. For all $m \geq 1$, every indecomposable socle-projective representation of $\P$ is stable with stability function $\phi_m$.
\end{theorem}

\begin{proof} 
Since $\P$ is a sincere poset of type $\mathbb{A}_n$, according to Table \ref{sincereposets}, $\P$ is viewed as a Dynkin quiver $Q$ of type $\mathbb{A}_n$. We are going to construct specific coordinates for the $(n+1)$-gon $P'(Q)$ associated with $Q$ (see Section~\ref{section3.1}), and then Theorem \ref{theototal} will imply our result.

The vertices of $Q$ are labeled from left to right, from $1$ to $n$, as shown in the following picture:

\begin{center}
\begin{tabular}{|c|c|c|}
\hline
\begin{tikzpicture}
 \node (1) at (0,0) {$1$};
 
\node (A) at (1.55,0.5) {Case $\P=\mathcal{S}^{(r)}_{1}$:};
\node [right of=1, node distance=0.7cm](2) {$3$};

\node [right of=2, node distance=0.7cm](3) {$5$};
\node [right of=3, node distance=0.7cm] (4)  {};
\node [right of=4, node distance=0.7cm] (5) {$_{n}$};

\node [below  of=2, node distance=0.7cm] (2')  {$2$};
\node [below  of=3, node distance=0.7cm] (3')  {$4$};
\node [below  of=4, node distance=0.35cm] (4')  {$\dots$};
\node [below  of=5, node distance=0.7cm] (5')  {$_{n-1}$};

   
    \draw [blue, thick, ->, shorten <=-2pt, shorten >=-2pt] (2') -- (2);
    
\draw [blue, thick, ->, shorten <=-2pt, shorten >=-2pt] (3') -- (3);
    \draw [blue, thick, ->, shorten <=-2pt, shorten >=-2pt] (5') -- (5);    
    
    \draw [blue, thick, ->, shorten <=-2pt, shorten >=-2pt] (2') -- (1);
    
 \draw [blue, thick, ->, shorten <=-2pt, shorten >=-2pt] (3') -- (2);
    
\end{tikzpicture} & 

\begin{tikzpicture}
 \node (1) at (0,0) {$1$};
 
\node (A) at (1.55,0.5) {Case $\P=\mathcal{S}^{(r)}_{2}$:};
\node [right of=1, node distance=0.7cm](2) {$3$};

\node [right of=2, node distance=0.7cm](3) {$5$};
\node [right of=4, node distance=0.7cm] (5) {$_{n-1}$};

\node [below  of=2, node distance=0.7cm] (2')  {$2$};
\node [below  of=3, node distance=0.7cm] (3')  {$4$};
\node [below  of=4, node distance=0.35cm] (4')  {$\dots$};
\node [below  of=5, node distance=0.7cm] (5')  {$_{n-2}$};
\node [right of=5', node distance=0.7cm] (6')  {$_{n}$};

  
    \draw [blue, thick, ->, shorten <=-2pt, shorten >=-2pt] (2') -- (2);
    
\draw [blue, thick, ->, shorten <=-2pt, shorten >=-2pt] (3') -- (3);
    \draw [blue, thick, ->, shorten <=-2pt, shorten >=-2pt] (5') -- (5);    
    
    \draw [blue, thick, ->, shorten <=-2pt, shorten >=-2pt] (2') -- (1);
    
 \draw [blue, thick, ->, shorten <=-2pt, shorten >=-2pt] (3') -- (2);
    
\draw [blue, thick, ->, shorten <=-2pt, shorten >=-2pt] (6') -- (5);

\end{tikzpicture}

 & 
 
\begin{tikzpicture}
 \node (1) at (0,0) {$2$};
 
\node (A) at (1.55,0.5) {Case $\P=\mathcal{S}^{(r)}_{3}$:};
\node [right of=1, node distance=0.7cm](2) {$4$};

\node [right of=2, node distance=0.7cm](3) {$6$};
\node [right of=4, node distance=0.7cm] (5) {$_{n-1}$};

\node [below of=1, node distance=0.7cm] (1')  {$1$};
\node [below  of=2, node distance=0.7cm] (2')  {$3$};
\node [below  of=3, node distance=0.7cm] (3')  {$5$};
\node [below  of=4, node distance=0.7cm] (4')  {};
\node [below  of=4, node distance=0.35cm]   {$\cdots$};

\node [below  of=5, node distance=0.7cm] (5')  {$_{n-2}$};
\node [right of=5', node distance=0.7cm] (6')  {$_{n}$};

 \draw [blue, thick, ->, shorten <=-2pt, shorten >=-2pt] (1') -- (1);
    \draw [blue, thick, ->, shorten <=-2pt, shorten >=-2pt] (2') -- (2);
    
\draw [blue, thick, ->, shorten <=-2pt, shorten >=-2pt] (3') -- (3);
    \draw [blue, thick, ->, shorten <=-2pt, shorten >=-2pt] (5') -- (5);    
    
    \draw [blue, thick, ->, shorten <=-2pt, shorten >=-2pt] (2') -- (1);
    
 \draw [blue, thick, ->, shorten <=-2pt, shorten >=-2pt] (3') -- (2);
    
\draw [blue, thick, ->, shorten <=-2pt, shorten >=-2pt] (6') -- (5);
\end{tikzpicture} 
\\
\hline
 \end{tabular}

\end{center}

\subsection*{Positioning Vertices and Reflection of Polygons}

In the cases of $\mathcal{S}^{(r)}_{1}$ or $\mathcal{S}^{(r)}_{2}$, we position the vertices of $P'(Q)$ as follows: start by placing a vertex labeled $0$ at the point $(0,0)$; for each $1 \leq k \leq n$, place a vertex labeled $k$. Let $M$ be the sincere socle-projective representation of $\P$. We consider the weights  
$\wpv=\sum^{n}_{i=1}w_i e_i$ and $\kpv=\sum^{n}_{i=1}\kappa_i e_i$ from equations \eqref{equation2.6} and \eqref{equation2.7}. Assuming that  $w_{x_{0}} = \kappa_{x_{0}} = 0$, we have the following:
\begin{itemize}
\item[--] If the vertex $k$ in $Q$ corresponds to a maximal element $z_{i}$ of $\P$, the vertex $k$ of $P'(Q)$ is positioned at the coordinates:
  \[
  \left( (2j-1)m+\sum_{j=1}^{i} w_{z_{j}} + w_{x_{j-1}}, (2j-1)m + \sum_{j=1}^{i} \kappa_{z_{j}} + \kappa_{x_{j-1}} \right),
  \]
  where $j$ is the index such that $z_{j}$ is the $j$-th maximal element.
  
\item[--] If the vertex $k$ in $Q$ corresponds to a minimal element $x_{i}$ of $\P$, the vertex $k$ of $P(Q)$ is positioned at the coordinates:
  \[
  \left( 2jm +\sum_{j=1}^{i} w_{z_{j}} + w_{x_{j}}, 2jm+\sum_{j=1}^{i} \kappa_{z_{j}} + \kappa_{x_{j}} \right),
  \]
  where $j$ is the index such that $x_{j}$ is the $j$-th minimal element.
\end{itemize}
Note that when $\P = \mathcal{S}^{(r)}_{3}$, the polygon $P'(Q)$ is the reflection across the line $y=x$ of the polygon obtained in the case of $\mathcal{S}^{(r)}_{1}$. The reflection of the polygon $P'(Q)$ on the line $y=x$ gives us the geometric model considered in the categorical equivalence studied in Theorem \ref{categoricalequivalence}. 

Furthermore, observe that for all $x$ and $y$ in the positive real numbers,
\[
\text{arg}(x+yi) = \frac{\pi}{2} - \text{arg}(y+xi).
\]
Similar arguments in the proof of Theorem \ref{theototal} let us conclude that every indecomposable socle-projective representation of $\P$ is $\phi_m$-stable.
\end{proof}

\begin{example} \label{PolygonP}
Let $M$ be the sincere socle-projective representation  of  the poset $\P$ given by

\begin{center}
\begin{tikzpicture}
\node (1) at (0,0) {$\star$};
\node [left  of=1, node distance=0.15cm] (l2)  {$_{1}$};
\node [right  of=1, node distance=0.6cm] (2)  {$\star$};
\node [right  of=2, node distance=0.18cm] (l2)  {$_{3}$};
\node [right  of=1, node distance=1.2cm] (3)  {$\star$};
\node [right  of=2, node distance=0.78cm] (l2)  {$_{5}$};

\node [below  of=2, node distance=0.6cm] (1')  {$\circ$};
\node [left  of=1', node distance=0.18cm] (l2)  {$_{2}$};
\node [below  of=3, node distance=0.6cm] (2')  {$\circ$};
\node [right  of=2', node distance=0.18cm] (l2)  {$_{4}$};

\draw [shorten <=-2pt, shorten >=-2pt] (1) -- (1');
\draw [shorten <=-2pt, shorten >=-2pt] (2) -- (2');
\draw [shorten <=-2pt, shorten >=-2pt] (1') -- (2);
\draw [shorten <=-2pt, shorten >=-2pt] (2') -- (3);

\end{tikzpicture}
\end{center}

According to Equation \ref{equation2.6}, we have:
$$
\wpv = (0, -1, 0, 1, 1) \quad \text{and} \quad \kpv = (1, 1, 1, -1, -1).
$$
Figure \ref{Polygon1} shows the polygon $P'(Q)$ associated with $\P$ for $m = 1$, $m = 2$, and $m = 3$. Note that as $m$ increases, the polygon undergoes dilation.

\begin{figure}[ht]
\begin{center} 

\tikzset{every picture/.style={line width=0.75pt}} 

\begin{tikzpicture}[x=0.75pt,y=0.75pt,yscale=-1,xscale=1]

\draw    (136.67,230.17) -- (146.33,210.5) ;
\draw    (146.33,210.5) -- (166,189.83) ;
\draw    (166,189.83) -- (196,170.5) ;
\draw    (136.67,230.17) -- (166.33,209.83) ;
\draw    (166.33,209.83) -- (186,190.17) ;
\draw    (186,190.17) -- (196,170.5) ;
\draw    (136.67,230.17) -- (195.29,230.23) -- (213.5,230.25) ;
\draw [shift={(215.5,230.25)}, rotate = 180.06] [color={rgb, 255:red, 0; green, 0; blue, 0 }  ][line width=0.75]    (10.93,-3.29) .. controls (6.95,-1.4) and (3.31,-0.3) .. (0,0) .. controls (3.31,0.3) and (6.95,1.4) .. (10.93,3.29)   ;
\draw    (136.67,230.17) -- (136.5,152.25) ;
\draw [shift={(136.5,150.25)}, rotate = 89.88] [color={rgb, 255:red, 0; green, 0; blue, 0 }  ][line width=0.75]    (10.93,-3.29) .. controls (6.95,-1.4) and (3.31,-0.3) .. (0,0) .. controls (3.31,0.3) and (6.95,1.4) .. (10.93,3.29)   ;
\draw    (236,230.5) -- (236,101.25) ;
\draw [shift={(236,99.25)}, rotate = 90] [color={rgb, 255:red, 0; green, 0; blue, 0 }  ][line width=0.75]    (10.93,-3.29) .. controls (6.95,-1.4) and (3.31,-0.3) .. (0,0) .. controls (3.31,0.3) and (6.95,1.4) .. (10.93,3.29)   ;
\draw    (236,230.5) -- (363.5,229.76) ;
\draw [shift={(365.5,229.75)}, rotate = 179.67] [color={rgb, 255:red, 0; green, 0; blue, 0 }  ][line width=0.75]    (10.93,-3.29) .. controls (6.95,-1.4) and (3.31,-0.3) .. (0,0) .. controls (3.31,0.3) and (6.95,1.4) .. (10.93,3.29)   ;
\draw    (236,230.5) -- (256.5,201.25) ;
\draw    (256.5,201.25) -- (296,160.75) ;
\draw    (296,160.75) -- (345,119.75) ;
\draw    (236,230.5) -- (286.5,189.75) ;
\draw    (286.5,189.75) -- (326.5,149.75) ;
\draw    (326.5,149.75) -- (345,119.75) ;
\draw    (386.33,230.33) -- (386.33,51.67) ;
\draw [shift={(386.33,49.67)}, rotate = 90] [color={rgb, 255:red, 0; green, 0; blue, 0 }  ][line width=0.75]    (10.93,-3.29) .. controls (6.95,-1.4) and (3.31,-0.3) .. (0,0) .. controls (3.31,0.3) and (6.95,1.4) .. (10.93,3.29)   ;
\draw    (386.33,230.33) -- (564.33,230.33) ;
\draw [shift={(566.33,230.33)}, rotate = 180] [color={rgb, 255:red, 0; green, 0; blue, 0 }  ][line width=0.75]    (10.93,-3.29) .. controls (6.95,-1.4) and (3.31,-0.3) .. (0,0) .. controls (3.31,0.3) and (6.95,1.4) .. (10.93,3.29)   ;
\draw    (386.33,230.33) -- (416,190.75) ;
\draw    (416,190.75) -- (476.33,131) ;
\draw    (476.33,131) -- (547,71.25) ;
\draw    (386.33,230.33) -- (456,171.25) ;
\draw    (456,171.25) -- (516,110.25) ;
\draw    (516,110.25) -- (547,71.25) ;
\draw    (416,227) -- (416,233.25) ;
\draw    (457.5,226.5) -- (457.5,232.75) ;
\draw    (515.5,227.35) -- (515.5,233.6) ;
\draw    (546,226.25) -- (546,232.5) ;
\draw    (256.5,226.75) -- (256.5,233) ;
\draw    (287,226.75) -- (287,233) ;
\draw    (325.5,227.25) -- (325.5,233.5) ;
\draw    (346.5,226.75) -- (346.5,233) ;
\draw    (146.14,226.96) -- (146.14,233.21) ;
\draw    (165.93,227.18) -- (165.93,233.43) ;
\draw    (186.36,226.96) -- (186.36,233.21) ;
\draw    (476.5,227.25) -- (476.5,233.5) ;
\draw    (196.21,227.18) -- (196.21,233.43) ;
\draw    (296.21,226.89) -- (296.21,233.14) ;
\draw    (139.29,190.43) -- (133.29,190.43) ;
\draw    (139.29,209.86) -- (133.29,209.86) ;
\draw    (139.29,170.14) -- (133.29,170.14) ;
\draw    (238.71,200.43) -- (232.71,200.43) ;
\draw    (239,190.71) -- (233,190.71) ;
\draw    (238.71,160.43) -- (232.71,160.43) ;
\draw    (239,149.86) -- (233,149.86) ;
\draw    (239,120.49) -- (233,120.49) ;
\draw    (389.14,190.43) -- (383.14,190.43) ;
\draw    (388.86,170.14) -- (382.86,170.14) ;
\draw    (389.14,130.43) -- (383.14,130.43) ;
\draw    (389.14,110.49) -- (383.14,110.49) ;
\draw    (389.71,69.92) -- (383.71,69.92) ;

\draw (142.5,237.4) node [anchor=north west][inner sep=0.75pt]  [font=\fontsize{0.29em}{0.35em}\selectfont]  {$1$};
\draw (162.5,237.4) node [anchor=north west][inner sep=0.75pt]  [font=\fontsize{0.29em}{0.35em}\selectfont]  {$3$};
\draw (182.5,237.4) node [anchor=north west][inner sep=0.75pt]  [font=\fontsize{0.29em}{0.35em}\selectfont]  {$5$};
\draw (192.5,237.4) node [anchor=north west][inner sep=0.75pt]  [font=\fontsize{0.29em}{0.35em}\selectfont]  {$6$};
\draw (252.5,237.4) node [anchor=north west][inner sep=0.75pt]  [font=\fontsize{0.29em}{0.35em}\selectfont]  {$2$};
\draw (282.5,237.4) node [anchor=north west][inner sep=0.75pt]  [font=\fontsize{0.29em}{0.35em}\selectfont]  {$5$};
\draw (292.5,237.4) node [anchor=north west][inner sep=0.75pt]  [font=\fontsize{0.29em}{0.35em}\selectfont]  {$6$};
\draw (321.5,237.4) node [anchor=north west][inner sep=0.75pt]  [font=\fontsize{0.29em}{0.35em}\selectfont]  {$9$};
\draw (340.5,237.4) node [anchor=north west][inner sep=0.75pt]  [font=\fontsize{0.29em}{0.35em}\selectfont]  {$11$};
\draw (412,237.4) node [anchor=north west][inner sep=0.75pt]  [font=\fontsize{0.29em}{0.35em}\selectfont]  {$3$};
\draw (453.5,237.4) node [anchor=north west][inner sep=0.75pt]  [font=\fontsize{0.29em}{0.35em}\selectfont]  {$7$};
\draw (472.5,237.4) node [anchor=north west][inner sep=0.75pt]  [font=\fontsize{0.29em}{0.35em}\selectfont]  {$9$};
\draw (510,237.4) node [anchor=north west][inner sep=0.75pt]  [font=\fontsize{0.29em}{0.35em}\selectfont]  {$13$};
\draw (540,237.4) node [anchor=north west][inner sep=0.75pt]  [font=\fontsize{0.29em}{0.35em}\selectfont]  {$16$};
\draw (125,206) node [anchor=north west][inner sep=0.75pt]  [font=\fontsize{0.29em}{0.35em}\selectfont]  {$2$};
\draw (125,186) node [anchor=north west][inner sep=0.75pt]  [font=\fontsize{0.29em}{0.35em}\selectfont]  {$4$};
\draw (125,166) node [anchor=north west][inner sep=0.75pt]  [font=\fontsize{0.29em}{0.35em}\selectfont]  {$6$};
\draw (225,196) node [anchor=north west][inner sep=0.75pt]  [font=\fontsize{0.29em}{0.35em}\selectfont]  {$3$};
\draw (225,186) node [anchor=north west][inner sep=0.75pt]  [font=\fontsize{0.29em}{0.35em}\selectfont]  {$4$};
\draw (225,156) node [anchor=north west][inner sep=0.75pt]  [font=\fontsize{0.29em}{0.35em}\selectfont]  {$7$};
\draw (225,146) node [anchor=north west][inner sep=0.75pt]  [font=\fontsize{0.29em}{0.35em}\selectfont]  {$8$};
\draw (222,116) node [anchor=north west][inner sep=0.75pt]  [font=\fontsize{0.29em}{0.35em}\selectfont]  {$11$};
\draw (375.5,186) node [anchor=north west][inner sep=0.75pt]  [font=\fontsize{0.29em}{0.35em}\selectfont]  {$4$};
\draw (375.5,166) node [anchor=north west][inner sep=0.75pt]  [font=\fontsize{0.29em}{0.35em}\selectfont]  {$6$};
\draw (371,126) node [anchor=north west][inner sep=0.75pt]  [font=\fontsize{0.29em}{0.35em}\selectfont]  {$10$};
\draw (371,106) node [anchor=north west][inner sep=0.75pt]  [font=\fontsize{0.29em}{0.35em}\selectfont]  {$12$};
\draw (371,66) node [anchor=north west][inner sep=0.75pt]  [font=\fontsize{0.29em}{0.35em}\selectfont]  {$16$};
\draw (142,200) node [anchor=north west][inner sep=0.75pt]  [font=\fontsize{0.35em}{0.42em}\selectfont,color={rgb, 255:red, 208; green, 2; blue, 27 }  ,opacity=1 ]  {$1$};
\draw (252.5,190) node [anchor=north west][inner sep=0.75pt]  [font=\fontsize{0.35em}{0.42em}\selectfont,color={rgb, 255:red, 208; green, 2; blue, 27 }  ,opacity=1 ]  {$1$};
\draw (412,180) node [anchor=north west][inner sep=0.75pt]  [font=\fontsize{0.35em}{0.42em}\selectfont,color={rgb, 255:red, 208; green, 2; blue, 27 }  ,opacity=1 ]  {$1$};
\draw (129,230) node [anchor=north west][inner sep=0.75pt]  [font=\fontsize{0.35em}{0.42em}\selectfont,color={rgb, 255:red, 208; green, 2; blue, 27 }  ,opacity=1 ]  {$0$};
\draw (229,230) node [anchor=north west][inner sep=0.75pt]  [font=\fontsize{0.35em}{0.42em}\selectfont,color={rgb, 255:red, 208; green, 2; blue, 27 }  ,opacity=1 ]  {$0$};
\draw (379,230) node [anchor=north west][inner sep=0.75pt]  [font=\fontsize{0.35em}{0.42em}\selectfont,color={rgb, 255:red, 208; green, 2; blue, 27 }  ,opacity=1 ]  {$0$};
\draw (162.5,213.5) node [anchor=north west][inner sep=0.75pt]  [font=\fontsize{0.35em}{0.42em}\selectfont,color={rgb, 255:red, 208; green, 2; blue, 27 }  ,opacity=1 ]  {$2$};
\draw (282.5,193.5) node [anchor=north west][inner sep=0.75pt]  [font=\fontsize{0.35em}{0.42em}\selectfont,color={rgb, 255:red, 208; green, 2; blue, 27 }  ,opacity=1 ]  {$2$};
\draw (162,179.5) node [anchor=north west][inner sep=0.75pt]  [font=\fontsize{0.35em}{0.42em}\selectfont,color={rgb, 255:red, 208; green, 2; blue, 27 }  ,opacity=1 ]  {$3$};
\draw (182.5,193) node [anchor=north west][inner sep=0.75pt]  [font=\fontsize{0.35em}{0.42em}\selectfont,color={rgb, 255:red, 208; green, 2; blue, 27 }  ,opacity=1 ]  {$4$};
\draw (197,162) node [anchor=north west][inner sep=0.75pt]  [font=\fontsize{0.35em}{0.42em}\selectfont,color={rgb, 255:red, 208; green, 2; blue, 27 }  ,opacity=1 ]  {$5$};
\draw (292,150) node [anchor=north west][inner sep=0.75pt]  [font=\fontsize{0.35em}{0.42em}\selectfont,color={rgb, 255:red, 208; green, 2; blue, 27 }  ,opacity=1 ]  {$3$};
\draw (322,153) node [anchor=north west][inner sep=0.75pt]  [font=\fontsize{0.35em}{0.42em}\selectfont,color={rgb, 255:red, 208; green, 2; blue, 27 }  ,opacity=1 ]  {$4$};
\draw (347,112) node [anchor=north west][inner sep=0.75pt]  [font=\fontsize{0.35em}{0.42em}\selectfont,color={rgb, 255:red, 208; green, 2; blue, 27 }  ,opacity=1 ]  {$5$};
\draw (452.5,174.5) node [anchor=north west][inner sep=0.75pt]  [font=\fontsize{0.35em}{0.42em}\selectfont,color={rgb, 255:red, 208; green, 2; blue, 27 }  ,opacity=1 ]  {$2$};
\draw (472,121) node [anchor=north west][inner sep=0.75pt]  [font=\fontsize{0.35em}{0.42em}\selectfont,color={rgb, 255:red, 208; green, 2; blue, 27 }  ,opacity=1 ]  {$3$};
\draw (512.5,115) node [anchor=north west][inner sep=0.75pt]  [font=\fontsize{0.35em}{0.42em}\selectfont,color={rgb, 255:red, 208; green, 2; blue, 27 }  ,opacity=1 ]  {$4$};
\draw (548.5,63) node [anchor=north west][inner sep=0.75pt]  [font=\fontsize{0.35em}{0.42em}\selectfont,color={rgb, 255:red, 208; green, 2; blue, 27 }  ,opacity=1 ]  {$5$};

\filldraw[color=black] (146,210) circle (1pt);
\filldraw[color=black] (166,210) circle (1pt);
\filldraw[color=black] (166,190) circle (1pt);
\filldraw[color=black] (186,190) circle (1pt);
\filldraw[color=black] (196,170) circle (1pt);
\filldraw[color=black] (257,200) circle (1pt);
\filldraw[color=black] (286,190) circle (1pt);
\filldraw[color=black] (296,160) circle (1pt);
\filldraw[color=black] (326,150) circle (1pt);
\filldraw[color=black] (345,120) circle (1pt);
\filldraw[color=black] (416,190) circle (1pt);
\filldraw[color=black] (456,171) circle (1pt);
\filldraw[color=black] (476,131) circle (1pt);
\filldraw[color=black] (516,111) circle (1pt);
\filldraw[color=black] (546,72) circle (1pt);

\end{tikzpicture}

\end{center}
\caption{Polygons $P'(Q)$ associated with a poset $\P$ (left $m=1$, center $m=2$, right $m=3$)} \label{Polygon1}

\end{figure}
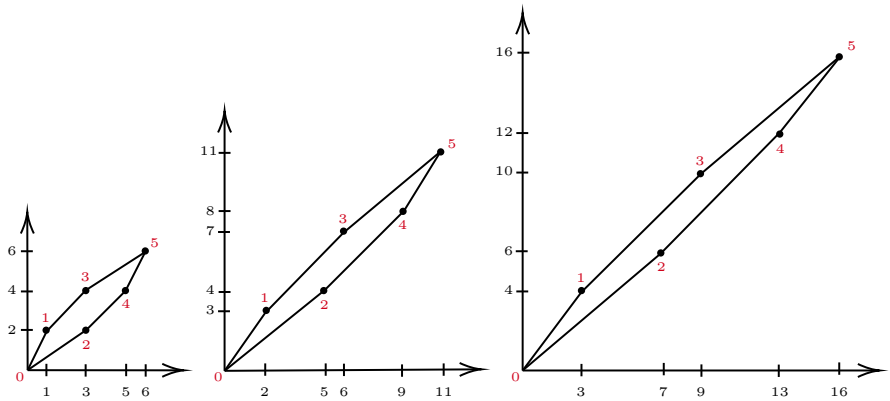

For $m = 2$, let $M'$ be the socle-projective representation  with $\supp M' = \{2, 3, 4, 5\}$. Using Lemmas \ref{Supp} and \ref{bijectionproper}, we find that the proper socle-projective subrepresentations of $M'$ are $N$ and $L$, where $\supp N = \{1\}$, and $\supp L = \{3,4\}$. In this case, the left part of Figure \ref{complex} illustrates the sp-segments corresponding to the complex numbers:
\begin{align*}
 Z_{2}[M']&= 5 + 4i \qquad \mbox{(blue)},\\
 Z_{2}[N]&= 1 + 2i \qquad \mbox{(purple)},\\
Z_{2}[L]&= 2 + 2i \qquad \mbox{(green)}.
\end{align*}
The sp-segments are represented in the right part of Figure \ref{complex}. There, we have used dashed lines for the respective vectors obtained by reflection over the line $ y = x $. According to Theorem \ref{Theorempeak}, it follows that $\phi_{2}(N)<\phi_{2}(M')$ and $\phi_{2}(L) < \phi_{2}(M')$.

\begin{figure} [ht]
    \begin{center}
        \tikzset{every picture/.style={line width=0.75pt}} 

\begin{tikzpicture}[x=0.75pt,y=0.75pt,yscale=-1,xscale=1]

\draw  (101.25,300.26) -- (119.8,261) ;
\draw    (159.8,220.2) -- (219.8,179.8) ;
\draw    (101.25,300.26) -- (159.8,260.6) ;
\draw    (159.8,260.6) -- (200.2,220.6) ;
\draw    (101.25,300.26) -- (217,300.42) -- (248.6,300.21) ;
\draw [shift={(250.6,300.2)}, rotate = 179.62] [color={rgb, 255:red, 0; green, 0; blue, 0 }  ][line width=0.75]    (10.93,-3.29) .. controls (6.95,-1.4) and (3.31,-0.3) .. (0,0) .. controls (3.31,0.3) and (6.95,1.4) .. (10.93,3.29)   ;
\draw    (101.25,300.26) -- (100.61,151.8) ;
\draw [shift={(100.6,149.8)}, rotate = 89.75] [color={rgb, 255:red, 0; green, 0; blue, 0 }  ][line width=0.75]    (10.93,-3.29) .. controls (6.95,-1.4) and (3.31,-0.3) .. (0,0) .. controls (3.31,0.3) and (6.95,1.4) .. (10.93,3.29)   ;
\draw    (119.96,293.36) -- (119.8,305.8) ;
\draw    (159.03,293.65) -- (159,305) ;
\draw    (200.57,293.23) -- (200.2,305.4) ;
\draw    (219.63,293.38) -- (219.8,305.8) ;
\draw    (107.08,220.78) -- (95.24,220.78) ;
\draw    (106.42,260.86) -- (94.57,260.86) ;
\draw    (106.95,180.77) -- (95.1,180.77) ;
\draw    (300.85,301.07) -- (300.6,172.2) ;
\draw [shift={(300.6,170.2)}, rotate = 89.89] [color={rgb, 255:red, 0; green, 0; blue, 0 }  ][line width=0.75]    (10.93,-3.29) .. controls (6.95,-1.4) and (3.31,-0.3) .. (0,0) .. controls (3.31,0.3) and (6.95,1.4) .. (10.93,3.29)   ;
\draw    (300.85,301.07) -- (427.4,299.82) ;
\draw [shift={(429.4,299.8)}, rotate = 179.43] [color={rgb, 255:red, 0; green, 0; blue, 0 }  ][line width=0.75]    (10.93,-3.29) .. controls (6.95,-1.4) and (3.31,-0.3) .. (0,0) .. controls (3.31,0.3) and (6.95,1.4) .. (10.93,3.29)   ;
\draw    (320.65,293.9) -- (321,306.6) ;
\draw    (340.4,294.37) -- (340.6,306.6) ;
\draw    (399.87,293.5) -- (400.2,305) ;
\draw    (306.15,259.86) -- (294.3,259.86) ;
\draw    (306.31,219.98) -- (294.46,219.98) ;
\draw [->-=0.66,color={rgb, 255:red, 74; green, 144; blue, 226 }  ,draw opacity=1 ]   (119.8,261) -- (219.8,179.8) ;
\draw [->-=0.66,color={rgb, 255:red, 189; green, 16; blue, 224 }  ,draw opacity=1 ]   (200.2,220.6) -- (219.8,179.8) ;
\draw [->-=0.66,color={rgb, 255:red, 126; green, 211; blue, 33 }  ,draw opacity=1 ]   (119.8,261) -- (159.8,220.2) ;
\draw [->,color={rgb, 255:red, 74; green, 144; blue, 226 }  ,draw opacity=1 ]   (300.85,301.07) -- (400.6,219.8) ;
\draw [->,color={rgb, 255:red, 126; green, 211; blue, 33 }  ,draw opacity=1 ]   (300.85,301.07) -- (340.2,261) ;
\draw [->,color={rgb, 255:red, 189; green, 16; blue, 224 }  ,draw opacity=1 ]   (300.85,301.07) -- (319.8,260.6) ;
\draw  [dotted]  (300.85,301.07) -- (429.4,169.8) ;
\draw [->,dashed, color={rgb, 255:red, 74; green, 144; blue, 226 }  ,draw opacity=1 ]   (300.85,301.07) -- (379.8,200.2) ;
\draw [->,dashed, color={rgb, 255:red, 189; green, 16; blue, 224 }  ,draw opacity=1 ]   (300.85,301.07) -- (340.6,279.8) ;
\draw    (381.2,294.57) -- (381.4,306.8) ;
\draw    (306.71,279.98) -- (294.86,279.98) ;
\draw    (305.91,200.38) -- (294.06,200.38) ;

\draw (116,312) node [anchor=north west][inner sep=0.75pt]  [font=\fontsize{0.59em}{0.71em}\selectfont]  {$1$};
\draw (116,243) node [anchor=north west][inner sep=0.75pt]  [font=\fontsize{0.59em}{0.71em}\selectfont,color={rgb, 255:red, 208; green, 2; blue, 27 }  ,opacity=1 ]  {$1$};
\draw (88.47,299.73) node [anchor=north west][inner sep=0.75pt]  [font=\fontsize{0.59em}{0.71em}\selectfont,color={rgb, 255:red, 208; green, 2; blue, 27 }  ,opacity=1 ]  {$0$};
\draw (158,264.93) node [anchor=north west][inner sep=0.75pt]  [font=\fontsize{0.59em}{0.71em}\selectfont,color={rgb, 255:red, 208; green, 2; blue, 27 }  ,opacity=1 ]  {$2$};
\draw (151.94,207.09) node [anchor=north west][inner sep=0.75pt]  [font=\fontsize{0.59em}{0.71em}\selectfont,color={rgb, 255:red, 208; green, 2; blue, 27 }  ,opacity=1 ]  {$3$};
\draw (201.26,228.71) node [anchor=north west][inner sep=0.75pt]  [font=\fontsize{0.59em}{0.71em}\selectfont,color={rgb, 255:red, 208; green, 2; blue, 27 }  ,opacity=1 ]  {$4$};
\draw (219.8,168.42) node [anchor=north west][inner sep=0.75pt]  [font=\fontsize{0.59em}{0.71em}\selectfont,color={rgb, 255:red, 208; green, 2; blue, 27 }  ,opacity=1 ]  {$5$};
\draw (155.17,312) node [anchor=north west][inner sep=0.75pt]  [font=\fontsize{0.59em}{0.71em}\selectfont]  {$3$};
\draw (196.37,312) node [anchor=north west][inner sep=0.75pt]  [font=\fontsize{0.59em}{0.71em}\selectfont]  {$5$};
\draw (215.71,312) node [anchor=north west][inner sep=0.75pt]  [font=\fontsize{0.59em}{0.71em}\selectfont]  {$6$};
\draw (81,255.8) node [anchor=north west][inner sep=0.75pt]  [font=\fontsize{0.59em}{0.71em}\selectfont]  {$2$};
\draw (81,216.48) node [anchor=north west][inner sep=0.75pt]  [font=\fontsize{0.59em}{0.71em}\selectfont]  {$4$};
\draw (81,176.88) node [anchor=north west][inner sep=0.75pt]  [font=\fontsize{0.59em}{0.71em}\selectfont]  {$6$};
\draw (317.17,312) node [anchor=north west][inner sep=0.75pt]  [font=\fontsize{0.59em}{0.71em}\selectfont]  {$1$};
\draw (336.77,312) node [anchor=north west][inner sep=0.75pt]  [font=\fontsize{0.59em}{0.71em}\selectfont]  {$2$};
\draw (377.17,312) node [anchor=north west][inner sep=0.75pt]  [font=\fontsize{0.59em}{0.71em}\selectfont]  {$4$};
\draw (396.37,312) node [anchor=north west][inner sep=0.75pt]  [font=\fontsize{0.59em}{0.71em}\selectfont]  {$5$};
\draw (282,276.08) node [anchor=north west][inner sep=0.75pt]  [font=\fontsize{0.59em}{0.71em}\selectfont]  {$1$};
\draw (282,256.48) node [anchor=north west][inner sep=0.75pt]  [font=\fontsize{0.59em}{0.71em}\selectfont]  {$2$};
\draw (282,216.88) node [anchor=north west][inner sep=0.75pt]  [font=\fontsize{0.59em}{0.71em}\selectfont]  {$4$};
\draw (282,196.88) node [anchor=north west][inner sep=0.75pt]  [font=\fontsize{0.59em}{0.71em}\selectfont]  {$5$};

\end{tikzpicture}

    \end{center}
    \caption{Some line segments on the polygon $P'(Q)$, along with their associated vectors and their reflections across the line 
$y=x$}
    \label{complex}
\end{figure}
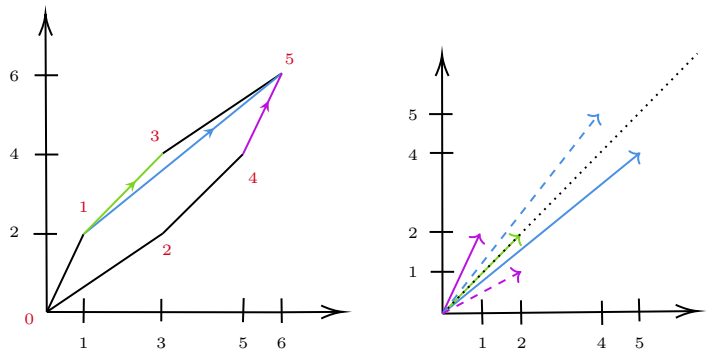




\end{example}

\bibliographystyle{plain}
\bibliography{references}

\end{document}